\documentclass[12pt,twoside]{amsart}
\usepackage[latin1]{inputenc}
\usepackage{times}
\usepackage{amsthm}
\usepackage{amssymb, verbatim}
\usepackage{amsmath}
\usepackage{eucal}
\usepackage{mathrsfs}

\topmargin = - 40 pt %distanza dal margine superiore 0= un bel po` negativo diminuisce
\usepackage[all]{xy}

\theoremstyle{definition}

\numberwithin{equation}{section}

\begin{document}

\title[$C^*$--quantum groupoids]
{  Quasi-coassociative $C^*$-quantum groupoids of   type $A$ and   modular $C^*$-categories}  
\author[S.~Ciamprone]{Sergio Ciamprone }
 \author[C.~Pinzari]{Claudia Pinzari}
 \email{ciamprone@mat.uniroma1.it, pinzari@mat.uniroma1.it}
\address{Dipartimento di Matematica, Sapienza Universit\`a   di
Roma \\ P.le Aldo Moro, 5 -- 00185 Rome, Italy}
\maketitle

\centerline{\it In memory of Joe Jerome P\' erez}

\medskip

\begin{abstract} We construct a new class of finite-dimensional $C^*$-quantum groupoids at   roots of unity $q=e^{i\pi/\ell}$,   with  limit the discrete dual of the classical ${\rm SU}(N)$ for large orders. The representation category of our groupoid   turns out to be tensor equivalent to   the well known quotient $C^*$-category of the category of tilting modules of the non-semisimple   quantum group $U_q({\mathfrak sl}_N)$ of Drinfeld, Jimbo and Lusztig.

As an algebra, the $C^*$-groupoid is a quotient of $U_q({\mathfrak sl}_N)$. As a coalgebra, it naturally reflects the categorical quotient construction. In particular, it is not coassociative, but   satisfies   axioms of the weak quasi-Hopf
$C^*$-algebras:
  quasi-coassociativity and non-unitality of the coproduct. There are  also a multiplicative counit, an   antipode, and  an $R$-matrix.

 For this,    we  give a general construction of   quantum groupoids for   complex simple Lie algebras ${\mathfrak g}\neq E_8$ and certain   roots of unity.
 Our main  tools  here are Drinfeld's coboundary associated to the $R$-matrix, related to the algebra involution, and certain canonical    projections introduced by Wenzl, which yield the coproduct and Drinfeld's associator in an explicit way. Tensorial properties of the negligible modules   reflect in a rather special nature of the associator. 
We next reduce the proof  of  the   categorical equivalence
  to the problems of establishing   semisimplicity and computing  dimension of the groupoid.
    In the case ${\mathfrak g}={\mathfrak sl}_N$ we construct a  (non-positive) Haar-type   functional on an associative version of the dual groupoid satisfying key non-degeneracy properties. This   enables us to complete the proof.

\end{abstract}

%%\maketitle

\section{Introduction}

Let    ${\mathfrak g}$ be a complex simple Lie algebra, $q$  a primitive root of unity,   and let $U_q({\mathfrak g})$ be Lusztig's  restricted quantum group \cite{CP}. It is well known that this is a non-semisimple algebra,  and so is its representation category. However, a semisimple tensor category, ${\mathcal F}$,  called the fusion category, can be obtained  as  a quotient of the category of tilting modules by the ideal of the negligible  ones  \cite{Lusztig, Lusztig_book, Andersen,  AP, GK, Turaev}.   
These categories play a prominent role   in  conformal field theory   \cite{BK}, subfactor theory \cite{Jones}, topological quantum field theory  and the study of invariants of $3$-manifolds \cite{Witten, RT}.

   For the roots of unity $q=e^{i\pi/d\ell}$ with $d$ the ratio of the square lengths of a long root to a  short root,
    Kirillov  introduced   a  $^*$--involution and an associated inner product  on the arrow spaces of  ${\mathcal F}$
 and  conjectured positivity \cite{Kirillov}.   Wenzl   proved   the conjecture, and derived  existence of  a $C^*$--structure     making the braiding unitary \cite{Wenzl}. See  also \cite{Xu} for an independent related result and    \cite{Sawin, Rowell1, Rowell2, Rowell3} and references therein  for a  complete discussion on unitarity and modularity.  In particular these results are of interest in the algebraic  approach to low dimensional QFT where tensor  
 $C^*$--categories with unitary braided symmetries arise \cite{Haag},  and constituted  our original motivation.
    
     The aim of this paper is  to        construct   $C^*$--quantum groupoids with representation category tensor equivalent to   the fusion $C^*$--categories of ${\mathfrak sl}_N$. They arise as  linear duals of quantised function algebras, and we thus denote them by   $\widehat{{\mathcal C}({\rm SU}(N),\ell)}$.     Every such  groupoid  turns out to be described by a finite dimensional $C^*$-algebra which is 
   a quotient of $U_q({\mathfrak sl}_N)$ determined by the irreducible  $^*$--representations with  dominant weights in the open Weyl alcove.  Hence the dimensions of the representations are classical. Our constructions may be interpreted as a natural dual of the quotient construction at the  categorical level   as recalled above.

  The problem of constructing quantum groupoids describing  a given tensor category via   representation theory originates in the physics literature (see \cite{Haag} and references therein) and  subsequently it has  been vastly considered   also with   different motivations. 
 For the   $C^*$--categories ${\mathcal F}$,   the fusion rules   become close to those of  ${\mathfrak g}$ for $\ell$ large. Physical considerations lead to look for  quantisations of the compact real form $G$ of ${\mathfrak g}$ \cite{G}. 
The quotient construction defining ${\mathcal F}$ allows a more precise formulation of the problem.
Indeed,   it may   be understood as taking a tilting module $T$ to a maximal non-negligible summand $\overline{T}$ and a morphism to its compression   to corresponding summands. Hence   objects of ${\mathcal F}$
 can be viewed  as  Hilbert spaces of classical dimension, and   arrows as  linear maps between them. We may regard this association as an approximation of the usual embedding functor of $\text{Rep}(G)$ into the Hilbert spaces for $\ell$ large, and ask
  whether the former can rigorously   be  interpreted as  
 the embedding functor of a $C^*$--quantum groupoid. Since the non-negligible modules are after all the tilting representations of   $U_q({\mathfrak g})$ with positive quantum dimension, to construct this groupoid one is led to start with the semisimple quotient of $U_q({\mathfrak g})$ determined by   such representations.

  The first construction  of this kind  has been done in the physics literature
    by Mack and Schomerus in the early 90s, who were motivated by  certain models of rational conformal field theory  \cite{MS0, MS1, MS}. 
  Specifically, they      started with the quotient of 
  $U_q({\mathfrak sl}_2)$ alluded to above   and showed that   truncation of tensor products of non-negligible modules 
     leads to a non-unital and  quasi-coassociative coproduct.  Furthermore, they   introduced the general notion of a  weak quasi-Hopf $C^*$--algebra, as a generalisation  to the non-unital case of the     quasi-Hopf algebras previously introduced by  Drinfeld  \cite{Drinfeld_coboundary}.

 To the best of our knowledge  in the last two decades there has not been any progress in extending   Mack  and Scomerus construction to pairs $({\mathfrak g}, \ell)$ for other Lie algebras. Notice that such a generalisation would not be obvious, 
 as the construction of the coproduct in the ${\mathfrak sl}_2$ case relies on the special property of that Lie algebra that tensor products of all the irreducible representations are multiplicity free.
Furthermore, not much is known about  
  the relation between   Mack and Schomerus   groupoid and the quotient category of ${\mathfrak sl}_2$, although such a relation is certainly expected.

This may be partially explained by the fact that most later studies   have focused on the coassociative weak Hopf algebras of  \cite{BS1,
 BNS, Nill}. Indeed, soon after their introduction such algebras    were   shown to cover  several  physical models, including orbifold models   previously described by genuine quasi-Hopf structures, see e.g.  \cite{DPR, BS1, EN, PZ, Rehren,  NV}. 
 Furthermore, every semisimple  fusion category of a rather general kind was shown to correspond to a weak Hopf algebra \cite{Hayashi, Szlachanyi, Ostrik}.
  We should however remark that  when applied to 
  the fusion categories ${\mathcal F}$,
   the relation between  the reconstructed weak Hopf algebra and the original Lie algebra is not apparent.
  For example,  the set of dimensions of   the base coefficient algebras   is unbounded if the order of the root of unity becomes large.
  This is due to the fact that  Hayashi-Szlachanyi duality relies on the construction of an (associative) embedding  functor  which does not   reproduce the above truncation procedure.

Our groupoids turn out  to satisfy 
  axioms of the weak quasi-Hopf $C^*$-algebras, and we shall interchangeably refer to them in this way or as   quasi-coassociative $C^*$--quantum groupoids. More in detail, 
they  are finite dimensional $C^*$-algebras endowed with  a   quasi-coassociative  and non-unital coproduct and      a multiplicative counit.   They also have  an invertible  antimultiplicative antipode, and an $R$-matrix. Coassociativity failure is described by a Drinfeld associator explicitly derived from an iterative truncation procedure associated to powers of the vector representation. Tensorial properties of the negligible modules in the   tilting category reflect in  rather strong properties 
of the associator that we describe in Sect. 9. For example, it is an idempotent.
 Counit multiplicativity has the advantage  that representations act on vector spaces rather than bimodules. Furthermore, since our groupoids    are built from the classical group using only representations of the non-semisimple   counterpart  of positive quantum dimensions,   they naturally approximate the classical Lie group in the sense described above.

We next explain our approach.   While there is in general no canonical choice of  truncated tensor products, Wenzl was able to make one for   tensor products $V_\lambda\otimes V$ of an irreducible $V_\lambda$ of $U_q({\mathfrak g})$  in the open alcove by a suitable 
generating representation $V$, 
  specifically  chosen for each Lie type (fundamental representation). He then proved Kirillov conjecture by showing that  this subspace is a Hilbert space  under the restriction of an hermitian form
  obtained as a  deformation of the usual tensor product structure via the $R$-matrix.  The corresponding projection  then becomes selfadjoint.
  We are interested in a coherent iterative choice of a sequence of such projections $p_n: V^{\otimes n}\to V^{\underline{\otimes}n}$ 
  onto  truncated tensor powers of $V$. One thus obtains a natural   functor
$$W: {\mathcal F}\to{\mathcal H}$$ from    the fusion    $C^*$--category to the Hilbert spaces taking a tensor power $V^{\otimes n}$  to the range   of $p_n$.

  We adapt   Tannakian methods for    quantum groups    to the functor $W$. In other words, we pass to the dual viewpoint  and talk of comodules.   
  Notice   that   Wenzl's functor 
   is not a tensor $^*$--functor into the Hilbert spaces in the usual sense \cite{NT}. Even compared to more general Tannakian reconstructions (e.g. for weak Hopf algebras \cite{Szlachanyi}, or ergodic actions of compact quantum groups \cite{PR_ergodic}), there are   independent
    new obstructions.   First, the usual tensor product of Hilbert space representations is not a representation on the tensor product Hilbert space, already if $q$ is not a root of unity. This is due to the anticomultiplicativity property of the involution  of $U_q({\mathfrak g})$. This means that the construction of the adjoint can not be treated along  lines similar, e.g.,  to those for  compact quantum groups 
 \cite{W}, as the involution is comultiplicative for the latter.   At roots  of unity, a second obstacle is associativity failure of   $W$   and it is due to   the truncation procedure.   Indeed, Wenzl's   projections $p_n$     
  fail associativity already for ${\mathfrak sl}_2$ at the smallest root, 
 in agreement with   Mack and Schomerus observation, cf. Ex. 5.2.      In essence,  our contribution may be regarded  as the construction of a tensor structure of a weak type on Wenzl's functor. Specifically, it turns out to be a weak quasi-tensor structure in the sense of \cite{HO}, but of a rather special kind, in that it satisfies suitable weak tensoriality relations which are generalizations of an ordinary notion of weak tensor functor relying on coassociativity relations. This viewpoint will be clarified  in    \cite{CP2} where we put the results of this paper in an abstract framework. 
  
To obtain this, we replace associativity failure of $W$ by certain remarkable properties of the negligible modules of the tilting category, due to Andersen \cite{Andersen} and  emphasised by Gelfand and Kazhdan \cite{GK} in their quotient construction. We show that these properties play a   role also in our case: they allow us to  construct a coassociative but non-associative  bi-algebra ${\mathcal D}(V,\ell)$ acting as  a universal algebra of a vector space.
We endow it with an  involution induced   by Drinfeld's coboundary associated to the so called unitarized $R$-matrix $\overline{R}$ of \cite{Drinfeld_coboundary}.
  We next introduce the `function algebra' ${\mathcal C}(G,\ell)$ as a   quotient of ${\mathcal D}(V,\ell)$ by   a coideal induced by   the fusion category which, by our coherent choice of the $p_n$ is also an ideal, but only of one-sided type. Hence ${\mathcal C}(G,\ell)$ is naturally only a coalgebra with involution.  The   problem becomes that  of making it into an algebra. Notice that this corresponds to the most delicate point in Mack and Schomerus approach, namely the construction of a coproduct.

We observe  that our problem can be reduced to the question of cosemisimplicity of ${\mathcal C}(G,\ell)$. Specifically, we show that if we know that  ${\mathcal C}(G,\ell)$ is cosemisimple
 then one can derive a non-associative algebra structure on it  as a pull back of    the product of ${\mathcal D}(V,\ell)$. But more is true:  the dual groupoid $\widehat{{\mathcal C}(G,\ell)}$ can be made into a weak quasi-Hopf $C^*$--algebra with $R$-matrix with representation category equivalent to ${\mathcal F}$.  The associator is explicitly described by associativity failure of Wenzl's projections $p_n$. It has the virtue of being an idempotent.   Our second main result  is that ${\mathcal C}(G,\ell)$ is indeed cosemisimple for $G={\rm SU}(N)$. This is done via the construction of a suitable   filtration of ${\mathcal C}(G,\ell)$ in the general case, for which an (even associative) algebra structure is naturally defined, and of a    Haar functional for that filtration in the type $A$ case. This functional     is in turn
  achieved via an analysis on    the Weyl filtrations of the negligible summands   of the tensor powers $V^{\otimes n}$ of the fundamental up to a specific    value of $n$ specified by   the conjugation structure of the representation category of ${\rm SU}(N)$.  
   
     We should remark that because of space limitations,     cosemisimplicity   of ${\mathcal C}(G,\ell)$ for   Lie types other than $A$ has not been considered in this paper, but will be taken up in the future. A positive answer would allow us to extend our main result to   such groups 
  provided ${\mathfrak g}\neq E_8$.

   The paper is organised as follows. Sect. 2 is dedicated to preliminaries on quantum groups at roots of unity and  the associated fusion categories. In Sect. 3 we recall Drinfeld's unitarized $R$-matrix and the associated coboundary.  In Sect. 4 we recall Kirillov-Wenzl's theory mostly following Wenzl's approach, while in Sect. 5 we construct Wenzl's functor. We dedicate Sect. 6 to  Tannakian reconstruction of the function algebra quantum group  in the generic case (i.e. $q$ is not a root of unity), while in Sect. 7 we perform the construction of the universal non-associative   algebra ${\mathcal D}(V,\ell)$ at roots of unity   and we show   multiplicativity of the coproduct for  all ${\mathfrak g}\neq E_8$. Section 8 deals with the construction of  the involutive function coalgebra ${\mathcal C}(G,\ell)$ and an associated  filtration endowed with an associative product.     In Sect. 9 we construct quasi-coassociative dual $C^*$-groupoids $\widehat{{\mathcal C}(G,\ell)}$ under the assumption of cosemisimplicity of ${\mathcal C}(G,\ell)$,   we study Drinfeld's associator, verify
   the quasi-Hopf $C^*$-algebra axioms, and show the mentioned equivalence of tensor $C^*$--categories.
   Finally, in Sect. 10 we develop a cosemisimplicity condition for  ${\mathcal C}(G,\ell)$ involving a Haar functional on    the   filtration, and we verify its validity for ${\mathfrak g}={\mathfrak sl}_N$. We include an appendix where we explicit the generator for the fusion $C^*$--categories of  ${\mathfrak sl}_N$  and the conjugation structure.

\section{Quantum groups at roots of unity}

Let ${\mathfrak g}$ be a complex simple Lie algebra, ${\mathfrak h}$ a Cartan subalgebra, $\alpha_1,\dots,\alpha_r$ a set
of simple roots, and $A=(a_{ij})$ the associated Cartan matrix. Consider the unique invariant symmetric and bilinear form on ${\mathfrak h}^*$ such that $\langle\alpha,\alpha\rangle=2$
for a short root $\alpha$. Let
$E$ be the real vector space generated by the roots endowed with its euclidean structure $\langle x, y\rangle$. Let
 $\Lambda$ be the weight lattice of $E$ and $\Lambda^+$ the cone of dominant weights.

Consider  the complex $^*$--algebra 
${\mathbb C}[x,x^{-1}]$ 
of   Laurent polynomials with involution
$x^*=x^{-1}$, and let ${\mathbb C}(x)$ be the associated quotient field, endowed with the involution naturally induced from ${\mathbb C}[x,x^{-1}]$. 
 We consider Drinfeld-Jimbo quantum group     $U_x({\mathfrak g})$, i.e. the algebra over ${\mathbb C}(x)$ defined by generators
$E_i$, $F_i$, $K_i$, $K_i^{-1}$, $i=1,\dots, r$, and relations  
$$K_iK_j=K_jK_i,\quad K_iK_i^{-1}=K_i^{-1}K_i=1,$$
$$K_iE_jK_i^{-1}=x^{\langle\alpha_i,\alpha_j\rangle}E_j,\quad K_iF_jK_i^{-1}=x^{-\langle\alpha_i,\alpha_j\rangle}F_j,$$
$$E_iF_j-F_jE_i=\delta_{ij}\frac{K_i-K_i^{-1}}{x^{d_i}-x^{-d_i}},$$
$$\sum_0^{1-a_{ij}}(-1)^kE_i^{(1-a_{ij}-k)}E_j E_i^{(k)}=0,\quad\sum_0^{1-a_{ij}}(-1)^kF_i^{(1-a_{ij}-k)}F_j F_i^{(k)}=0,\quad i\neq j,$$
where $d_i=\langle\alpha_i,\alpha_i\rangle/2$, and, for $k\geq0$, $E_i^{(k)}=E_i^k/[k]_{d_i}!$,  $F_i^{(k)}=F_i^k/[k]_{d_i}!$. Note that $d_i$ is an integer, hence so is
every inner product $\langle\alpha_i,\alpha_j\rangle$.  Quantum integers and factorials are defined in the usual way, $[k]_x=\frac{x^k-x^{-k}}{x-x^{-1}}$; $[k]_x!=[k]_x\dots[2]_x$, $[k]_{d_i}=[k]_{x^{d_i}}$, and result selfadjoint scalars of ${\mathbb C}(x)$. There is a unique $^*$--involution on $U_x({\mathfrak g})$ making it into a $^*$--algebra over ${\mathbb C}(x)$ such that     $$K_i^*=K_i^{-1},\quad E_i^*=F_i.$$  

This algebra becomes a Hopf algebra, i.e. a coassociative coalgebra with coproduct $\Delta$, counit $\varepsilon$ and antipode $S$ defined, e.g, as in \cite{Wenzl}, where his $\tilde{K}_i$ corresponds to our $K_i$.

One has the following relations between coproduct, antipode and involution for $a\in U_x({\mathfrak g})$,
$$\Delta(a^*)= \Delta^{\text{op}}(a)^*, \eqno(2.1)$$
$$\varepsilon(a^*)=\overline{\varepsilon(a)}. \eqno(2.2)$$
$$S(a^*)=S(a)^*,\eqno(2.3)$$
$$S^2(a)=K_{2\rho}^{-1}aK_{2\rho},\eqno(2.4)$$
where $ \Delta^{\text{op}}$ is the coproduct opposite to $\Delta$, $2\rho$  the sum of the positive roots, and, for an element $\alpha=\sum_i k_i\alpha_i$ of the root lattice, $K_\alpha:=K_1^{k_1}\dots K_r^{k_r}$.\medskip

\noindent{\bf 2.1. Remark} 
Notice  that $U_x({\mathfrak g})$  is not a Hopf $^*$--algebra in the   sense of, e.g.,  \cite{CP}, where    $(2.1)$--$(2.3)$ are  replaced by $\Delta(a^*)=\Delta(a)^*$ and  
$S(a^*)=S^{-1}(a)^*$.\medskip

We next  consider Lusztig's integral form \cite{Lusztig}. 
In order to construct braided tensor categories, it is well known that we   need to embed the original algebra into a larger algebra, via a procedure involving, among other things, extension of scalars, see \cite{Sawin} for details. For the purposes of this paper (see the next section) we shall actually need  a further extension, and the correct polynomial ring for the integral form will be 
 $${\mathcal A}:={\mathbb C}[x^{1/2L}, x^{-1/2L}],$$ with $L$ the smallest positive integer such that
$L\langle\lambda, \mu\rangle\in{\mathbb Z}$ for all dominant weights $\lambda$, $\mu$.
The explicit values are listed in \cite{Sawin} for all Lie types. For example,   $L=N$ for ${\mathfrak g}={\mathfrak sl}_N$.

 We define the integral form   ${\mathcal U}_{\mathcal A}$ as
  the ${\mathcal A}$--subalgebra generated by the elements $E_i^{(k)}$, $F_i^{(k)}$ and $K_i$.
This is known to be a $^*$--invariant Hopf ${\mathcal A}$--algebra with the structure inherited from $U_x({\mathfrak g})$. Notice that, in connection with the mentioned extension needs, ${\mathcal U}_{\mathcal A}$ is not quasi triangular, even topologically, but the representation categories that we next consider will be braided tensor categories, and this will suffice 
  for our purposes.

We fix $q\in{\mathbb T}$, and consider the $^*$--homomorphism ${\mathcal A}\to {\mathbb C}$ which evaluates every polynomial in $q$, and form the tensor product
  $^*$--algebra,
$$U_q({\mathfrak g}):={\mathcal U}_{\mathcal A}\otimes_{\mathcal A}{\mathbb C},$$ which becomes
a complex Hopf algebra with a $^*$--involution, and properties $(2.1)$--$(2.4)$ still hold for $U_q({\mathfrak g})$.\medskip

 Given a   dominant weight $\lambda$    of ${\mathfrak g}$, we can associate various modules    $V_\lambda(x)$, $V_\lambda({\mathcal A})$, and $V_\lambda(q)$ of    $U_x({\mathfrak g})$,  ${\mathcal U}_{\mathcal A}$ and $U_q({\mathfrak g})$,  
respectively, usually called Weyl modules, and thus form corresponding representation categories as follows. We shall mostly be interested in $V_\lambda(q)$ that we   shall   eventually     simply denote by $V_\lambda$ as well.
\medskip

\noindent{\it a) The category $\text{Rep}({\mathcal U}_{\mathcal A})$}\medskip

\noindent Let $V_\lambda(x)$ be
the  irreducible representation of $U_x({\mathfrak g})$ with highest weight $\lambda$ and let $v_\lambda$ be the highest weight  vector of $V_\lambda(x)$ and form the  cyclic  module 
        of ${\mathcal U}_{\mathcal A}$ generated by     $v_\lambda$,
$$V_\lambda({\mathcal A}):=u_\lambda({\mathcal U}_{\mathcal A})v_\lambda.$$
It is known to be  a free ${\mathcal A}$--module  satisfying 
$$V_\lambda({\mathcal A})\otimes_{\mathcal A}{\mathbb C}(x)=V_\lambda(x).$$
We denote by $\text{Rep}({\mathcal U}_{\mathcal A})$   the linear category over ${\mathcal A}$ with objects finite tensor products of modules
$V_{\lambda}({\mathcal A})$. It becomes a tensor ${\mathcal A}$--category in the natural way.
\medskip

\noindent{\it b) The tilting category ${\mathcal T}({\mathfrak g}, \ell)$ }\medskip

\noindent Every module $V_{\lambda}({\mathcal A})$ gives rise to
the complex  ${ U}_q({\mathfrak g})$-modules via specialisation at a complex number $q$:
$$V_\lambda(q):=V_\lambda({\mathcal A})\otimes_{\mathcal A}{\mathbb C}.$$ 
This is obviously a cyclic module, but it is not always
  irreducible if $q$ is a root of unity. 
    The   {\it linkage principle} gives information on irreducibility  at primitive roots of unity. Briefly, 
one needs to consider the   {\it  affine} Weyl group,  and   its  {\it translated } action on the real vector space $E$ spanned by the roots, defined by $w.x=w(x+\rho)-\rho$.  The structure of this group depends on the parity and divisibility by $d:=\max{d_i}$ of the order  $m$ of $q$ \cite{Sawin}.
 In this paper we are interested in the case  where  $m$   is divisible by $2d$, and more specifically we take 
 $q=e^{\pi i/d\ell}$, as in this case one obtains tensor $C^*$--categories  \cite{Wenzl, Xu}. 
 The affine Weyl group,   $W_\ell$, is generated by the ordinary Weyl group $W$ and translations by $\ell\theta$, where $\theta$ is the highest root. The translated action   admits a fundamental domain, called the principal Weyl alcove, that intersects $\Lambda^+$   in the set 
 $$\overline{\Lambda}_{\ell}:=\{ \lambda\in\Lambda^+: \langle   \lambda+\rho, \theta\rangle \leq d\ell\}.$$ 
 The linkage principle then    implies that $V_\lambda(q)$ are  irreducible   for  
  $\lambda\in \overline{\Lambda}_{\ell}$. Moreover, they are pairwise inequivalent.
 We refer to \cite{Andersen, Sawin, CP} for   complete  explanations.

We follow \cite{Wenzl} for the construction of the tilting category. Namely, we fix, for each Lie type,
a   representation $V$ of ${\mathfrak g}$ taken from a specific list, that we call {\it fundamental}. For example, in the type $A$ case, $V$ is the vector representation. Each fundamental representation 
has in particular the   property that every irreducible of ${\mathfrak g}$ is contained in some tensor  power of $V$.

One can form the  category ${\mathcal T}({\mathfrak g}, \ell)$, also denoted ${\mathcal T}_\ell$ for brevity,  with objects finite tensor powers of $V(q)$
and arrows   the intertwining operators,   
completed with subobjects (i.e. summands) and direct sums.

 This is a strict tensor but non-semisimple category. It is known \cite{Wenzl} that 
the objects of ${\mathcal T}_\ell$ are tilting modules in the sense of Andersen \cite{Andersen}, and that conversely, for $\ell$ large enough
every tilting module 
is isomorphic to an object of ${\mathcal T}_\ell$. More precisely, one needs an order $\ell$ such that   $\kappa\in \Lambda_\ell$, where $\kappa$ is the dominant weight of $V$ and
$$\Lambda_{\ell}:=\{ \lambda\in\Lambda^+: \langle   \lambda+\rho, \theta \rangle < d\ell\}.$$
Tilting modules were originally  defined  as those modules $W$  admitting, together with their dual,
a Weyl filtration, i.e. a sequence of modules $\{0\}\subset W_1\subset\dots \subset W$ such that
$W_{i+1}/W_i$ is isomorphic to a Weyl module $V_{\lambda_i}(q)$ with $\lambda_i\in\Lambda^+$.
  Weyl filtrations are non-unique, but for all filtrations of $W$ the number of successive factors isomorphic to a given 
$V_\lambda(q)$  is unique, and it is in fact given by the multiplicity of $V_\lambda(x)$ in $W(x)$ if $W$ is obtained from a specialisation  $x\to q$ of a module 
$W(x)$ of $U_x({\mathfrak g})$, see Prop. 3 and Remark 2 in  \cite{Sawin} for a  precise statement. This in particular implies  that   the multiplicities of the dominant weights of the factors in the Weyl filtrations of   tensor products $V_{\lambda_1}(q)\otimes\dots\otimes V_{\lambda_n}(q)$   with   $\lambda_i\in\overline{\Lambda}_\ell$ (or more generally of tilting modules) are the same as those determined by decomposition into irreducibles of the corresponding tensor product in the classical (or generic) case.  
\bigskip

\noindent{\it c) The quotient category ${\mathcal F}_\ell$}\bigskip

\noindent  We next briefly recall the construction of  the semisimple quotient of ${\mathcal T}_\ell$,   henceforth   written  ${\mathcal F}_\ell$, following the approach of Gelfand and Kazhdan \cite{GK}.
Every object of ${\mathcal T}_\ell$  
  decomposes as a direct sum of indecomposable submodules, and this decomposition is unique up to isomorphism.
 One can   form two full linear (non-tensorial) subcategories, ${\mathcal T}^0$, and ${\mathcal T}^\perp$ of ${\mathcal T}_\ell$, with objects, respectively, 
  those representations which can be written as direct sums of $V_\lambda(q)$, with $\lambda\in\Lambda_\ell$ only,
   and those which have no such $V_\lambda(q)$ as a direct summand.\medskip
   
   \noindent{\bf 2.2. Definition} The objects of ${\mathcal T}^\perp$ and ${\mathcal T}^0$ are called  negligible and non-negligible, respectively. An arrow $T:W\to W'$ of ${\mathcal T}_\ell$ is called negligible 
if it factors through $W\to N\to W'$ with $N$ negligible. \medskip

The category ${\mathcal T}^\perp$ of negligible modules satisfies the following properties, first shown by  Andersen \cite{Andersen}, and abstracted in \cite{GK}  \medskip
 
 \noindent (1) Any object $W\in{\mathcal T}_\ell$ is isomorphic to a direct sum $W\simeq W_0\oplus N$ with $W_0\in{\mathcal T}^0$ and
  $N\in{\mathcal T}^\perp$.\medskip
  
  \noindent(2) For any pair of arrows $T:W_1\to N$, $S:N\to W_2$ of ${\mathcal T}_\ell$, with $N\in{\mathcal T}^\perp$, $W_i\in{\mathcal T}^0$, then $$S T=0.$$ 
  
  \noindent (3) For any pair of objects   $W\in{\mathcal T}_\ell$, $N\in{\mathcal T}^\perp$, then   
  $W\otimes N$ and $N\otimes W\in{\mathcal T}^\perp$.\medskip

Property $(1)$ follows immediately from the mentioned decomposition of objects of ${\mathcal T}_\ell$, while property $(2)$ means that no non-negligible module can be a summand of a negligible one (however, it can be a factor of a Weyl filtration of a negligible).

Let $\text{Neg}(W, W')$ be the subspace of negligible arrows of   $(W, W')$.
Then the quotient category, ${\mathcal F}_\ell$, is the category with the same objects as ${\mathcal T}_\ell$ and arrows between the objects $W$ and $W'$ the quotient space,
$$(W, W')_{{\mathcal F}_\ell}:=(W, W')/\text{Neg}(W, W').$$
 Gelfand and Kazhdan endow ${\mathcal F}_\ell$ with the  unique structure of a tensor category   such that the quotient map ${\mathcal T}_\ell\to{\mathcal F}_\ell$ is a tensor functor.  The   tensor product  of  objects and  arrows of ${\mathcal F}_\ell,$ is usually denoted   by 
  $W\underline{\otimes} W'$
 and
 $S\underline{\otimes} T$ respectively,
 and referred to  as the truncated tensor product  in the physics literature. This is now a
  semisimple tensor category and 
  $\{V_\lambda, \, \lambda\in\Lambda_\ell\}$ is a complete set of irreducible objects.
 
 \medskip
 
 \noindent{\bf 2.3. Remark} 
 By Lemma 1.1 in \cite{GK}, composition of  inclusion   ${\mathcal T}^0\to{\mathcal T}_\ell$ with   projection   ${\mathcal T}_\ell\to{\mathcal F}_\ell$ is an equivalence of linear categories. Hence ${\mathcal T}^0$ becomes
a semisimple tensor category as well tensor equivalent to ${\mathcal F}_\ell$.
 This  procedure may be understood as a categorification of the following visualisation of truncated tensor products at the level of Grothendieck rings.  
 For $\lambda$, $\mu\in\Lambda_{\ell}$, one can decompose $V_\lambda\otimes V_\mu$ uniquely up to isomorphism in ${\mathcal T}_\ell$,
 $$V_\lambda\otimes V_\mu\simeq \oplus_{\nu\in\Lambda_\ell} m^\nu_{\lambda, \mu} V_\nu \oplus N,$$
 with $N$ negligible. Then in ${\mathcal F}_\ell$,
 $$V_\lambda\underline{\otimes } V_\mu\simeq \oplus_{\nu\in\Lambda_\ell} m^\nu_{\lambda, \mu} V_\nu.$$
 Notice that, although unique up to isomorphism, the decomposition of $V_\lambda\otimes V_\nu$  described in ${\mathcal T}_\ell$ is not canonical
 (cf. also Sect. 11.3C in \cite{CP}  and references therein.)
 \medskip

\bigskip

\section{Ribbon and coboundary structures}

The main topic of   this section       is the coboundary structure of the category $\text{Rep}({\mathcal U}_{\mathcal A})$ associated to      the ribbon structure.  It was first introduced by Drinfeld \cite{Drinfeld_coboundary} for different purposes and used by Wenzl to obtain a variant of Kirillov inner product
\cite{Wenzl}. We shall later need it as well for the construction of the $^*$--involution of the groupoid.

In order to  obtain ribbon Hopf algebras, one needs to pass to a suitable extension of ${\mathcal U}_{\mathcal A}$.
We  refrain from explicitly  recalling the construction of $R$ and we refer to \cite{Sawin} for a detailed treatment, but our notation will agree with \cite{Wenzl}. 
The $R$ matrix then lies in  a suitable topological and algebraic completion of the of 
the tensor product   algebra by itself.
  As anticipated in the previous section, we slightly modify this construction in that we further extend the scalars to ${\mathcal A}={\mathbb C}[x^{1/2L}, x^{-1/2L}]$.
This will allow the construction of a central square root  of the ribbon element in the completed algebra which will be useful later for  the hermitian structure.  

We recall the main algebraic properties of $R$ \cite{Drinfeld_quasi_triangular},
$$ \Delta^{\text{op}}(a)=R\Delta(a)R^{-1},\eqno(3.1)$$
$$1\otimes\Delta(R)=R_{13}R_{12},\eqno(3.2)$$
$$\Delta\otimes 1(R)=R_{13}R_{23}.\eqno(3.3)$$
Relations $(3.1)$--$(3.3)$ mean  that for any pair of representations $u$, $v$ of $\text{Rep}({\mathcal U}_{A})$, the   operators $\varepsilon_{u, v}:=\Sigma u\otimes v(R)$, with $\Sigma$ the flip map, are intertwiners of the category satisfying    naturality in $u$ and $v$. Hence    $\varepsilon_{u, v}$ is a braided symmetry for the    category. In particular,
the Yang-Baxter relation
$$R_{12}R_{13}R_{23}=R_{23}R_{13}R_{12}$$
follows,
see e.g. \cite{Turaev}  for details.

There is an associated  ribbon element which is
  a central  invertible element $v$ in the completed Hopf algebra such that
$$R_{21}R=v\otimes v\Delta(v^{-1}),$$ and it is 
  given by  
$v=K_{2\rho}u,$
where $u=m(S\otimes 1(R_{21})),$ and often called  the quantum Casimir operator.
The element $u$ was originally introduced in
\cite{Drinfeld_quasi_triangular}
  the notion of ribbon Hopf algebras and ribbon tensor category is due to Reshetikhin and Turaev \cite{RT2, RT}.
For the following important result, see Theorem 3 in \cite{Sawin}  and references therein.\medskip

\noindent{\bf 3.1. Theorem} {\sl The category $\text{Rep}({\mathcal U}_{\mathcal A})$ is a ribbon tensor category. Hence so are ${\mathcal T}_\ell$ and ${\mathcal F}_\ell$.}\medskip

We next consider the  relation between the $R$ matrix and the   $^*$--involution. Wenzl proves the following    \cite{Wenzl},
$$R^*={R_{21}}^{-1}.$$
On the other hand it is known that the inverse of $u$ can be computed as 
$$u^{-1}=m(S^{-1}\otimes 1((R^{-1})_{21})).$$
It follows that
$$u^*=m(1\otimes S(R^*))=m(S^{-1}\otimes 1(R^*))=u^{-1},$$
for the second equality, see Lemma 2.1.1 in \cite{Turaev}.  
Therefore
$$v^*=v^{-1}$$ as well.
The action of $v$ on an irreducible module is derived from the expression of the $R$ matrix, and it is given by   scalar multiplication 
by $x^{-\langle\lambda,\lambda+2\rho\rangle}$ on a highest weight module  of weight $\lambda$, 
\cite{Drinfeld_quasi_triangular}, see also see \cite{Wenzl, CP}. By \cite{Drinfeld_coboundary}
one can construct   a central square root $w$ lying in the topological completion of the extension of ${\mathcal U}_{\mathcal A}$. It acts as the scalar $x^{-\frac{1}{2}\langle\lambda,\lambda+2\rho\rangle}$
on the same module. One thus  has:
$$w^2=v,$$
$$w^*=w^{-1}.$$
It follows that one can compute a   square root of $R_{21}R$,
$$(R_{21}R)^{1/2}=w\otimes w \Delta(w^{-1}).$$
Drinfeld used the element $w$   to construct  the so called        unitarized $R$ matrix
$\overline{R}$ as follows.  Set 
$$\Theta:=(R_{21}R)^{-1/2}=w^{-1}\otimes w^{-1}\Delta(w),$$
and define 
   $$\overline{R}:=R\Theta=Rw^{-1}\otimes w^{-1}\Delta(w).$$
 We recall the following     relations.\medskip
   
   \noindent{\bf 3.2. Lemma} {\sl 
   \begin{itemize}
\item[{\rm a)}] $\Theta^*=(\Theta)_{21}^{-1},$ 
   \item[{\rm b)}] $\Theta_{21}R=R\Theta$.
   \end{itemize}
      
   }\medskip
   
\noindent{\sl proof}   
 a) $\Theta^*=\Delta(w)^*w\otimes w= \Delta^{\text{op}}(w^*)w\otimes w=w\otimes w \Delta^{\text{op}}(w^{-1})=\Theta_{21}^{-1}.$
b) By $(3.1)$  and centrality of $w$,
 $$\Theta_{21}R=w^{-1}\otimes w^{-1} \Delta^{\text{op}}(w)R=w^{-1}\otimes w^{-1}R\Delta(w)=Rw^{-1}\otimes w^{-1}\Delta(w)=R\Theta.$$

  \medskip

 \noindent{\bf  3.3. Proposition} {\sl One has:
 $$\overline{R}^*=(\overline{R})_{21}^{-1}=\overline{R}.$$ In particular, $\overline{R}$ is selfadjoint.}\medskip

\noindent{\sl Proof} By a) of the previous lemma,
$$\overline{R}^*=(R\Theta)^*=\Theta^*R^*=(\Theta)_{21}^{-1}R_{21}^{-1}=(R_{21}\Theta_{21})^{-1}=(\overline{R})_{21}^{-1}.$$ Furthermore, by b) of the same lemma, 
$$\overline{R}_{21}\overline{R}=(R\Theta)_{21}R\Theta=R_{21}\Theta_{21}R\Theta=R_{21}R\Theta^2=1.$$\medskip

  The Yang-Baxter  equation can  equivalently be read as an associativity property for the $R$ matrix:
$$R_{12}\Delta\otimes 1(R)=R_{23}1\otimes\Delta(R).$$
 One can thus unambiguously define elements in the iterated tensor powers of the ribbon Hopf algebra,
$${R}^{(n+1)}:=R^{(n)}\otimes 1\Delta^{(n-1)}\otimes 1(R)=1\otimes R^{(n)}1\otimes\Delta^{(n-1)}(R)$$
where ${R}^{(2)}:={R}$.
We set, for $n>2$,
$$\Theta^{(n)}:=w^{-1}\otimes\dots \otimes w^{-1}\Delta^{(n-1)}(w)$$
and note that it
is associative as well:
$$\Theta^{(n)}\otimes 1\Delta^{(n-1)}\otimes 1(\Theta)=w^{-1}\otimes\dots\otimes w^{-1}\otimes 1\Delta^{(n-1)}(w)\otimes 1\Delta^{(n-1)}\otimes 1(w^{-1}\otimes w^{-1}\Delta(w))=$$
$$\Theta^{(n+1)}=1\otimes\Theta^{(n)}1\otimes\Delta^{(n-1)}(\Theta).$$
We also set:
$${\overline{R}}^{(n)}:=R^{(n)}\Theta^{(n)}.$$\medskip

   \noindent{\bf 3.4. Lemma} {\sl 
The following pairs of operators commute:
   \begin{itemize}
   \item[{\rm a)}] $\Delta^{(n)}\otimes1(R)$ and  $\Theta_{12}$,
        \item[{\rm b)}] $1\otimes\Delta^{(n)}(R)$ and  $\Theta_{23}$,
        \item[{\rm c)}] $\Delta^{(n-1)}\otimes 1(R)$ and $\Theta^{(n)}\otimes 1$,
        \item[{\rm d)}] $1\otimes \Delta^{(n-1)}(R)$ and $1\otimes\Theta^{(n)}$.
   \end{itemize}
   }\medskip
   
 \noindent{\sl proof}    a) we need to show that $\Delta^{(n)}\otimes 1(R)$ commutes with $R_{21}R=\Theta^{-2}$. We explicit $\Delta^{(n)}$ by  leaving $\Delta$ always on the left: $\Delta^{(2)}=\Delta\otimes 1\circ\Delta$,   $\Delta^{(3)}=\Delta\otimes1\otimes 1 \circ\Delta\otimes 1\circ\Delta$, and so on. 
An iteration of $(3.3)$ gives
$$\Delta^{(n)}\otimes 1(R)=R_{1\, n+2}R_{2\, n+2}\dots R_{n+1\, n+2}.$$
$R_{21}R_{12}$ obviously commutes with all the factors in this expression except the first two. For them, we use the 
 twice the Yang-Baxter equation, the first time   with $1$ and $2$ exchanged, in both cases   $3$ is replaced by $n+2$:
$$(R_{1\, n+2}R_{2\, n+2})(R_{21}R)=(R_{1\, n+2}R_{2\, n+2}R_{21})R=(R_{21}R_{2\, n+2}R_{1\, n+2})R_{12}=$$
$$R_{21}(R_{2\, n+2}R_{1\, n+2}R_{12})=R_{21}(R_{12}R_{1\,n+2}R_{2\, n+2})=(R_{21}R)(R_{1\,n+2}R_{2\, n+2})$$
 c) follows from a) and the fact that $$\Theta^{(n)}=\Theta\Delta\otimes 1(\Theta)\dots\Delta^{(n-2)}\otimes 1(\Theta).$$
 b) and d) can be proved in a similar way.\medskip

\noindent{\bf 3.5. Proposition} {\sl
  $\overline{R}^{(n)}$ is associative:
 $${\overline{R}}^{(n+1)}=
 \overline{R}^{(n)}\otimes 1\Delta^{(n-1)}\otimes 1(\overline{R})=1\otimes \overline{R}^{(n)}1\otimes\Delta^{(n-1)}(\overline{R}).$$}\medskip
 
 \noindent{\sl Proof} This follows from properties c) and d) of the previous lemma.\medskip
 
\noindent{\bf 3.6. Proposition} {\sl $\overline{R}^{(n)}$ is selfadjoint for all $n$.}\medskip

\noindent{\sl Proof} For all $n$, $$\Delta^{(n)}(a)^*={ \Delta^{\text{op}}}^{(n)}(a^*).$$
Furthermore, being $w$ central, we may replace $R$ with $\overline{R}$  in $(2.6)$:
$$ \Delta^{\text{op}}(a)\overline{R}=\overline{R}\Delta(a).$$
It follows that, for all $n$,
$${ \Delta^{\text{op}}}^{(n)}(a)\overline{R}_{n+1}=\overline{R}_{n+1}\Delta^{(n)}(a),$$
where
$\overline{R}_{n+1}:={ \Delta^{\text{op}}}^{(n-1)}\otimes 1(\overline{R})\dots \Delta^{\text{op}}\otimes 1(\overline{R})\overline{R}$.
We thus have, being $\overline{R}$ selfadjoint, 
$$(\overline{R}^{(n)})^*=(\overline{R}\Delta\otimes 1(\overline{R})\dots\Delta^{(n-2)}\otimes 1(\overline{R}))^*=$$
$${ \Delta^{\text{op}}}^{(n-2)}\otimes 1(\overline{R})\dots \Delta^{\text{op}}\otimes 1(\overline{R})\overline{R}=\overline{R}_{n}=$$
$${ \Delta^{\text{op}}}^{(n-2)}\otimes 1(\overline{R}){\overline{R}}_{n-1}={\overline{R}}_{n-1}{\Delta}^{(n-2)}\otimes 1(\overline{R}).$$
By induction on $n$, $\overline{R}^{(n-1)}=(\overline{R}^{(n-1)})^*={\overline{R}}_{n-1}.$ Inserting this information in the above computation gives
$$(\overline{R}^{(n)})^*=\overline{R}^{(n-1)}{\Delta}^{(n-2)}\otimes 1(\overline{R})={\overline{R}}^{(n)}.$$
\medskip

We next introduce Drinfeld's coboundary of $\text{Rep}({\mathcal U}_{\mathcal A})$.
First consider, for any pair of objects $u$, $v\in\text{Rep}({\mathcal U}_{\mathcal A})$, the braiding operators:
$$\varepsilon(u, v)=\Sigma u\otimes v(R)\in (u\otimes v, v\otimes u).$$
As recalled above, these are natural isomorphisms in $u$ and $v$ and define a braided symmetry, in that
$$\varepsilon(u, w)\otimes 1_v\circ 1_u\otimes\varepsilon(v, w)=\varepsilon(u\otimes v, w),$$
$$1_v\otimes\varepsilon(u, w)\circ\varepsilon(u, w)\otimes 1_w=\varepsilon(u, v\otimes w).$$
Correspondingly, we consider the associated modified form   
$$\sigma(u, v)=\Sigma u\otimes v(\overline{R})\in(u\otimes v, v\otimes u),$$
 obviously still arrows of the  category,
However, they do not define a braided symmetry. Rather, one has the following properties, which are immediate consequences  of the previous propositions.
\medskip

\noindent{\bf 3.7. Proposition} {\sl  The arrows 
$\sigma(u, v)\in(u\otimes v, v\otimes u)$ of $\text{Rep}({\mathcal U}_{\mathcal A})$
are natural isomorphisms satisfying
$$\sigma(v, w)\circ\sigma(w,v)=1_{w\otimes v},\eqno(3.4)$$
$$\sigma(v\otimes u, w)\circ\sigma(u,v)\otimes 1_w=\sigma(u,  w\otimes v)\circ 1_u\otimes\sigma(v, w).\eqno(3.5)$$}\medskip

An abstract tensor category admitting arrows $\sigma(u, v)$ satisfying the properties stated in the last proposition, is called a coboundary category
\cite{Drinfeld_coboundary}.

Notice that both sides of  $(3.5)$ define an intertwiner of the category which reverses the order in   triple tensor products, 
$$\sigma_3: v_1\otimes v_2 \otimes v_3\to v_3\otimes v_2\otimes v_1.$$ 
More generally, we consider arrows
$$\sigma_n\in(v_1\otimes\dots\otimes v_n,  v_n\otimes\dots\otimes v_1)$$
inductively  defined by 
$$\sigma_n:=\sigma(v_{n-1}\otimes\dots\otimes v_1, v_n)\circ\sigma_{n-1}\otimes 1_{v_n}.$$
\medskip
We explicit $\sigma_n$ dropping the involved representations in the notation. 
Let $\Sigma_n$ denote the permutation that reverses the order in a tensor product space with $n$ factors.
\medskip

\noindent{\bf 3.8. Proposition} {\sl 

\begin{itemize}
   \item[{\rm a)}] 
   $\sigma_n=\Sigma_n \overline{R}_n$,
      \item[{\rm b)}] $\sigma_n=
      \sigma_{n-1}\otimes 1_{v_1}\circ\sigma(v_1, v_2\otimes\dots\otimes v_{n})=$
      
 \noindent     $\sigma(v_1, v_n\otimes\dots\otimes v_2)\circ1_{v_1}\otimes\sigma_{n-1}=$
 
 \noindent$1_{v_n}\otimes\sigma_{n-1}\circ\sigma(v_1\otimes\dots\otimes v_{n-1}, v_n)$,
      \item[{\rm c)}]  $\sigma_n^2=1$.
    \end{itemize}

      }\medskip

\noindent{\sl Proof} a) By induction on $n$,
$$\sigma_n=\Sigma_{n-1,1}\Delta^{(n-1)}\otimes 1(\overline{R})\Sigma_{n-1}\overline{R}_{n-1}=\Sigma_{n-1,1}\Sigma_{n-1}{ \Delta^{\text{op}}}^{(n-1)}\otimes 1(\overline{R})\overline{R}_{n-1}\otimes 1=$$
$$\Sigma_n{ \Delta^{\text{op}}}^{(n-1)}\otimes 1(\overline{R})\overline{R}_{n-1}\otimes 1=\Sigma_n\overline{R}_n,$$
where $\Sigma_{m,n}$ is the permutation on a tensor product of $m+n$ factors  that exchanges the first $m$ factors with the remaining $n$, and $\Sigma_n$ is the permutation that reverses the order in a tensor product space with $n$ factors.
b) One similarly shows that the right hand side of b) equals $\Sigma_n\overline{R}_n$. c) follows from b) and induction on $n$.

\bigskip

\section{Kirillov-Wenzl theory}

\noindent{\it a) Hermitian structures on $\text{Rep}({\mathcal U}_{\mathcal A})$}\bigskip

\noindent In this subsection we recall Kirillov's $^*$--involution in $\text{Rep}({\mathcal U}_{\mathcal A})$ making it into a tensor $^*$--category, following Wenzl's approach.

By an hermitian space we shall mean 
  a free finitely generated module $H$ over the involutive algebra ${\mathcal A}$ endowed with a   non degenerate, sesquilinear  ${\mathcal A}$-valued form $(\xi,\eta)$  on $H$,   hermitian with respect to the involution of ${\mathcal A}$:
  $$(\xi,\eta)^*=(\eta,\xi).$$
Given a linear map $T: H\to H'$ between   hermitian spaces, one can   define the adjoint $T^*:H'\to H$. 
The category ${\mathcal H}$ of  hermitian spaces   over   ${\mathcal A}$  is  a $^*$--category,
in the sense of \cite{DR1}.

By a $^*$--representation of ${\mathcal U}_{\mathcal A}$ we mean a  representation $u$ of it on an hermitian space which preserves the involutions, $u(x)^*=u(x^*)$.
  Weyl modules can be made  into     $^*$--representations as follows.
Given  a dominant weight $\lambda$, and the associated Weyl module $V_\lambda(x)$ at the generic level,   we can form both     the dual module
$V_\lambda(x)^*$ and the conjugate module $\overline{V_\lambda(x)}$.
 In the  former case, $U_x({\mathfrak g})$ acts   by transposition of the right action
$\xi a:=S^{-1}(a)\xi$  while  in the latter  the   action is given  on the conjugate vector space  by $a\overline{\xi}:=\overline{S^{-1}(a^*)\xi}$. These representations are still irreducible and have the same lowest weight, $-\lambda$.
Hence
a unique invertible intertwiner $\Phi: \overline{W_\lambda}\to W_\lambda^*$ arises from the identification  of their   natural lowest weight vectors, thereby defining the  form:
$$(\xi,\eta):= \Phi(\overline{\xi})(\eta).$$
This form turns out hermitian and, when restricted to   $V_\lambda({\mathcal A})$, takes values in the corresponding base   ring, hence it does make  that module into a $^*$--representation , see Sect. 2 in \cite{Wenzl}.

The main crux here    is that if  $u$ and $v$ are  $^*$--representations of ${\mathcal U}_{\mathcal A}$ then   their tensor product representation $u\otimes v$ defined in the usual way by means of the coproduct, 
$u\otimes v:=u\otimes v\circ\Delta$,
is {\it not} a $^*$--representation with respect to the natural  product form, $$(\xi\otimes\eta,\xi'\otimes\eta')_p:=(\xi,\xi')(\eta,\eta').$$
This is due to the fact that   the $\Delta$ is not $^*$--preserving. However, the ribbon structure fixes the problem. Not only this, but it is the first key step towards the construction of   unitary braid group representations,  in the sense that unitarity is first achieved at the    algebraic level,  as we next recall.  

\medskip

\noindent{\bf 4.1.  Proposition} {\sl For any pair of $^*$--representations $u$ and $v$ of ${\mathcal U}_{\mathcal A}$,   the following  form:
$$(\xi\otimes\eta, \xi'\otimes\eta'):=(\xi\otimes\eta, u\otimes v(\overline{R})\xi'\otimes\eta')_p\eqno(4.1)$$
is hermitian and makes
makes $u\otimes v$  into a $^*$--representation.   }\medskip

\noindent{\sl Proof} Recall that $\overline{R}$   is an invertible and selfadjoint    element of the completed tensor product of the extended algebra by itself, and this  implies  that $u\otimes v(\overline{R})$ is a selfadjoint operator with respect to the product form. Hence the right hand side does define an hermitian form on the tensor product space. Moreover, the adjoints $A^*$ and $A^{+}$  of an operator $A: U\otimes V\to U\otimes V$, 
with $U$ and $V$ the spaces corresponding to $u$ and $v$, and computed
 with respect to the modified  and  product  forms respectively, are related by
 $$A^*=u\otimes v(\overline{R})^{-1}A^+{u}\otimes { v}(\overline{R}).$$
 It follows that for $a\in{\mathcal U}_{\mathcal A}$,
 $$u\otimes v(a)^*=u\otimes v(\overline{R})^{-1}(u\otimes v\circ\Delta(a))^+u\otimes v(\overline{R})=$$
 $$u\otimes v(\overline{R})^{-1}u\otimes v\circ \Delta^{\text{op}}(a^*)u\otimes v(\overline{R})=u\otimes v(a^*).$$
 
 \bigskip

 In order to simplify notation, we shall switch   from   representation   to   module notation when this will cause no confusion. Hence we identify a representation $u$ with its space $U$, and simply denote by $a\xi$ the action of the operator $u(a)$ on a vector $\xi$.

 We now make all objects of   $\text{Rep}({\mathcal U}_{\mathcal A})$ into $^*$--representations 
 in the   way just explained. In other words, we endow a tensor product module  $V_{\lambda_1}\otimes\dots\otimes V_{\lambda_n}$ 
 with the   form  
 defined by the action of the matrix $\overline{R}^{(n)}$,
 $$(\xi,\eta):=(\xi, \overline{R}^{(n)}\eta)_p,$$  which is indeed hermitian by selfadjointness of $\overline{R}^{(n)}$, cf. Proposition 3.6.
  We can thus compute the adjoint of every arrow of $\text{Rep}({\mathcal U}_{\mathcal A})$.\medskip

 \noindent{\bf  4.2. Theorem} {\sl The hermitian forms so defined on objects of $\text{Rep}({\mathcal U}_{\mathcal A})$ make    it  into a ribbon   tensor  $^*$--category. Furthermore, both the braiding  operators $\varepsilon({U,V})=\Sigma R\in(U\otimes V, V\otimes U)$
 and the coboundary operators $\sigma(U, V)=\Sigma\overline{R}\in(U\otimes V, V\otimes U)$ are unitary arrows of $\text{Rep}({\mathcal U}_{\mathcal A})$.
 
 }\medskip
 
 \noindent{\sl Proof} We need to verify that the tensor  product makes the set of objects into a (unital) associative semigroup and that the tensor product of arrows satisfies $(S\otimes T)^*=S^*\otimes T^*$.  In regards to the former, a computation shows that the hermitian form of $(U\otimes V)\otimes W $   is defined by the action of  $\overline{R}\otimes1\Delta\otimes 1(\overline{R})$, while that for $U\otimes (V\otimes W) $    corresponds to  $1\otimes\overline{R} 1\otimes\Delta(\overline{R})$, and these two coincide by Proposition 3.5. Let now $S\in (U, U')$, $T\in(V, V')$. Using   $(4.2)$ we get
  $$(S\otimes T)^*=\Theta^{-1} R^{-1}S^*\otimes T^* R\Theta=S^*\otimes T^*.$$
  We next check unitarity of $\varepsilon({U,V})$,
 $$\varepsilon({U,V})^*=
 \overline{R}^{-1}(\Sigma R)^+ \overline{R}=\Theta^{-1} {R}^{-1} R^+\Sigma {R} \Theta=$$
$$ \Theta^{-1}{R}^{-1} R_{21}^{-1}\Sigma{R}\Theta=\Theta^{-1} {R}^{-1}\Sigma\Theta= (R\Theta)^{-1}\Theta_{21}\Sigma=$$
$$(\Theta_{21}R)^{-1}\Theta_{21}\Sigma=\varepsilon({U,V})^{-1}$$
by Lemma 3.2 b).
Finally, since $\sigma(U, V)=\varepsilon(U, V)\Theta$, it suffices to show that $\Theta\in(U\otimes V, U\otimes V)$ is unitary.
By Lemma 3.2 again,
$$\Theta^*=\overline{R}^{-1}(\Theta)^+\overline{R}=\overline{R}^{-1}(\Theta)_{21}^{-1}\overline{R}=$$
$$\Theta^{-1}(\Theta_{21}R)^{-1}R\Theta=\Theta^{-1}(R\Theta)^{-1}R\Theta=\Theta^{-1}.$$
 \bigskip

\noindent{\it  b) $C^*$--structures at roots of unity}\bigskip

\noindent We next recall the main results of \cite{Wenzl}.
The  hermitian form previously considered on  Weyl modules $V_\lambda({\mathcal A})$     specialises to any $q\in{\mathbb T}$,
hence induces  a complex   sesquilinear and  hermitian  form  on $V_\lambda(q)$.

  \medskip

\noindent{\bf 4.3. Theorem} {\sl   Let $q\in{\mathbb T}$, then the hermitian form of $V_\lambda({\mathcal A})$ specialises to a  positive definite form on $V_\lambda(q)$ in the following cases,
 \begin{itemize}
   \item[{\rm a)}] for   $\lambda\in\Lambda^+$ satisfying $\langle \lambda+\rho, \theta\rangle< d+\frac{1}{|t|}$  if 
   $q=e^{i\pi t}$ is not a root of unity,
  \item[{\rm b)}] for   $ \lambda\in\overline{\Lambda}_\ell$ if  $q=e^{{\pi i}/{d\ell}}$.
   \end{itemize}
}\medskip

Let $V=V_\kappa$ denote the fundamental representation of ${\mathfrak g}$ as in \cite{Wenzl}, with $\kappa$ its dominant weight.
 For example, if ${\mathfrak g}$ is of type $A$ or $C$, $V$ is the vector representation,     for types $B$ is the spin representation, for type $D$ is one of the two   spin representations,   $V_{\kappa_1}$, $V_{\kappa_2}$.  The order  of the root of unity is chosen in such a way that $V(q)$  is non-negligible,  i.e., $\kappa\in\Lambda_\ell$.
  This   representation  ($V_{\kappa_1}\oplus V_{\kappa_2}$ in type $D$)
      has the property  that   its    tensor powers contain  every irreducible of ${\mathfrak g}$.

     We summarise  properties  concerning  fusion of tensor products with the fundamental.      We add a few remarks in the proof aiming  to connect  various  results known in the literature.  
   \medskip
 
 \noindent{\bf 4.4. Theorem} {\sl    Let $V$ be the fundamental representation of ${\mathfrak g}\neq E_8$,  ($V=V_{\kappa_i}$ in the type $D$ case)
 and pick $\lambda\in\Lambda_\ell$. Then 
   \begin{itemize}
 \item[{\rm a)}]  all irreducible submodules $V_\mu$ of $V_\lambda\otimes V$ in   $\text{Rep}({\mathfrak g})$ have weights $\mu\in\overline{\Lambda_\ell}$,
  \item[{\rm b)}]  the maximal negligible and non-negligible summands $V_\lambda(q)\otimes V(q)$ in ${\mathcal T}_\ell$ are unique and given by 
  $$N_\lambda=\bigoplus_{\mu\in\overline{\Lambda_\ell}\setminus\Lambda_\ell} m_\mu V_\mu,\quad\quad  V_\lambda\underline{\otimes} V=\bigoplus_{\mu\in\Lambda_\ell} m_\mu V_\mu,$$
  with multiplicities as in the classical case. Specifically,
   both decompositions  are
   multiplicity free for  ${\mathfrak g}\neq F_4$;  while $N_\lambda$ is so for   ${\mathfrak g}= F_4$,
   \item[{\rm c)}] if ${\mathfrak g}$ is a classical Lie algebra, $V$ is minuscule hence fusion rules are given by
    $$m_\mu=1 \iff \mu\in  E_\lambda:=\{\lambda+\gamma: \gamma \text{ weight of } V\}\cap\Lambda^+.$$
     \end{itemize}
 }\medskip
 
 \noindent{\sl Proof}  a)  and b) are Theorem 3.5. in \cite{Wenzl}.  Notice that
   on one  hand,
the dominant weights of the Weyl modules  appearing as factors in the Weyl filtration of an indecomposable submodule $T_\mu$ with maximal weight $\mu$ of $V_\lambda(q)\otimes V(q)$  must belong to the set of dominant weights arising from the classical decomposition of $V_\lambda\otimes V$ see e.g. Prop. 3 in \cite{Sawin}. On the other, the Weyl modules $V_\mu(q)$ with $\mu\in\overline{\Lambda_\ell}$ are irreducible and tilting  
\cite{Andersen,  AP, Wenzl}. Since decomposition into isotypic components of irreducibles  is unique, this shows   b).
 c) Every summand $V_\mu$ of $V_\lambda\otimes V$ has weight   of the form $\mu=\lambda+\gamma$,
  where $\gamma$ is a weight of $V$. But if in addition   ${\mathfrak g}$ is of type $ABCD$,  then $V$ is minuscule by, e.g., table A. 2.3 in \cite{LR}. This implies the previous statement has a converse:    for any weight $\gamma$ of $V$ such that $\lambda+\gamma$ is dominant, $V_{\lambda+\gamma}$ does appear  in $V_\lambda\otimes V$, see, e.g., Lemma 3.1 in \cite{LZ} .  \medskip

\noindent{\it Notation.} From now on we shall mostly work    in the tilting category ${\mathcal T}_\ell$, or its quotient ${\mathcal F}_\ell$, hence we shall simply write $V_\lambda$ for the Weyl module $V_\lambda(q)$.
 
 The specialised coboundary matrices $\overline{R}^{(n)}$ are still   invertible and selfadjoint, hence the corresponding  hermitian form of  $V_{\lambda_1}\otimes\dots\otimes V_{\lambda_r}$ is non degenerate, provided  $\lambda_i\in\overline{\Lambda}_\ell$, by Theorem 4.3.
However, it   may  degenerate on subspaces. Wenzl used the above properties  of the fundamental representation to show the following important result.

\medskip

\noindent{\bf 4.5. Theorem } {\sl Under the same assumptions of the previous theorem,  for $q=e^{{\pi i}/{d\ell}}$ and $\lambda\in\Lambda_\ell$, Wenzl's hermitian form is positive definite on the submodule  $V_\lambda\underline{\otimes} V$. Furthermore, for  any
  $\gamma\in\Lambda_\ell$,   the canonical projection
  $$p_{\lambda,\gamma}: V_\lambda\otimes V\to m_\gamma V_\gamma$$
   is selfadjoint under the same form  and one  has
   $p_{\lambda,\gamma}p_{\lambda,\mu}=0$ for $\gamma\neq\mu$.
   }\medskip
   
   \noindent{\bf  4.6. Remark } Notice that Wenzl includes the $E_8$ case, which we exclude   in this paper,  as, beyond being the most delicate case, it may give rise to indecomposable but reducible  summands of $V_\lambda\otimes V$, cf. Theorem 3.5 c) of \cite{Wenzl}, while our later constructions need complete reducibility  of such tensor products.\medskip

\section{Wenzl's functor}

In this section we construct   a sequence of projections $p_n$ on powers $V^{\otimes n}$ onto Hilbert subspaces and describe their main properties.

Consider, for each $\lambda\in\Lambda_\ell$, the (selfadjoint) projection
$p_\lambda: V_\lambda\otimes V\to V_\lambda\underline{\otimes} V,$  
  given by $p_\lambda=\sum p_{\lambda,\gamma}$,
where the sum is made over all possible summands $V_\gamma$ of $V_\lambda\otimes V$
with $\gamma\in\Lambda_\ell$.
  Let $p_0$ and $p_1$  denote the identity maps on the trivial module ${\mathbb C}$, and $V$ respectively.  
Set, for ${\mathfrak g}\neq D_n$, 
$p_2=p_\kappa$, with $\kappa$   the dominant weight of $V$. In the type $D_n$ case, if  $\kappa_1$ and $\kappa_2$ denote the dominant weights of the two half-spin irreducible  subrepresentations of $V$, we     set $p_2=p_{\kappa_1}+p_{\kappa_2}$.
Define inductively canonical   projectons $p_n$ on  $V^{\otimes n}$, as follows. Given    $$p_{n}: V^{\otimes n}\to\oplus_{\alpha\in A_n}V_{\alpha, n}$$   iteratively   projecting onto the non-negligible part of a canonical orthogonal decomposition into submodules,
$$V^{\otimes n}=\bigoplus_{\alpha\in A_n} V_{\alpha, n}\,\bigoplus\, N_{n},$$ with $N_n$ negligible,  consider  unitaries $U_{\alpha,\mu}: V_\mu\to 
V_{\alpha, n}
$ with $\mu\in\Lambda_\ell$,  and
 then   set $$p_{n+1}=\sum_{\alpha\in A_{n}}            
 U_{\alpha,\mu}\otimes 1_V\circ p_\mu\circ  U_{\alpha,\mu}^{-1}\otimes 1_V\circ p_n\otimes 1_V.$$
   The next   lemmas   will often  turn out useful.\medskip 
    
  \noindent{\bf 5.1. Lemma} {\sl 
   \begin{itemize}
   \item[{\rm a)}] $p_n\circ p_m\otimes 1_{V^{\otimes r}}=p_n= p_m\otimes 1_{V^{\otimes r}}\circ p_n,$
    \item[{\rm b)}]  
 $A\otimes 1_{V^{\otimes r}}\circ p_{m+r}=p_{n+r}\circ A\otimes 1_{V^{\otimes r}}=A\underline{\otimes} 1_{V^{\underline{\otimes} r}}$,\quad  \text{for} \quad  $A\in(V^{\underline{\otimes} m}, V^{\underline{\otimes} n})$ . 
     \end{itemize} }\medskip
     
 \noindent{\sl Proof} The case $r=1$ holds by construction,  a simple iteration shows the general case.  b) It suffices to show the first equality. Assume $r=1$. For $A=U_{\beta,\mu} U_{\alpha,\mu}^{-1}: V_{\alpha, m} \to V_\mu\to V_{\beta, n}$
 $$A\otimes 1_V\circ p_{m+1}=U_{\beta,\mu}\otimes 1_V\circ p_\mu\circ U_{\alpha,\mu}^{-1}\otimes 1_V=p_{n+1}\circ A\otimes 1_V=A\underline{\otimes} 1_V.$$
The equality also holds  for $A\in(V^{\underline{\otimes} m}, V^{\underline{\otimes} n})$, by complete reducibility of $V^{\underline{\otimes} m}$
and $V^{\underline{\otimes} n}$.
For $r>1$, by a) and induction on $r$,
$$A\otimes 1_{V^{{\otimes}r}}\circ p_{m+r}= (A\otimes 1_{V^{{\otimes}r-1}}\circ p_{m+r-1})\otimes 1_V\circ p_{m+r}=$$
$$(p_{n+r-1}\circ A\otimes 1_{V^{{\otimes}r-1}})\otimes 1_V\circ p_{m+r}=p_{n+r}\circ(p_{n+r-1}\circ A\otimes 1_{V^{{\otimes}r-1}})\otimes 1_V=$$
$$p_{n+r}\circ A\otimes 1_{V^{\otimes r}}.$$
\medskip

One inductively derives from a)
 that $p_n$ is selfadjoint w.r.t. Kirillov-Wenzl form, for all $n$.
In general, the projections $p_n$ do not satisfy right associativity, i.e. $p_n\circ 1_{V^{\otimes r}}\otimes p_m\neq p_n$.
 We see this with an example. \medskip

\noindent{\bf 5.2. Example} Set ${\mathfrak g}={\mathfrak sl}_2$, $\ell=3$. 
We denote by $\alpha$  the  unique  simple root, so  $d=1$,  $\theta=\alpha$, $\rho=\frac{1}{2}\alpha$, and $\left\langle\alpha,\alpha\right\rangle=2$.  These data imply that the open Weyl alcove contains   two non-negligible irreducible representations $V_0$ and $V_1$, with $V_0$ the trivial  and $V_1=V$   the vector representation, respectively. Furthermore, $V_2$ is a    negligible irreducible representation. We have the following decompositions    $V\otimes V\simeq V_0\oplus V_2$, implying    $V\underline{\otimes}V\simeq V_0$, and $V_0\otimes V\simeq V\simeq V_0\underline{\otimes}V$. The quantum determinant element
$S=\psi_1\otimes\psi_2-q\psi_2\otimes \psi_1$ is a morphism in $(V_0,V\otimes V)$ 
satisfying $S^*\otimes 1_V\circ 1_V\otimes S=-1_V$, $S^*S=[2]_q=q+q^{-1}$, see Prop. A.4 and A.5 in the appendix for details.
 Now $q+q^{-1}=1$ for $\ell=3$, hence $S$ is an isometry, so  we have $p_2=SS^*$.  There is no further truncation at the next power, so  $p_3=p_2\otimes 1_V=SS^*\otimes 1_V.$ We next compute 
$$p_3\circ 1_{V}\otimes p_2=p_2\otimes 1_{V}\circ 1_{V}\otimes p_2=$$
$$SS^*\otimes 1_{V}\circ 1_{V}\otimes SS^*= - S\otimes 1_{V}\circ 1_{V}\otimes S^*$$
and it is now easy to see that it differs from $p_3$.
\medskip

However, the following property will function as a replacement of right associativity failure of the $p_n$. It is the way how   Andersen-Gelfand-Kazhdan properties of Sect. 2   will often manifest themselves in this paper.
\medskip

\noindent{\bf 5.3. Lemma} {\sl For any pair of morphisms $S\in(V^{\otimes m}, V^{\otimes n})$, $T\in(V^{\otimes n}, V^{\otimes r})$ of the tilting category,
$$p_r\circ T\circ 1_{V^{\otimes s}}\otimes p_t\otimes 1_{V^{\otimes u}}\circ S\circ p_m=p_r\circ T\circ   S\circ p_m$$}\medskip

\noindent{\sl Proof} The range of $1_{V^{\otimes s}}\otimes (1-p_t)\otimes 1_{V^{\otimes u}}$
is a negligible module by $(3)$, while   $p_m$ and $p_r$ have  support on non-negligible modules, hence 
$p_r\circ T\circ 1_{V^{\otimes s}}\otimes (1-p_t)\otimes 1_{V^{\otimes u}}\circ S\circ p_m=0$
by $(2)$.\medskip

 As  an illustration,  in our   example 5.2,
$$p_3\circ 1\otimes p_2\circ p_3= - S\otimes 1_{V}\circ 1_{V}\otimes S^*\circ SS^*\otimes 1_V=$$
$$ - S\otimes 1_{V}\circ (S^*\otimes 1_V\circ 1_V\otimes S)^*\circ S^*\otimes 1_V=SS^*\otimes 1_V=p_3.$$\medskip

  We  next  describe the first important consequence. Consider the category  ${\mathcal G}_\ell$   with objects $$V^{\underline{{\otimes}} n}:=p_nV^{{\otimes} n}$$
and morphisms
 $$(V^{\underline{{\otimes}} n}, V^{\underline{{\otimes}} m})=\{S\in (V^{{{\otimes}} n}, V^{{{\otimes}} m}): Sp_n=p_mS=S\}.$$
We introduce a   tensor product in ${\mathcal G}_\ell$,
 $$V^{\underline{{\otimes}} m}\underline{\otimes}V^{\underline{{\otimes}} n}:=V^{\underline{{\otimes}} m+n},\quad\quad  S\underline{\otimes} T:=p_{m'+n'}\circ S\otimes T\circ p_{m+n},$$
 where $S\in (V^{\underline{{\otimes}} m}, V^{\underline{{\otimes}} m'})$, $T\in (V^{\underline{{\otimes}} n}, V^{\underline{{\otimes}} n'})$.
  \bigskip
 
The following theorem is  due to  \cite{Wenzl}.  We briefly comment on  those aspects that will play a role in the sequel.  \bigskip
 
 \noindent{\bf  5.4. Theorem} {\sl This tensor product makes ${\mathcal G}_\ell$ into a strict tensor $C^*$--category with a unitary braided symmetry given by 
 $$\underline{\varepsilon}(b):= p_n\varepsilon(b)p_n, \quad b\in{\mathbb B}_n.$$ Furthermore the composition of
 the inclusion ${\mathcal G}_\ell\to{\mathcal T}_\ell$ with the quotient   ${\mathcal T}_\ell\to{\mathcal F}_\ell$ is an equivalence of braided tensor $^*$--categories.}\medskip
 
 \noindent{\sl Proof} 
 The statement about the equivalence of the two categories can be proved in a way similar to the proof of Lemma 1.1 in \cite{GK},  once the tensor structure 
 of ${\mathcal G}_\ell$ has been verified.  The identity arrow on $V^{\underline{\otimes} n}$ is $p_n$.  We need to verify the following properties: $p_n\underline{\otimes } p_m=p_{n+m}$; $(S\underline{\otimes }T)^*=S^*\underline{\otimes} T^*$, associativity, $(S\underline{\otimes} T)\underline{\otimes} U= S\underline{\otimes}( T\underline{\otimes} U)$, and  the exchange rule  between composition and tensor product, $(S\circ T)\underline{\otimes}(S'\circ T')=(S\underline{\otimes} S')\circ(T\underline{\otimes} T')$. The second property follows from the fact that the projections $p_n$ are selfadjoint arrows of ${\mathcal T}_\ell$, while  the remaining properties follow   from  Lemma 5.3. For example,
  $$(S\underline{\otimes} T)\underline{\otimes} U=p_{m'+n'+p'}\circ p_{m'+n'}\otimes 1_{p'}\circ S\otimes T\otimes U\circ p_{m+n}\otimes 1_{p}\circ p_{m+n+p}=$$
 $$p_{m'+n'+p'}\circ S\otimes T\otimes U\circ  p_{m+n+p}=p_{m'+n'+p'}\circ 1_{m'}\otimes p_{n'+m'}\circ S\otimes T\otimes U\circ1_{m}\otimes p_{n+m}\circ   p_{m+n+p}=$$
 $$S\underline{\otimes} (T\underline{\otimes} U).$$
 Similar computations show that $\underline{\varepsilon}$ is a braided symmetry for ${\mathcal G}_\ell$.  The $C^*$--structure and unitarity of the braided symmetry are consequences of the results recalled in the previous section.

   \bigskip

    This realisation  of the quotient category is useful to obtain an analogue of a fibre  functor, the starting point of our reconstruction, 
  $W: {\mathcal G}_\ell\to{\mathcal H}$ taking $V^{\underline{\otimes} n}$ to its Hilbert space, $p_nV^{{\otimes} n}$, and acting identically on arrows.
In the following we shall identify the abstract quotient category ${\mathcal F}_\ell$ with its concrete realisation ${\mathcal G}_\ell$ and we shall not distinguish the notation.

\medskip

    \section{Algebraic compact quantum groups ${\mathcal C}(G_q)$ for $q\in{\mathbb T}$   generic}
    
    In this section  we pause on roots of unity and Wenzl's projections. Thus we fix $q\in{\mathbb T}$  not a root of unity, and ${\mathfrak g}$ is a complex simple Lie algebra. We apply Tannaka-Krein duality to  the tensor category of $^*$--representations of 
    $U_q({\mathfrak g})$ and construct   a    quantum group, denoted ${\mathcal C}={\mathcal C}(G_q)$, which may be viewed as a quantization of     the function algebra over the compact group $G$ corresponding to the  compact  real form  of  ${\mathfrak g}$.  This construction        conveys a general strategy that we will extend at roots of unity in later sections.

The main novelty of our construction, in comparison with the corresponding construction for the compact quantum groups $G_q$ for $q\in{\mathbb R}$   \cite{W},
 is the algebra involution, which now relies on 
 the coboundary.
        The main result of this section is cosemisimplicity of   the associated function algebra,
    which is derived from existence of a Haar functional, which, however, 
    is not positive.  
            
Let ${\mathcal V}_{\overline{R}}$ denote the category  with objects the hermitian spaces arising as representation spaces of the objects of ${\mathcal T}_\ell$
and arrows all the linear maps. The adjoint still makes sense for the arrows of this category, and one easily verifies that in this way
${\mathcal V}_{\overline{R}}$ satisfies all the axioms of a tensor $^*$--category, except the rule of the adjoint of a tensor product arrow, which, unless $S$ and $T$ are morphisms of representations, it is    replaced by 
the   property
$$(S\otimes T)^*=\overline{R}^{-1}\circ S^*\otimes T^*\circ \overline{R}.\eqno(6.1)$$
We consider  the spaces of the tensor powers of the fundamental representation of $U_q({\mathfrak g})$,   simply denoted       by $V^n$,    $n\in{\mathbb N}$. Same notation applies to simple tensors $\phi_1\dots\phi_n$ of  $V^n$. Elements of $V^n$ are regarded as arrows of $(\iota, V^n)$ in ${\mathcal V}_{\overline{R}}$  in the natural way. Similarly,  elements of the dual space $(V^n)^*$ will be regarded as arrows of $(V^n,\iota)$.

 We form the diagonal subalgebra ${\mathcal D}={\mathcal D}(V)$ of the mixed tensor algebra associated to $V$,
  $${\mathcal D}:=\oplus_{n\in{\mathbb N}}( {V^n}){^*}\otimes V^n.$$
  In other words, the multiplication in ${\mathcal D}$ is given by
  $$(\varphi\otimes\psi)(\varphi'\otimes\psi'):=\varphi\varphi'\otimes\psi\psi',\eqno(6.2)$$
  $\varphi\in ({V^n}){^*}$, $\varphi'\in ({V^m}){^*}$,  $\psi\in V^n$, $\psi' \in V^m$. Notice that the subset of simple tensors  $$\{v_{\phi,\psi}:=\phi\otimes\psi;\  \psi\in V, \phi\in V^*\},$$ generates ${\mathcal D}$ as an algebra.
   A non-involutive quantum group can be reconstructed via usual Tannaka-Krein duality.
 Specifically, let ${\mathcal I}$ be the linear subspace of ${\mathcal D}$ generated by 
$${\mathcal I}:=\text{l.s.}\{\varphi\otimes A\circ\psi-\varphi\circ A\otimes\psi; \ \psi\in V^n, \varphi\in{V^*}{^m}, A\in(V^n, V^m)\}.$$

 \medskip

\noindent{\bf 6.1. Proposition} {\sl ${\mathcal I}$ is a two-sided ideal of ${\mathcal D}$, hence
the quotient space
$${\mathcal C}:={\mathcal D}/{\mathcal I}$$
is an associative and unital  algebra.}\medskip

\noindent{\sl Proof}  It suffices to show that   ${\mathcal I}$ is stable under left and right multiplication by the generators $v_{\phi,\psi}$. This follows from the fact that   the arrows $A\in(V^n, V^m)$  
are stable under left and right tensoring  by $1_V$.
 
\medskip
 
 We next make ${\mathcal C}$ into an involutive  Hopf algebra.  The adjoint can equivalently be     introduced 
 in two different ways, one relying on   the $n$-th tensor power of  Hermitian structure of $V$
 and the other on
 Kirillov-Wenzl
    hermitian structure of $V^n$ (i.e. a deformation of the former via $\overline{R}_n$). Both turn out useful.  
  
  For $\psi\in V$, $\psi^*$ is the linear functional 
on $V$ defined by $\psi^*(\psi')=(\psi,\psi')$.
 We can thus identify
the dual  space $({V^n})^*$ with ${V^*}{^n}$, denoting by
 $\psi_1^*\dots\psi_n^*$ the tensor product functional $\psi_1^*\otimes\dots\otimes\psi_n^*$. 
We introduce the $^*$--involution as follows
  $$(\phi_1^*\dots\phi_n^*\otimes\psi_1\dots\psi_n)^*=\psi_n^*\dots\psi_1^*\otimes\phi_n\dots\phi_1,\eqno(6.3)$$
 for $\phi_i$, $\psi_i\in V$.
 This is clearly a well defined map on ${\mathcal D}$. We next  express this adjoint  
    in terms of the hermitian structure
 of $V^n$ and the (specialised) coboundary operators $\sigma_n=\Sigma_n\overline{R}^{(n)} \in(V^{ n}, V^{ n})$   of Sect. 3, \medskip
 
 \noindent{\bf 6.2. Lemma} {\sl    For $\psi_1,\dots\psi_n\in V$,
 $$(\psi_1\dots\psi_n)^*= \psi_n^*\dots\psi_1^*\sigma_n.$$

   }\medskip

 \noindent{\sl Proof}
 Taking into account    $(6.1)$ for arrows $(\psi_1\dots\psi_n)^*$ and 
 $\psi_1^*\dots\psi_n^*$,
 $$(\psi_1\dots\psi_n)^*=\psi_1^*\dots\psi_n^*\circ\overline{R}^{(n)}=\psi_n^*\dots\psi_1^*\sigma_n.$$ \medskip

 We can thus  alternatively represent simple tensors   of ${\mathcal D}$ in the form $\psi^*\sigma_n^{-1}\otimes\phi$,    and define the  adjoint of ${\mathcal D}$ by  
$$(\psi^*\sigma_n^{-1}\otimes\phi)^*=\phi^*\sigma_n^{-1}\otimes\psi,\quad\quad \phi,\psi\in V^n.\eqno(6.4)$$\medskip

 \noindent{\bf  6.3. Proposition} {\sl This involution makes ${\mathcal D}$ into a $^*$--algebra and ${\mathcal I}$ into a $^*$--ideal. Hence ${\mathcal C}$ is a $^*$--algebra.}\medskip
 
 \noindent{\sl Proof}   If $A\in (V^m, V^n)$, $\psi\in V^n$, $\phi\in V^m$,
 we compute the following adjoints by means of $(6.4)$, $$(\psi^*\sigma_n^{-1}A\otimes\phi)^*=((\sigma_n^{-1}A^*\sigma_n\psi)^*\sigma_n^{-1}\otimes\phi)^*=\phi^*\sigma_n^{-1}\otimes\sigma_n^{-1}A^*\sigma_n\psi,$$
since $\sigma_n$ is unitary. Similarly,
$$(\psi^*\sigma_n^{-1}\otimes A\phi )^*=(A\phi)^*\sigma_n^{-1}\otimes\psi=\phi^*\sigma_n^{-1}(\sigma_n A^*\sigma_n^{-1})\otimes\psi.$$
Taking into account the important relation $\sigma_n^2=1$,  this   computation shows that
 the adjoint of an element of the form 
$ \psi^*\sigma_n^{-1}A\otimes\phi-\psi^*\sigma_n^{-1}\otimes A\phi$ is of the same form. \medskip

\noindent{\bf 6.4. Remark} Notice that, for class elements of ${\mathcal C}$, for $\psi,\psi'\in V^n$,
$[\psi^*\otimes\psi']^*=[\psi'^*\otimes\psi]$, where the adjoints of $\psi$ and $\psi'$ refer to the hermitian form of $V^n$. 
On the other hand, up to a scalar multiple, this form restricts to the hermitian form of $V_\alpha$ if $V_\alpha$ is 
a submodule of $V^n$. Hence, for  $ \psi, \psi'\in V_\alpha$, we may regard $\psi^*$ and $\psi'^*$
as adjoints 
  relative to the hermitian form of $V_\alpha$. \medskip

We next  introduce in ${\mathcal C}$ the structure of  a Hopf algebra by first   endowing ${\mathcal D}$ with   coproduct defined by:
$$\Delta(\phi\otimes\psi)=(\phi\otimes\eta_r)\otimes(\eta^r\otimes\psi),$$
where $\eta_r\in V^n$ is a basis and $\eta^r\in(V^n)^*$ is a  dual basis. A routine computation shows that $\Delta$ does not depend 
on the choice of the basis.

\medskip

\noindent{\bf 6.5. Proposition} {\sl  The coproduct $\Delta$ is   unital and coassociative, and satisfies
$\Delta(a^*)= \Delta^{\text{op}}(a)^*$. Furthermore ${\mathcal I}$ is also a coideal, hence $\Delta$ induces a coproduct on ${\mathcal C}$, still denoted $\Delta$, satisfying the same properties.}\medskip

\noindent{\sl Proof} It is straightforward to check that $\Delta$ is   unital, multiplicative and coassociative. Since
$a\to\Delta(a^*)$ and $a\to \Delta^{\text{op}}(a)^*$ are both antimultiplicative maps on ${\mathcal D}$, it suffices, and it is easy, to check that $\Delta(a^*)= \Delta^{\text{op}}(a)^*$ for $a=\phi^*\otimes\psi$, $\phi$, $\psi\in V$.
  Finally, if $A\in(V^n, V^m)$ is an intertwiner, and $\eta_i\in V^n$, $\xi_j\in V^m$ are linear bases with dual bases $\eta^i$, $\xi^j$ respectively, for $\psi\in V^n$, $\phi\in(V^m)^*$,
$$\Delta(\phi\circ A\otimes\psi)=\phi\circ A\otimes \eta_i\otimes\eta^i\otimes\psi,\quad\Delta(\phi\otimes A\circ\psi)=\phi\otimes\xi_j\otimes\xi^j\otimes A\circ\psi.$$
An easy computation gives
$$\Delta(\phi A\otimes\psi-\phi\otimes A\psi)=(\phi A\otimes \eta_i-\phi\otimes A\eta_i)\otimes\eta^i\otimes\psi
-\phi\otimes\xi_j\otimes(\xi^j A\otimes\psi-\xi^j\otimes A\psi),
$$
implying $\Delta({\mathcal I})\subset {\mathcal I}\otimes {\mathcal D}+{\mathcal D}\otimes{\mathcal I}$. 
\medskip

We next introduce the functional  
$$h:{\mathcal C}\to{\mathbb C},$$  which corresponding  to the Haar measure in the classical case by,
$$h(\phi\otimes\psi)= \phi(e_n\psi), \quad \psi\in V^n, \phi\in {V^n}^*,$$
where $e_n$ is the projection onto the   isotypical component of the trivial subrepresentation of $V^n$. Every $A\in(V^n, V^m)$ satisfies $Ae_n=e_mA$, and this shows that $h$ annihilates ${\mathcal I}$, thus defining a linear functional on ${\mathcal C}$.

This functional turns out useful to show cosemisimplicity of   ${\mathcal C}$.
Consider a complete set, $V_\lambda$, $\lambda\in \Lambda^+$, of irreducible  representations.
An isometric intertwiner $S\in(V_\lambda, V^{\otimes n})$ induces a   linear inclusion
$$V_\lambda^*\otimes_{\mathbb C} V_\lambda\to{\mathcal C}$$
which takes a simple tensor $\phi\otimes\psi$ to the class of $\phi S^*\otimes S\psi$. 
\medskip

\noindent{\bf 6.6. Proposition} {\sl The above   inclusion $V_\lambda^*\otimes_{\mathbb C} V_\lambda\to{\mathcal C}$ does not depend on $n$ and   $S\in(V_\lambda, V^{\otimes n})$.
Furthermore, the image is a subcoalgebra.
}\medskip

\noindent{\sl Proof} If $T\in(V_\lambda, V^{\otimes m})$ is another isometry then in the quotient,
$$[ \phi S^*\otimes S\psi]=[\phi S^*\otimes ST^*T\psi]=[\phi S^*ST\otimes T\psi]=[\phi T^*\otimes T\psi].$$
For the last statement we notice that a similar computation shows that if $\psi$ and $\psi'$ lie in orthogonal invariant subspaces of some $V^{\otimes n}$ then $[\psi^*\otimes\psi']=0$.
Therefore if $\phi$ and $\psi$ lie in the same irreducible component of $V^{\otimes n}$ then
$\Delta([\phi^*\otimes \psi])$ can be simply expressed by means of an orthonormal basis of that submodule, rather than of the whole $V^{\otimes n}$.
  \medskip

 \noindent{\bf 6.7. Theorem} {\sl For every $\lambda\in \Lambda^+$, the natural inclusion 
$V_\lambda^*\otimes_{\mathbb C}V_\lambda \to {\mathcal C}$ is faithful.
Therefore, as a coalgebra,
$${\mathcal C}=\bigoplus_{\lambda\in \Lambda^+} V_\lambda^*\otimes_{\mathbb C}V_\lambda.$$
and ${\mathcal C}$ is cosemisimple.}\medskip
   
\noindent{\sl Proof} Consider isometries $S_{\lambda, i}\in(V_\lambda, V^n)$ such that $\sum_{\lambda, i}S_{\lambda, i}S_{\lambda, i}^*=1_{V^n}$. An element $\phi\otimes\psi\in(V^n)^*\otimes V^n$, regarded as an element of ${\mathcal C}$, can be written in the form  
$$\phi\otimes\psi=\phi\otimes \sum S_{\lambda,i}S_{\lambda, i}^*\psi=\sum\phi S_{\lambda,i}\otimes  S_{\lambda, i}^*\psi.$$
Hence ${\mathcal C}$ is linearly generated by the various $V_\lambda^*\otimes V_\lambda$. We next show that the linear inclusion of 
$V_\lambda^*\otimes_{\mathbb C} V_\lambda$ in ${\mathcal C}$ is faithful.   
We evaluate $h(ab)$, for $a=\psi^*\otimes\psi'\in V_\lambda^*\otimes V_\lambda$, $b=\phi^*\otimes\phi'\in V_\mu^*\otimes V_\mu$.
  If $V_\lambda$ is not conjugate to $V_\mu$, $h(ab)=0$, as $V_\lambda\otimes V_\mu$ does not contain the trivial representation. If $V_\lambda$ is conjugate to $V_\mu$, we identify $V_\lambda$ with the conjugate hermitian representation $\overline{V_\mu}$. If    $\psi_i$ is an orthonormal basis
  of $V_\mu$, then $r:=\sum_i\overline{\psi_i}\otimes K_{-2\rho}\psi_i\in(\iota, V_\lambda\otimes V_\mu)$. Lemma 3.4 of \cite{Wenzl} easily implies that
 the adjoint of $r$ with respect to the hermitian form is given by $r^*(\overline{\psi}\otimes\phi):=(\psi,\phi).$
 Hence for $a=\overline{\phi}^*\otimes\overline{\psi}$, $b=\xi^*\otimes\eta$, and a non-zero scalar $d(\mu)$, the quantum dimension of $V_\mu$,
 $$d(\mu) h(ab)= \overline{\phi}^*\xi^*\otimes rr^*\overline{\psi}\eta=\overline{\phi}^*\xi^*r\otimes r^*\overline{\psi}\eta=$$$$\sum_i(\psi_i,\phi)(\xi, K_{-2\rho}\psi_i)(\psi,\eta)=(\xi, K_{-2\rho}\phi)(\psi,\eta).$$ 
\medskip

Let us fix a complete set of irreducibles parametrised by $\Lambda^+$. Thus   for any $\mu\in\Lambda^+$, the conjugate of $V_\mu$ is $V_\lambda$, with $\lambda=-w_0\mu$. The composition
of the complex conjugation $J_\mu: V_\mu\to\overline{V_\mu}$   with a unitary intertwiner $U_\mu:\overline{V_\mu}\to V_\lambda$ is an antiunitary map $j_\mu: V_\mu\to V_\lambda$, unique up to scalar multiples by $z\in{\mathbb C}$ with $|z|=1$.
We can thus define a linear map, the antipode, 
 $S: {\mathcal C}\to {\mathcal C},$ by
$$S(\phi^*\otimes\psi)=({j_\mu{\psi}})^*\otimes j_\mu{\phi},\quad \phi^*\otimes\psi\in V_\mu^*\otimes V_\mu.$$ 
Notice that $S$ does not depend on the choice of $j_\mu$ or of the set of irreducibles.
We also define the counit $$\varepsilon: {\mathcal C}\to{\mathbb C}, \quad 
\varepsilon(\phi^*\otimes\psi)=(\phi, \psi).$$
\medskip

\noindent{\bf  6.8. Proposition} {\sl    Antipode $S$ and counit $\varepsilon$  make  ${\mathcal C}$  into a Hopf algebra. Furthermore, $S$ and $\varepsilon$   commute with the adjoint map
and $S^2(\phi^*\otimes\psi)=\phi^*K_{2\rho}\otimes K_{2\rho}^{-1}\psi$.}\medskip

\noindent{\sl Proof} The relations $S(a)^*=S(a^*)$ and $\overline{\varepsilon(a)}=\varepsilon(a^*)$ are   easy to check.
If $\eta_r$ is an orthonormal basis of $V_\mu$, with  $\mu$ and $j$ as above, then 
$\overline{r}=\sum_r \eta_r\otimes j{\eta_r}\in(\iota, V_\mu\otimes{V_\lambda})$. We verify the relation $m\circ 1\otimes S\circ\Delta(a)=\varepsilon(a)$ for $a=\phi\otimes\psi\in V_\mu^*\otimes V_\mu$.
$$m(1\otimes S(\Delta(\phi\otimes\psi))=$$
$$m(1\otimes S(\phi\otimes\eta_r\otimes\eta_r^*\otimes\psi))=$$
$$m((\phi\otimes\eta_r)\otimes ((j_\mu\psi)^*\otimes j_\mu\eta_r))=$$
$$\phi({j_\mu\psi})^*\otimes\overline{r}(1_{\mathbb C})=
\phi({j_\mu\psi})^*\circ\overline{r}\otimes 1_{\mathbb C}=$$
$$\sum_r\phi(\eta_r)(j_\mu{\psi},j_\mu{\eta_r})=\sum_r\phi(\eta_r)({\eta_r},\psi)=$$
$$\phi(\psi)=\varepsilon(\phi\otimes\psi).$$
 The analogous relation with $S$ on the left follows from this after taking the adjoint.
 We finally compute $S^2$. Keeping the same notation used in the definition of  $S$,  we have that $S^2(\phi^*\otimes\psi)=(j_\lambda j_\mu\phi)^*\otimes j_\lambda j_\mu\psi$. 
We identify $\overline{\overline{V_\mu}}$   with $V_\mu$ as an hermitian space. In other words, the corresponding conjugation identifies with ${J_\mu}^{-1}$.
 Observe that ${J_\mu}^{-1} U_\mu^*J_\lambda^{-1}$ is a unitary intertwiner of $(\overline{V_\lambda}, \overline{\overline{V_\mu}})$.
The action of $U_q({\mathfrak g})$ on $\overline{\overline{V_\mu}}$ is $a \xi=S^{-2}(a)\xi$, and we know that $S^{-2}(a)=K_{2\rho}a K_{2\rho}^{-1}$, hence $K_{2\rho}^{-1}\in( \overline{\overline{V_\mu}}, V_\mu)$. We can thus choose $U_\lambda=K_{2\rho}^{-1}{J_\mu}^{-1} U_\mu^*J_\lambda^{-1}$, implying $j_\lambda=K_{2\rho}^{-1}{J_\mu}^{-1} U_\mu^*$, which together with $j_\mu=U_\mu J_\mu$ gives $j_\lambda j_\mu=K_{2\rho}^{-1}$.
\medskip

\noindent We finally pass to the dual space,
${\mathcal C}'$
and we endow it with the usual dual algebra   structure  given by
$$\omega\omega':=\omega\otimes\omega'\circ\Delta$$
  coproduct
$\Delta:{\mathcal C}'\to({\mathcal C}\otimes{\mathcal C})'\supset{\mathcal C}'\otimes{\mathcal C}'$,
 $$\Delta(\omega)(a,b)=\omega(ab), \quad\quad a,b\in{\mathcal C},\quad\omega\in{\mathcal C}',$$
and involution given by duality with the involution of ${\mathcal C}$, 
$$\omega^*(a):=\overline{\omega(a^*)}, \quad \quad a\in{\mathcal C},\quad  \omega\in{\mathcal C}'.$$
Notice that duality for the involution differs from the case of a ordinary Hopf $^*$--algebra  (which is defined via  duality with    $a\to S(a^*)$.)
The   counit
$\varepsilon:{\mathcal C}'\to{\mathbb C}$ and   antipode $S:{\mathcal C}'\to{\mathcal C}'$ are defined as usual by
$$\varepsilon(\omega)=\omega(I),\quad\quad S(\omega)=\omega\circ S.$$

 \medskip

\noindent{\bf 6.9. Theorem } {\sl ${\mathcal C}'$ is
 isomorphic, as a $^*$-algebra, to the direct product  
of full matrix algebras
$${\mathcal C}'\simeq\prod M_{n_\lambda}({\mathbb C}), \quad n_\lambda=\text{dim}(V_\lambda).$$
 The coproduct $\Delta$
 is a homomorphism satisfying again  $\Delta(\omega^*)= \Delta^{\text{op}}(\omega)^*$.  
  }\medskip

Counit
  and   antipode   of  ${\mathcal C}'$ satisfy the same properties as for $U_q({\mathfrak g})$.
 We identify  the dual space of ${\mathcal C}\otimes {\mathcal C}$ with $\prod_{\lambda,\mu} M_{n_\lambda}\otimes M_{n_\mu}$ in the natural way, and in this   sense we understand the stated  homomorphism property of $\Delta$. Furthermore, the $R$ matrix can be found as
  an element of this algebra, see Sect. 1 in \cite{Sawin}, or  Sect. 2.5 in \cite{NT}. 

\section{The universal algebra ${\mathcal D}(V, \ell)$ }

 In this section the deformation parameter is a fixed root of unity of the form $q=e^{{i\pi}/{d\ell}}$. Our aim is to construct  non-associative bi-algebras ${\mathcal D}(V,\ell)$ endowed with an involution,  for ${\mathfrak g}\neq E_8$,  playing a role similar to that of ${\mathcal D}(V)$ of the generic case. \medskip

 Let $V$ be Wenzl's fundamental representation of ${\mathfrak g}$ ($V=V_{\kappa_1}\oplus V_{\kappa_2}$ in the type $D$ case.) Consider    the   (infinite dimensional) linear space,
      $${\mathcal D}(V, \ell)=\bigoplus_{n=0}^{\infty}( {V^{\otimes n}}){^*}p_n\otimes p_nV^{\otimes n}.$$
Notice that ${\mathcal D}={\mathcal D}(V, \ell)$  depends  not only on $V$ but also on   the root of unity.  
  We define a    multiplication on ${\mathcal D} $  by
  $$\alpha\beta:=\phi\phi' p_{m+n}\otimes p_{m+n}\psi\psi',$$
for  $$\alpha=\phi\otimes\psi\in {V^m}{^*}p_m\otimes p_mV^m, \quad \beta=\phi'\otimes\psi'\in {V^n}{^*}p_n\otimes p_nV^n.$$
In this way ${\mathcal D}$ becomes a unital but   not associative algebra, as, if we   pick a third element $\gamma=\phi''\otimes\psi''\in {V^r}{^*}p_r\otimes p_rV^r$, and we take into account the relation $p_{m+n+r}\circ p_{m+n}\otimes 1_r=p_{m+n+r},$ of Lemma 5.1, we see that
$$(\alpha\beta)\gamma=
(\phi\phi'\phi'')p_{m+n+r}\otimes
p_{m+n+r}(\psi\psi'\psi'')$$
 but
$$ \alpha(\beta\gamma)=(\phi \phi'\phi'')1_m\otimes p_{n+r}\circ p_{m+n+r}\otimes
p_{m+n+r}\circ 1_m\otimes p_{n+r}(\psi \psi'\psi'')$$
that differs from the previous one as in general
   $p_{m+n+r}\circ 1_m\otimes p_{n+r}\neq p_{m+n+r}$. 
   Notice however that the elements
 $$v_{\xi,\eta}:=\xi^*\otimes\eta, \quad\quad \xi,\eta\in V$$ 
   still generate ${\mathcal D}$ as an algebra,   
 $$\xi_1\dots\xi_np_n\otimes p_n\eta_1\dots\eta_n=(\dots((v_{\xi_1,\eta_1}v_{\xi_2,\eta_2})v_{\xi_3,\eta_3})\dots v_{\xi_n,\eta_n}).$$
   
    We next introduce an involution in  ${\mathcal D}$ as suggested by the generic case.
 Specifically, we replace the   boundary operators $\sigma_n$ by their truncated version:
 $$\tau_n:=p_n\sigma_np_n\in(V^{\underline{\otimes} n},  V^{\underline{\otimes} n}),$$
 which   satisfy  properties similar to those of $\sigma_n$  (except associativity).\medskip
 
 \noindent{\bf 7.1. Proposition} {\sl $\tau_n^*=\tau_n$, $\tau_n^2=p_n$.   }\medskip
 
 \noindent{\sl Proof}  Lemma 5.3 gives $p_n\sigma_n(1-p_n)\sigma_np_n=0$, implying $\tau_n^2=p_n\sigma_np_n\sigma_np_n=p_n\sigma_n^2p_n=p_n.$ Furthermore, $\tau_n^*=p_n\sigma_n^*p_n=p_n\sigma_np_n=\tau_n$.
 \medskip
 
 We can   write any element $\phi\in (V^{\otimes n})^*p_n$ uniquely in  the form
 $\phi=\psi^*$ with $\psi\in V^{\underline{\otimes} n}$, the adjoint of $\psi$ being computed with respect to
 Wenzl's inner product of $V^{\underline{\otimes} n}$.
 We set for $\psi,\psi'\in V^{\underline{\otimes} n}$,
 $$(\psi^*\otimes\psi')^*:={(\tau_n\psi')}^* \otimes\tau_n^{-1}\psi.$$
 \medskip

\noindent{\bf 7.2. Proposition} {\sl The following properties hold for the involution of ${\mathcal D}$.  For any $a\in {\mathcal D}$ and $v_{\xi,\eta}\in V^*\otimes V$,
\begin{itemize}
\item[{\rm a)}] $a\to a^*$ is antilinear,
\item[{\rm b)}] $a^{**}=a$,
\item[{\rm c)}] ${v_{\xi,\eta}}^*=v_{\eta,\xi}$,
\item[{\rm d)}] $(v_{\xi\eta}a)^* =a^*v_{\xi,\eta}^*$.
\end{itemize}
} \medskip

\noindent{\sl Proof} a) is obvious, b) is a consequence of involutivity of $\tau_n$, and c) of $\tau_1=1$. We show 
d) for $a=\phi^*\otimes\psi\in ( {V^{\otimes n}}){^*}p_n\otimes p_nV^{\otimes n}$. We have
$$a^*v_{\xi,\eta}^*=(\psi^*{\tau_n}^{-1}\otimes \tau_{n}^{-1}\phi)(\eta^*\otimes\xi)=$$
$$[\psi^*\eta^*\circ \tau_{n}^{-1}\otimes 1_V\circ p_{n+1}] \otimes [p_{n+1}\circ \tau_n^{-1}\otimes 1_V\circ \phi\xi ]=  $$
$$[(\eta\psi)^* \circ \sigma(V, {V^{\underline{\otimes}n}})^{-1}\circ \tau_n^{-1}\otimes 1_V\circ p_{n+1}]\otimes  [p_{n+1}\circ \tau_n^{-1}\otimes 1_V\circ \phi\xi ].$$
Notice that $1_V\otimes p_n\circ (1-p_{n+1})\circ \sigma(V, {V^{\underline{\otimes}n}})^{-1}\circ \tau_n^{-1}\otimes 1_V\circ p_{n+1}=0$
by the following lemma.  Hence we can  
insert the projection $p_{n+1}$ after $(\eta\psi)^*$ in the above computation,   use once again   Lemma 5.3 and the iterative definition of $\sigma_n$
to conclude that
$$a^*v_{\xi,\eta}^*= (\eta\psi)^*\circ\tau_{n+1}^{-1}\otimes  [p_{n+1}\circ \tau_n^{-1}\otimes 1_V\circ \phi\xi ].$$
Similar arguments show that $(v_{\xi,\eta} a)^*$ yields   the same expression.
\medskip

\noindent{\bf 7.3. Lemma} {\sl Assume ${\mathfrak g}\neq E_8$, and let $T\in(V^{\otimes m}, V^{\otimes n})$ be a negligible arrow of the tilting category ${\mathcal T}_\ell$. Then
 $$p_n\circ T\circ 1_V\otimes p_{m-1}=0.$$}\medskip

 \noindent{\sl Proof} Set $Y=p_n\circ T\circ 1_V\otimes p_{m-1}$.
    We  take into account the basic  property of Wenzl's fundamental representation $V$ for ${\mathfrak g}\neq E_8$ recalled in Theorem 4.4 a) according to which for any $\lambda\in\Lambda_\ell$, $V\otimes V_\lambda$ is completely reducible  (in a way that, although with same multiplicities, may differ from the decomposition of $V_\lambda\otimes V$)
  and the dominant weights
  $\mu$ of the irreducible components $V'_\mu$ all lie in $\overline{\Lambda}_\ell$.
  Thus the space of such  a component  with $\mu\in\Lambda_\ell$ must be in the kernel of $Y$ by Lemma 5.3. On the other hand, if $\mu\in\overline{\Lambda_\ell}\setminus\Lambda_\ell$ then   $Y V'_{\mu}=\{0\}$ as otherwise it would be an irreducible submodule of $p_n V^{\otimes n}$ of weight $\mu$. 
    \medskip

We    introduce a    coproduct   
$$\Delta: {\mathcal D}\to {\mathcal D}\otimes {\mathcal D}$$ in a way similar to the generic case, i.e.
by means of a pair $\eta_r\in p_n V^n$, $\eta^r\in(V^n)^*p_n\simeq (p_n V^n)^*$ of  dual bases:
$$\Delta(\phi\otimes\psi)=(\phi\otimes\eta_r)\otimes(\eta^r\otimes\psi),\quad \phi\otimes\psi\in (p_n V^n)^*\otimes (p_n V^n),$$
again independent 
on the choice of the basis.

\medskip

\noindent{\bf  7.4. Theorem} {\sl   The coproduct 
 $\Delta$ is  unital, coassociative,     and satisfies  for $a, b\in {\mathcal D}$, 
\begin{itemize}
\item[{\rm a)}]  $\Delta(a^*)= \Delta^{\text{op}}(a)^*$,
\item[{\rm b)}] $\Delta(ab)=\Delta(a)\Delta(b)$.
\end{itemize}
}\medskip

\noindent{\sl Proof} 
   a) For $a=\psi^*\otimes \psi'\in (p_n V^n)^*\otimes(p_n V^n)$, let $\psi_r$ be an orthonormal basis of $p_n V^n$ with respect to the Wenzl's inner product, and consider the dual basis $\xi^r:=\psi_r^*$. We compute 
$\Delta(a^*)$ with respect to the dual pair $\tau_n^{-1}\psi_r$, $(\tau_n\psi_r)^*$,
$$\Delta(a^*)=\Delta((\tau_n\psi')^*\otimes\tau_n^{-1}\psi)=\sum_r ((\tau_n\psi')^*\otimes\tau_n^{-1}\psi_r)\otimes((\tau_n\psi_r)^*\otimes\tau_n^{-1}\psi)=$$
$$  \sum_r (\psi_r^*\otimes\psi')^*\otimes(\psi^*\otimes \psi_r)^*= \Delta^{\text{op}}(a)^*.$$
b)    Let $\eta_s$ be 
  an orthonormal basis of $p_m V^m$ and $b=\xi^*\otimes\eta\in (p_m V^m)^*\otimes(p_m V^m)$. We have
$$\Delta(a)\Delta(b)=$$
$$\sum_{r,s}(\psi^*\xi^*p_{n+m}\otimes p_{n+m}\psi_r\eta_s)\otimes(\psi_r^*\eta_s^*p_{n+m}\otimes p_{n+m}\psi'\eta).$$
We   explicit  the middle term 
$$\sum_{r,s}p_{n+m}\psi_r\eta_s\otimes\psi_r^*\eta_s^*p_{n+m}=\sum_{h} \zeta_h\otimes \zeta^h $$
using  an orthonormal basis $\zeta_h$ of $p_{n+m}V^{n+m}$, where
$$\zeta^h:=(\zeta_h, p_{n+m}\psi_r\eta_s)\psi_r^*\eta_s^*p_{n+m}.$$
We adopt the   same notation as before to distinguish between a tensor product of inner products $(\zeta, \zeta')_p$ and Wenzl's inner product $(\zeta, \zeta')$ on  the  subspace $(p_nV^n)\otimes p_mV^m$ of $V^{n+m}$.
We have,
 $$(\zeta^h)^*=\sum_{r,s} (p_{n+m}\psi_r\eta_s,\zeta_h) (\psi_r^*\eta_s^*p_{n+m})^*=\sum_{r,s} (p_{n+m}\psi_r\eta_s,\zeta_h)p_{n+m}\overline{R}^{-1}\psi_r\eta_s=$$
 $$\sum_{r,s} (\psi_r\eta_s, \zeta_h)p_{n+m}\overline{R}^{-1}\psi_r\eta_s=\sum_{r,s}  p_{n+m}\overline{R}^{-1}(\psi_r\eta_s, \overline{R}\zeta_h)_p\psi_r\eta_s=$$
 $$p_{n+m}\overline{R}^{-1} \circ p_n\otimes p_m\circ\overline{R}\zeta_h=p_{n+m}\varepsilon^{-1}\circ p_m\otimes p_n\circ\varepsilon\circ p_{n+m}\zeta_h=\zeta_h,$$
by Lemma 5.3, where $$\varepsilon=\varepsilon({ p_nV^n, p_mV^m}).$$
 This shows that $$\Delta(a)\Delta(b)=  \sum_h \psi^*\xi^*p_{n+m}\otimes\zeta_h\otimes\zeta_h^*\otimes p_{n+m}\psi'\eta=  \Delta(ab).$$

\medskip

\noindent{\bf 7.5. Proposition} {\sl The linear map $\varepsilon: \phi^*\otimes\psi\in{\mathcal D}(V,\ell)\to(\phi,\psi)\in{\mathbb C}$ is a counit for the coproduct $\Delta$ making ${\mathcal D}(V,\ell)$ into a coalgebra.
It  is compatible with
the involution, $\varepsilon(a^*)=\overline{\varepsilon(a)}$.}\medskip

\noindent{\sl Proof} The former statement is easy to check. The latter follows from Prop. 7.1.\medskip

\noindent{\bf  7.6. Remark} It is important to notice  that, unlike the generic case,  $\varepsilon$ is {\it not} multiplicative on ${\mathcal D}(V,\ell)$, as a consequence of $p_n\neq 1$ in general.\bigskip

 \section{ The quantum groupoid ${\mathcal C}(G,\ell)$ and a corresponding  associative filtration}
 
 In analogy to the generic case,   we   follow a Tannakian reconstruction from the quotient category ${\mathcal F}_\ell$ for which  ${\mathcal D}={\mathcal D}(V,\ell)$ plays the role of   a universal algebra
 and we obtain a quantum groupoid ${\mathcal C}(G, \ell)$.
 The new main difference with the generic case (beyond associativity failure of ${\mathcal D}(V,\ell)$) is  the fact    that  the    ideal of  ${\mathcal D}$ defining ${\mathcal C}(G, \ell)$ is only a right ideal.  
 Correspondingly, ${\mathcal C}(G, \ell)$ is naturally only a $^*$-coalgebra, 
 and new effort is needed
  to construct an algebra structure in ${\mathcal C}(G, \ell)$ compatible with the coproduct, which we consider in the next sections.
 
  Specifically, in this section we define ${\mathcal C}(G,\ell)$ and establish   the main properties. We then 
   introduce an associated non-trivial,   possibly nilpotent, but associative structure, described by a finite sequence $\widetilde{\mathcal C}_k$ of $^*$-coalgebras, which we regard as a generalised algebra filtration.
 We shall eventually be able to give positive answers  to the above questions for ${\mathcal C}(G,\ell)$ in the type $A$ case by analysing this   filtration.\bigskip
 
 \noindent{\it a) The coalgebra ${\mathcal C}(G,\ell)$}\bigskip

\noindent We introduce   identifications in ${\mathcal D}$ arising from    arrows of  ${\mathcal F}_\ell$ as follows.
 Consider the linear span ${\mathcal J}$ of  ${\mathcal D}$  of  elements of the form  
  $$[\phi, A, \psi]:=\phi^*\otimes A\psi-\phi^*\circ A\otimes\psi,$$ where $A\in(V^{\underline{\otimes} m}, V^{\underline{\otimes} n})$ and set
  $${\mathcal C}(G, \ell)={\mathcal D}(V,\ell)/{\mathcal J}.$$
  We start summarising the  structure  that ${\mathcal C}$ inherits from ${\mathcal D}$.\medskip

  \noindent{\bf 8.1. Proposition} {\sl  ${\mathcal C}(G,\ell)$ is a finite dimensional,  coassociative, counital coalgebra with involution. More precisely, 
  \begin{itemize}
\item[{\rm a)}] ${\mathcal J}$ is a right ideal, 
\item[{\rm b)}]  ${\mathcal C}$ is finite dimensional and linearly spanned by class tensors $v^\lambda_{\phi,\psi}=[\phi^*\otimes\psi]$, $\phi$, $\psi\in V_\lambda$, $\lambda\in\Lambda_\ell$,
\item[{\rm c)}] ${\mathcal J}$ is a  coideal annihilated by $\varepsilon$,
\item[{\rm d)}]  ${\mathcal J}$ is $^*$-invariant, hence the involution of ${\mathcal D}$ factors through ${\mathcal C}$ and satisfies
 $$\Delta(a^*)= \Delta^{\text{op}}(a)^*,\quad \varepsilon(a^*)=\overline{\varepsilon(a)}, \quad a\in{\mathcal C},$$
 $$\Delta(v_{\phi, \psi}^\lambda)=\sum_r v_{\phi, \eta_r}^\lambda\otimes v_{\eta_r, \psi}^\lambda,\quad (v_{\phi,\psi}^\lambda)^*=v_{\psi,\phi}^\lambda.$$

\end{itemize}}\medskip
   
  \noindent{\sl Proof} a) follows from   Lemma 5.1.  
  Specifically, for  
  $$x=\phi\otimes A\circ\psi-\phi\circ A\otimes\psi \in{\mathcal J}$$ and any simple tensor
 $$\zeta:=\xi\otimes\eta\in ( {V^{h}}){^*}p_h\otimes p_hV^h,$$
 we have
$$x\zeta=\phi\xi p_{n+h}\otimes p_{n+h}(A\psi)\eta-(\phi A)\xi p_{m+h}\otimes p_{m+h}\psi\eta=$$
$$\phi\xi p_{n+h}\otimes p_{n+h}\circ A\otimes 1_{V^h}(\psi\eta)-(\phi\xi )A\otimes 1_{V^h}\circ p_{m+h}\otimes p_{m+h}\psi\eta=$$
$$\phi\xi p_{n+h}\otimes  (A\underline{\otimes} 1) p_{m+h}
\psi\eta-\phi\xi p_{n+h}(A\underline{\otimes} 1 )\otimes p_{m+h}\psi\eta\in{\mathcal J}.$$
 To show b)  we may argue as in the generic case but now with a choice  of isometries
    $S_i\in (V_{\lambda_i}, V^{\underline{\otimes} n})$, with $\lambda_i\in\Lambda_\ell$,  in the   $C^*$--category ${\mathcal F}_\ell$ satisfying $\sum_i S_i{S_i}^*=p_n$.     The remaining statements can be proved   in analogy with the generic case, taking into account the results of the previous section.\bigskip
    
    \noindent{\it b) The filtration $\widetilde{{\mathcal C}}_k$}\bigskip
    
\noindent We filter ${\mathcal D}$ by the size of tensor products.   Set 
  $${\mathcal D}_k=\bigoplus_{n\leq k}( {V^{\otimes n}}){^*}p_n\otimes p_nV^{\otimes n}.$$
 \medskip
 
 \noindent{\bf 8.2. Proposition} {\sl \begin{itemize}
\item[{\rm a)}] ${\mathcal D}_k$ is a filtration of ${\mathcal D}$, i.e. it is an increasing sequence of subspaces satisfying
 $${\mathcal D}_0={\mathbb C},\quad {\mathcal D}_h{\mathcal D}_k\subset {\mathcal D}_{h+k},\quad \bigcup_{k=0}^\infty {\mathcal D}_k={\mathcal D}.$$
 \item [{\rm b)}] ${\mathcal D}_k$ are $^*$-invariant subcoalgebras:
 $({\mathcal D}_k)^*={\mathcal D}_{k},\quad \Delta({\mathcal D}_{k})\subset {\mathcal D}_{k}\otimes {\mathcal D}_{k}.$
 \end{itemize}
 }\medskip

   \noindent Set
 $${\mathcal J}_k:=\text{l.s.}\{[\phi, A, \psi],\quad A\in  (V^{\underline{\otimes}m}, V^{\underline{\otimes}n}), \psi\in V^{\underline{\otimes}m},\phi\in V^{\underline{\otimes}n},m,n\leq k \}$$
and notice that ${\mathcal J}_k\subset {\mathcal J}$ and ${\mathcal J}_k\subset {\mathcal J}_{k+1}$. Set ${\mathcal C}_k={\mathcal D}_k/{\mathcal J}_k$.  
Hence there are obvious maps  
${\mathcal C}_k\to{\mathcal C}_{k+1}$, and ${\mathcal C}_k\to{\mathcal C}$.
\medskip

 \noindent{\bf 8.3. Lemma} {\sl ${\mathcal J}\cap{\mathcal D}_k={\mathcal J}_k$. Furthermore, ${\mathcal J}_{k+1}\cap{\mathcal D}_k={\mathcal J}_k$. }\medskip
 
 \noindent{\sl Proof} It suffices to show the first relation. The inclusion ${\mathcal J}_k\subset{\mathcal J}\cap{\mathcal D}_k$ is obvious. Let now $X\in{\mathcal J}\cap{\mathcal D}_k$ be written as a finite sum of the spanning set of ${\mathcal J}$. We can assume that   the indices  $m$, $n$     appearing in it   satisfy 
    $\text{min}\{m, n\} \leq k$.
  After subtracting to $X$ some element of ${\mathcal J}_k$, we may assume that $X$ takes the form
 $$X=\sum (\phi\otimes A\circ\psi-\phi\circ A\otimes\psi)+\sum \xi\otimes A'\circ\eta-\xi\circ A'\otimes\eta$$
 where $A\in(V^{\underline{\otimes}m}, V^{\underline{\otimes}n})$, $A'\in(V^{\underline{\otimes}q}, V^{\underline{\otimes}r})$, $n, q\leq k$,  $m, r>k$. We reduce $m$ and $r$ as follows.
 Consider isometries $S\in (V_\lambda,  V^{\underline{\otimes}m})$, where for simplicity we have dropped indices to $S$ and $\lambda$, which have pairwise orthogonal ranges which sum up to $p_m$, and similarly
 consider $T\in(V_\mu, V^{\underline{\otimes}q})$. Notice that the indices of the domains of $S$ and $T$ satisfy the required bounds.
 We write the first sum in $X$ as
 $$\sum (\phi\otimes [AS]\circ[S^*\psi]-\phi\circ [AS]\otimes[S^*\psi])+\sum(\phi\circ[AS]\otimes [S^*\psi]-\phi A\otimes\psi)$$
 and similarly for the second sum.
 The first sum above now lies in the span of the claimed spanning set, hence we are left to show that the remaining part,
 now written  as
 $$Y=\sum([\phi A]S\otimes S^*\psi-\phi A\otimes\psi)+\sum(\xi\otimes A'\eta-\xi T\otimes T^*[A'\eta])$$
 vanishes. But $Y\in{\mathcal J}\cap{\mathcal D}_k$, and the domain of $A$ and the range of $A'$ have large indices, hence the sum of the terms in second and third position must vanish. We are thus reduced to show the 
  general statement that given  elements $\phi_i$ and $\psi_i$ of a Hilbert space $H$ and operators $Y_i:K_i\to H$, on  another Hilbert spaces $K_i$ such that   $\sum_i \phi_i^*\otimes\psi_i=0$ then $\sum_{i,r} \phi_i^*Y_i\otimes {Y_i}^*\psi_i=0$,
  and this can now be checked by means of orthonormal base decomposition of $\phi_i$, $\psi_i$.
 \medskip

\noindent{\bf 8.4. Proposition} {\sl  The natural maps
${\mathcal C}_k\to{\mathcal C}_{k+1}$ are faithful and form a finite increasing sequence of   inclusions of $^*$--invariant subcoalgebras
 $${\mathbb C}= {\mathcal C}_0\subset {\mathcal C}_1\subset{\mathcal C}_2\subset\dots\subset {\mathcal C}_m={\mathcal C}$$
 stabilizing to ${\mathcal C}$.}\medskip
 
 \noindent{\sl Proof} The natural maps ${\mathcal C}_k\to{\mathcal C}_{k+1}$, ${\mathcal C}_k\to{\mathcal C}$ are faithful, by the previous lemma.
On the other hand, there is an integer $m$ such that every irreducible $V_\lambda$ with $\lambda\in\Lambda_\ell$ is contained in ${V^{\otimes n}}$ for some $n\leq m$, hence ${\mathcal C}_k$ must stabilize to ${\mathcal C}_k={\mathcal C}$ for $k\geq m$.
  \medskip

We next pass from  ${\mathcal C}_k$ to a quotient $\widetilde{{\mathcal C}}_k$ and from inclusions ${\mathcal C}_k\to{\mathcal C}_{k+1}$ to (possibly non-faithful) linear maps
$\widetilde{{\mathcal C}}_k\to\widetilde{{\mathcal C}}_{k+1}$. The advantage  will be existence of natural, even associative,  multiplication maps $\widetilde{{\mathcal C}}_h\otimes \widetilde{{\mathcal C}}_k\to \widetilde{{\mathcal C}}_{h+k}$.
  
 More precisely,    we observe  that   associativity failure of ${\mathcal D}$ can be described in terms of certain negligible intertwiners $Z$
   of the tilting category, as follows.\medskip
   
   \noindent{\bf 8.5. Lemma} {\sl For any triple $\alpha=\phi\otimes\psi$, $\beta=\phi'\otimes\psi'$, $\gamma=\phi''\otimes\psi''$ of elements
   of ${\mathcal D}$, of grades $m$, $n$, $r$ respectively,
       $$(\alpha\beta)\gamma-\alpha(\beta\gamma)=$$
$$(\phi\phi'\phi'')p_{m+n+r}\otimes
Z(\psi\psi'\psi'')+(\phi\phi'\phi'')\circ Z^*\otimes
p_{m+n+r}(\psi p_{n+r}(\psi'\psi'')),$$
where
$$Z=p_{m+n+r}\circ1_m\otimes (1_{n+r}-p_{n+r}). $$}\medskip

   \noindent Consider
  the  following spaces of   negligible arrows of the tilting category,
  $$  {\mathcal Z}^{(k)}:=\{p_{q+j+r}\circ1_{V^q}\otimes (1-p_j)\otimes 1_{V^r};    q+j+r\leq k\}$$
      and then define   
$$\widetilde{\mathcal J}_k:=\text{l.s.}\{{\mathcal J}_k, \quad\phi\otimes Z\circ\psi', \quad\phi'\circ 
(Z')^*\otimes\psi \},$$
\noindent where 
  $Z$, $Z'$ vary in ${\mathcal Z}^{(k)},$
  $\psi$, $\phi$   in  the canonical truncated tensor powers of $V$, but   $\psi'$  and $\phi'$ belong to the full tensor powers.

\medskip

\noindent{\bf  8.6. Proposition} {\sl 
We have that 
${\mathcal D}_j\widetilde{\mathcal J}_k, \quad {\widetilde{\mathcal J}}_k {\mathcal D}_j\subset {\widetilde{\mathcal J}}_{j+k}.$  
 }\medskip

\noindent{\sl Proof}  Arguments similar to those of Prop. 8.1 a), but keeping track of the grades of homogeneous
elements, show that ${{\mathcal J}}_k {\mathcal D}_j\subset {{\mathcal J}}_{j+k}$, and hence ${{\mathcal J}}_k {\mathcal D}_j$ is a subspace of ${\widetilde{\mathcal J}}_{j+k}$.
  Similar considerations hold for products $y\zeta$  with  $y\in\widetilde{\mathcal J}_k$ of the form $y=\phi\otimes Z\psi'$ or
  $y=\phi\circ Z^*\otimes\psi$ and $\zeta:=\xi\otimes\eta\in ( {V^{h}}){^*}p_h\otimes p_hV^h,$, $h\leq j$, noticing that  in the first case for example
  $$y\zeta=\phi\xi p_{n+h}\otimes p_{n+h}\circ Z\otimes 1_{V^h}(\psi'\eta)$$ and that the map
   $Z\to p_{n+h}\circ Z\otimes 1_{h}$ for $Z=p_r\circ 1_u\otimes (1-p_s)\otimes 1_t$, $r=u+s+t\leq k,$ takes 
   ${\mathcal Z}^{(k)}$ to ${\mathcal Z}^{(k+j)}$.
  The left ideal property is more delicate to check, due
  to lack of associativity of the projections $p_n$. 
In order to check that $\zeta y\in{ \widetilde{\mathcal J}}_{j+k}$ we compute  
$$\zeta y=\xi\phi p_{h+n}\otimes p_{h+n}\circ 1_h\otimes Z(\eta\psi').$$
 Write $p_{h+n}\circ 1_h\otimes Z$ in the form
  $$p_{h+n}\circ 1_h\otimes Z=
  Z_1-Z_2\circ 1_{h+q}\otimes(1-p_j)\otimes 1_r$$
with 
  $Z_1=p_{h+n}\circ 1_{h+q}\otimes(1-p_j)\otimes 1_r$, $Z_2=p_{h+n}\circ 1_h\otimes(1-p_n)$
    and notice that $Z_i$ both lie in ${\mathcal Z}^{(j+k)}$ and this implies $\zeta y\in{ \widetilde{\mathcal J}}_{j+k}$.  
We are left to show that $\zeta x\in{ \widetilde{\mathcal J}}_{j+k}$,  for all $\zeta\in {\mathcal D}_j$ as above
and $x=\phi\otimes A\circ\psi-\phi\circ A\otimes\psi \in{\mathcal J}_k$.
We claim that it suffices to the take $j=1$ and $\zeta=v_{\xi, \eta}$  with $v_{\xi,\eta}:=\xi\otimes\eta\in V^*\otimes V$. Indeed, we have already noticed that  finite
 products  
 $ (\dots((v_{\xi_1,\eta_1}v_{\xi_2,\eta_2})v_{\xi_3,\eta_3})\dots v_{\xi_h,\eta_h}),$ are total in $(V^h)^*p_h\otimes p_h V^h$, 
 and   multiplication of ${\mathcal D}$ is associative up to summing elements of ${ \widetilde{\mathcal J}}_{j+k}$, by Lemma 8.5. We thus compute
 $$v_{\xi,\eta} x=\xi\phi p_{1+n}\otimes [p_{1+n}\circ  1_{V}\otimes A(\eta\psi)]-[(\xi\phi ) \circ 1_{V}\otimes A\circ p_{1+m}]\otimes p_{1+m}\eta \psi=$$
$$\xi\phi p_{1+n}\otimes    (1\underline{\otimes} A)(p_{1+m}\eta\psi)-(\xi\phi p_{1+n}) (1\underline{\otimes} A)\otimes p_{1+m}\eta \psi+$$
$$\xi\phi p_{1+n}\otimes Y_1(\eta\psi)-(\xi\phi )Y_2\otimes p_{1+m}\eta \psi$$
where
$$Y_1=p_{1+n}\circ  1_{V}\otimes A\circ (1-p_{1+m}),\quad Y_2=(1-p_{1+n})\circ 1_{V}\otimes A\circ p_{1+m}.$$
Hence $v_{\xi,\eta} x\in{ \widetilde{\mathcal J}}_{k+1}$, provided we show that
 the last two terms in the sum vanish, or, equivalently, that   
  $Y_1\circ 1_V\otimes p_m=0=1_V\otimes p_n\circ Y_2$, but this follows from Lemma 7.3,
  and the proof   is complete.

\medskip

 \noindent{\bf 8.7. Proposition} {\sl Each subspace   $\widetilde{{\mathcal J}}_k$ is $^*$--invariant. }\medskip
 
 \noindent{\sl Proof} As already noticed in the proof of Prop. 8.1, the proof of $^*$-invariance of ${\mathcal J}_k$   is similar to that of the generic case.
  For $Z\in{\mathcal Z}^{(k)}\cap(V^n, V^{\underline{\otimes} n})$, $\phi\in V^n$, $\psi\in p_nV^n$,  
$$(\psi^*\otimes Z\phi)^*=(\tau_nZ\phi)^*\otimes \tau_n^{-1}\psi=    $$
$$ [(\tau_nZ\phi)^*\otimes \tau_n^{-1}\psi-(\tau_nZ\phi)^*\tau_n^{-1}\otimes \psi]+\phi^*\circ Z^*\otimes \psi$$
 which thus  lies in $\widetilde{{\mathcal J}}_k$.
One can similarly show that 
$((\psi\circ Z')^*\otimes \phi )^*\in\widetilde{{\mathcal J}}_k$ as well for $Z'\in({{\mathcal Z}^{(k)}})^*$   and thus conclude that $\widetilde{{\mathcal J}}_k$ is $^*$-invariant.
 \medskip

\noindent{\bf  8.8. Proposition} {\sl  We have that  $\Delta(\widetilde{{\mathcal J}}_k)\subset \widetilde{{\mathcal J}}_k\otimes {\mathcal D}_k+{\mathcal D}_k\otimes\widetilde{{\mathcal J}}_k$.}\medskip
 
\noindent{\sl Proof} 
 This  can be proved similarly to  Prop. 6.5,   keeping track  of the grades of homogeneous elements.\medskip

 We   set $\widetilde{{\mathcal C}}_k:={\mathcal D}_k/\widetilde{{\mathcal I}}_k$ for $k\in{\mathbb N}$. The composition of 
 the natural  linear
 inclusion ${\mathcal D}_h\to{\mathcal D}_k$,  $k>h$, with projection ${\mathcal D}_k\to \widetilde{{\mathcal C}}_k$   factors
 through a linear map
 $$\widetilde{{\mathcal C}}_h\to\widetilde{{\mathcal C}}_k$$ and this is an inductive system.
 We have natural quotient maps
 $${\mathcal C}_k\to\widetilde{\mathcal C}_k$$
 and we denote by $e^\lambda_{i,j}\in \widetilde{{\mathcal C}}_k$ the image of the matrix coefficient $v^\lambda_{i,j}\in{{\mathcal C}}_k$ corresponding to an orthonormal basis of $V_\lambda$.
 Let $\Lambda_\ell^k$ denote the set of $\lambda\in\Lambda_\ell$ for which $V_\lambda$ is a summand of some $V^{\underline{\otimes} n}$ with $n\leq k$.
In analogy with the properties  of Prop. 8.1 for ${\mathcal C}$,  we also summarise the results of the last two sections for $\widetilde{\mathcal C}_k$, which now take a stronger form.
 \medskip
 
 \noindent{\bf 8.9. Theorem} {\sl Assume that ${\mathfrak g}\neq E_8$ and let $V$ be Wenzl's  fundamental representation of ${\mathfrak g}$  ($V=V_{\kappa_1}\oplus V_{\kappa_2}$ in the type $D$ case). 
 Then
 \begin{itemize}
\item[{\rm a)}] $\widetilde{{\mathcal C}}_k$ is a $^*$--coalgebra linearly spanned by elements $e^{\lambda}_{i,j}$  labelling matrix units corresponding to $V_\lambda$, for $\lambda\in\Lambda_\ell^k$,
\item[{\rm b)}]
coproduct and involution satisfy   
$$\Delta(e^\lambda_{i,j})=\sum_r e^{\lambda}_{i,r}\otimes e^\lambda_{r,j}$$
$$(e^\lambda_{i,j})^*=e^\lambda_{j,i},$$
in particular the involution is anticomultiplicative,
\item[{\rm c)}] there are  associative   multiplication maps $\widetilde{{\mathcal C}}_h\otimes\widetilde{{\mathcal C}}_k\to\widetilde{{\mathcal C}}_{h+k}$ 
and an   element $I\in\widetilde{{\mathcal C}}_0$ acting as the identity. The   involution is antimultiplicative and the coproduct is unital and multiplicative.
 \end{itemize}}

\section{Quasi-coassociative dual $C^*$--quantum groupoids   $\widehat{{\mathcal C}(G,\ell)}$}

 The aim of the present and the next section is to show the main result of the paper,  stating that if $G={\rm SU}(N)$ then    the dual groupoid $\widehat{{\mathcal C}(G,\ell)}$ can be made into a    $C^*$--quantum groupoid,
   satisfying the axioms of a weak quasi-Hopf $C^*$-algebra in the sense of \cite{MS}. Furthermore the representation category $\text{Rep}_V(\widehat{\mathcal C})$ of $\widehat{{\mathcal C}}$  generated by the fundamental representation turns out to be a tensor $C^*$--category equivalent to the original fusion category ${\mathcal F}_\ell$.
 
 We shall divide the proof in two parts. Throughout this  section ${\mathfrak g}$ is general (but $\neq E_8$) and we assume to know  that ${\mathcal C}(G,\ell)$ is cosemisimple with respect to the coalgebra structure introduced in Sect. 8.
We then show that  the above conclusions hold for  $\widehat{{\mathcal C}(G,\ell)}$.  More precisely, in Subsect. a) we upgrade  ${\mathcal C}(G,\ell)$ to a non-associative bi-algebra with antipode  associated to a fixed section of the quotient map ${\mathcal D}(V,\ell)\to{\mathcal C}(G,\ell)$, while in b) we pass to the dual groupoid $\widehat{{\mathcal C}(G,\ell)}$, in c) we construct the Drinfeld's associator of $\widehat{{\mathcal C}(G,\ell)}$ and discuss the main properties,   in  d) we explicit  quasi-invertible $R$ matrices for $\widehat{{\mathcal C}(G,\ell)}$,   in e) we briefly discuss the  relation between the groupoid structures associated to different sections, and finally in f) we show that $\text{Rep}_V(\widehat{\mathcal C})$ is a tensor $C^*$--category equivalent to   the fusion category ${\mathcal F}_\ell$.
  In the next section we verify cosemisimplicity in the type $A$ case.
   \bigskip

 \noindent{\it a) Algebra structure and antipode in ${\mathcal C}(G,\ell)$}\bigskip

\noindent  Let $V_\lambda$ be a copy of an irreducible representation of $U_q({\mathfrak g})$ with highest weight $\lambda\in\Lambda_\ell$ and contained in some $V^{\underline{\otimes} n}$, and let
$M_\lambda$ denote the image of $V_\lambda^*\otimes V_\lambda$ in ${{\mathcal C}}$ under the quotient map ${\mathcal D}\to{{\mathcal C}}$, which we already know to be subcoalgebras independent of the choice of $V_\lambda$ and spanning ${\mathcal C}$. 
  We shall say that ${\mathcal C}(G,\ell)$ is {\it cosemisimple} if 
the subcoalgebras $M_\lambda$   are   independent matrix coalgebras in ${{\mathcal C}}$ as $\lambda$ varies in $\Lambda_\ell$, always understood of full dimension $\text{dim}(V_\lambda)^2$.   \medskip

If we know that  ${\mathcal C}(G,\ell)$ is cosemisimple, we can   endow  it both with an invertible antipode and with a non-associative algebra structure. Let's start with the antipode, which we introduce in a way similar to  the generic case. Fix a complete set $V_\lambda$, $\lambda\in\Lambda_\ell$ of irreducibles contained in the various $V^{\underline{\otimes} n}$, and   set, for  $\phi^*\otimes\psi\in V_\lambda^*\otimes V_\lambda$ and $\lambda\in\Lambda_\ell$,
 $$ S(v^\lambda_{\phi,\psi})=v^{\overline{\lambda}}_{j\psi, j\phi}, $$
 still independent of the choice of $V_\lambda$.
It  satisfies the  relations
    $$S(a^*)=S(a)^*,\quad \quad S^2(v^\lambda_{\phi,\psi})=v^\lambda_{K\phi, K\psi},$$
$$\Delta \circ S=S\otimes S\circ  \Delta^{\text{op}}.$$ 
   where $K=K_{-2\rho}$.
 
As regards   the algebra structure,  we pull back the product of ${\mathcal D}(V,\ell)$ via   the choice of a   section $s: {\mathcal C}(G,\ell)\to{\mathcal D}(V,\ell)$ of the quotient map ${\mathcal D}(V,\ell)\to {\mathcal C}(G,\ell)$.
 Correspondingly, we have a  choice of irreducibles $V_\lambda$ and $s$ takes $v^\lambda_{\phi,\psi}$ to $\phi^*\otimes\psi$, for $\phi,\psi\in V_\lambda$.  We thus set
 $$v^\lambda_{\phi,\psi}v^\mu_{\xi,\eta}=[s(v^\lambda_{\phi,\psi})s(v^\mu_{\xi,\eta})].$$
 We always choose $V_0={\mathbb C}$, and $V_\kappa=V$ for
 the trivial and  fundamental representation, respectively.
 In this way, denoting as before by $v_{\xi,\eta}$ the coefficients of $V$, products of the form
$v^\lambda_{\phi,\psi}v_{\xi,\eta}$ encode fusion decomposition of the quotient category ${\mathcal F}_\ell$. However, the section is not unique, and the product of ${\mathcal C}(G,\ell)$ does depend on $s$ (but cf. Subsect. e).)

  \medskip

\noindent{\bf 9.1. Proposition} {\sl The product makes ${\mathcal C}(G,\ell)$ into a (non-associative) unital algebra and  the coproduct of  ${\mathcal C}(G,\ell)$ is a unital homomorphism, and furthermore the following relation holds, 
 $$m\circ 1\otimes S\circ\Delta=\varepsilon=m\circ S\otimes 1\circ\Delta.\eqno(9.1)$$

}\medskip

\noindent{\sl Proof} The first statement  follows from Prop. 7.4 b). The left hand side of $(9.1)$ can be proved with computations  similar those in the proof of Prop. 6.8, with the only variation that now  
$m\circ 1\otimes S\circ\Delta(v^\mu_{\xi,\eta})=\xi^*({j_\mu\eta})^*\circ p_{h+m}\otimes p_{h+m}\circ\overline{r}(1_{\mathbb C})$. It suffices to notice that  
  $\overline{r}$ always lies in the range of the Wenzl projection $p_{h+n}$, where
$h$ and $n$ correspond to the powers of $V$ containing $V_\lambda$ and $V_{\overline{\lambda}}$, respectively. More details on this argument can be found in the proof of the following Lemma 10.4.

\medskip

Hence ${\mathcal C}(G,\ell)$ satisfies all the axioms of a Hopf algebra except associativity of the product and multiplicativity of the counit. The antipode is not antimultiplicative.

\bigskip

 \noindent{\it b) The dual quantum groupoid $\widehat{{\mathcal C}(G,\ell)}$}\bigskip
 
\noindent We next study the algebra structure associated to  a fixed section $s: {\mathcal C}(G,\ell)\to{\mathcal D}(V,\ell)$.
 Notice that, as an algebra, ${\mathcal C}(G,\ell)$     is quite far from admitting   an interpretation   as a non commutative space, as compared, e.g., to the compact quantum groups of Woronowicz, in that it may lack some important properties, e.g. antimultiplicativity of the involution, or a $C^*$-norm, not to mention associativity. It is far more rewarding   to pass to the dual $\widehat{{\mathcal C}(G,\ell)}$, and  correspondingly consider its $^*$--representations.
 In this subsection we show that $\widehat{{\mathcal C}(G,\ell)}$ satisfies most properties of the weak quasi-Hopf $C^*$--algebras
 of \cite{MS}.  
  \medskip

   We identify elements  of tensor  power algebras $\widehat{\mathcal C}^{\otimes n}$    with functionals on   ${\mathcal C}^{\otimes n}$. We shall need various elements of these algebras, and we start with 
   $P\in\widehat{\mathcal C}^{\otimes 2}$, defined as follows
   $$P(v^{\lambda}_{\phi,\psi}, v^\mu_{\xi,\eta})=(\phi\otimes \xi,  p_{h+k}\psi\otimes \eta)_{p,2},$$ 
    where $h$, $k$  are such that $V_\lambda$, $V_\mu$   are summands of $V^{\underline{\otimes}h}$,  $V^{\underline{\otimes}k}$,     respectively.
 Notice that these integers are specified by the   section.  Furthermore, 
 the   forms defining $P$ is   understood  with respect to  the product form of 
 $V^{\underline{\otimes}h}\otimes V^{\underline{\otimes}k}$, each factor in turn endowed with Kirillov-Wenzl inner product.

   \medskip
  
   \noindent{\bf 9.2. Theorem } {\sl Endowed with dual bi-algebra structure, antipode and involution,
    $\widehat{\mathcal C}=\widehat{{\mathcal C}(G,\ell)}$ is isomorphic, as a $^*$--algebra,   to
  $$\widehat{\mathcal C}\simeq \bigoplus_{\lambda\in\Lambda_\ell}{\mathcal L}(V_\lambda),$$
  hence it is a $C^*$-algebra. Furthermore,
$(\widehat{\mathcal C}, \Delta, \varepsilon, S)$ satisfies all the axioms of a Hopf algebra except coassociativity and unitality of $\Delta$.
In particular, the following compatibility properties hold, for $\omega, \tau\in\widehat{\mathcal C}$,
 $$ \Delta(I)=P,\quad\quad\quad \Delta(\omega\tau)=\Delta(\omega)\Delta(\tau),\eqno(9.2)$$
  $$\varepsilon(\omega\tau)=\varepsilon(\omega)\varepsilon(\tau),\quad\quad\varepsilon(\omega^*)=\overline{\varepsilon(\omega)},\eqno(9.3)$$
 $$S(\omega\tau)=S(\tau)S(\omega),\quad\quad S(\omega^*)=S(\omega)^*,\eqno(9.4)$$

  }\medskip

 \noindent{\sl Proof} The first or second property stated in $(9.i)$ shall be referred to as $(9.i)_1$ or $(9.i)_2$ respectively.
 We  
 show   $(9.2)_1$. By duality, the identity $I$ of $\widehat{\mathcal C}$ is the counit $\varepsilon$ of ${\mathcal C}$, and $\Delta(I)$ is the two-variable functional given by
 $$\Delta(I)(v^\lambda_{\phi,\psi}, v^\mu_{\xi,\eta})=\varepsilon(v^\lambda_{\phi,\psi}v^\mu_{\xi,\eta})=\varepsilon([\phi^*\xi^*p_{h+k}\otimes p_{h+k}\psi\eta])=$$
 $$\phi^*\xi^*p_{h+k}(p_{h+k}\psi\eta)=(\phi\xi, p_{h+k}\psi\eta)_{p, 2}=P(v^{\lambda}_{\phi,\psi}, v^\mu_{\xi,\eta}).$$  
The remaining properties   follow from duality and   corresponding properties of ${\mathcal C}(G,\ell)$.  More precisely, 
the associative algebra structure and the non-coassociative but counital coalgebra structure of   $\widehat{\mathcal C}$ obviously follow from the coassociative coalgebra   and non-associative algebra structure    in ${\mathcal C}$, respectively. Antimultiplicativity of the involution of $\widehat{\mathcal C}$ follows from anticomultiplicativity of that of ${\mathcal C}$.
 The antipode axioms $m\circ S\otimes 1\circ\Delta=1=m\circ 1\otimes S\circ\Delta$,  $(9.2)_2$, and $(9.4)_2$ are selfdual and already shown in ${\mathcal C}$;  $(9.3)_1$ and $(9.4)_1$ correspond  respectively to unitality of the coproduct  and anticomultiplicativity of the antipode in ${\mathcal C}$ (Prop. 9.1; Sect. 9a). Property $(9.3)_2$ is obvious. \medskip

The following Prop. 9.9 explicits  a compatibility relation between $\Delta$ and $^*$-involution in terms of the $R$-matrix.
 For completeness we remark that it may be $S\otimes S\circ\Delta\neq  \Delta^{\text{op}}\circ S$ in $\widehat{{\mathcal C}(G,\ell)}$, as this corresponds to antimultiplicativity failure of the antipode of ${{\mathcal C}(G,\ell)}$. However, this fact plays no role in the representation theory of $\widehat{{\mathcal C}(G,\ell)}$.
\bigskip

\noindent{\it c) Drinfeld's associator}\bigskip

\noindent 
For a given weight $\lambda\in\Lambda_\ell$,  let $h_\lambda$  denote the truncated powers of $V$ containing $V_\lambda$, as prescribed by the choice of a  section $s$. 
It will be useful for later computations to have a multiplication rule for   
  elements of $\widehat{\mathcal C}^{\otimes n}$ of  the following form. Let $T=(T_n)$ be  a sequence of   linear maps $T_{n}: V^{\otimes n}\to V^{\otimes n}$
  and associate  
  the 
 element $\omega_T$ of $\widehat{\mathcal C}^{\otimes n}$ defined by
  $$\omega_T(v^{\lambda_1}_{\phi_1,\psi_1},\dots ,v^{\lambda_n}_{\phi_n,\psi_n})=(\phi_1\otimes \dots\otimes \phi_n, T_{h_{\lambda_1}+\dots+h_{\lambda_n}}\psi_1\otimes \dots\otimes \psi_n)_{p,n},$$
 where  the form is a tensor product form with $n$ factors.  The following lemma is  a convenient formulation of the generalisation to tensor powers of the direct sum decomposition of the algebra  $\widehat{\mathcal C}$ 
 into its simple components.
\medskip
 
 \noindent{\bf 9.3. Lemma} {\sl Given $S=(S_n)$  and $T=(T_n)$ as above, set $\omega=\omega_S\omega_T$.
 Then      $$\omega(v^{\lambda_1}_{\phi_1,\psi_1},\dots,v^{\lambda_n}_{\phi_n,\psi_n})=(\phi_1\dots\phi_n, S_{h_{\lambda_1}+\dots+h_{\lambda_n}}\circ p_{h_{\lambda_1}}\otimes\dots\otimes 
  p_{h_{\lambda_n}}\circ T_{h_{\lambda_1}+\dots+h_{\lambda_n}}\psi_1\dots\psi_n)_{p,n}.$$ }\medskip
  
\noindent  The following elements  
  $\Phi$, $\Psi\in\widehat{\mathcal C}^{\otimes 3}$ are  important examples. Set $$q_{h_\lambda, h_\mu, h_\nu}=\sum_{\gamma, i}1_{h_\lambda}\otimes S_{\gamma, i}\circ p_{h_\lambda+h_\gamma}\circ 1_{h_\lambda}\otimes S_{\gamma,i}^*,$$ 
  where    $S_{\gamma,i}\in(V_\gamma, V^{\underline{\otimes}(h_\mu+h_\nu)})$ are isometries of the fusion category satisfying $\sum_{\gamma, i}S_{\gamma, i}S_{\gamma, i}^*=p_{h_\mu+h_\nu}$. Observe that $q_{h_\lambda, h_\mu, h_\nu}$ does not depend on the choice of the isometries.
We set
   $$\Phi(v^{\lambda}_{\phi,\psi}, v^\mu_{\xi,\eta}, v^\nu_{\chi,\zeta})=(\phi\otimes \xi\otimes \chi, q_{h_\lambda, h_\mu, h_\nu}\circ p_{h_\lambda+h_\mu+h_\nu}\psi\otimes \eta\otimes \zeta)_{p,3},$$
$$\Psi(v^{\lambda}_{\phi,\psi}, v^\mu_{\xi,\eta}, v^\nu_{\chi,\zeta})=(\phi\otimes \xi\otimes \chi, p_{h_\lambda+h_\mu+h_\nu}\circ q_{h_\lambda, h_\mu, h_\nu}\psi\otimes \eta\otimes \zeta)_{p,3}.$$
 The lemma then shows that
  $\Phi$, $\Psi$  are idempotent   but not selfadjoint  elements of the corresponding tensor power $^*$-algebra. The meaning of these elements and computation of the adjoints will soon be apparent.

    \medskip

  \noindent{\bf 9.4. Theorem } {\sl The idempotents $\Phi$, $\Psi$ satisfy the following relations,
  $$\Psi\Phi=\Delta\otimes 1\circ\Delta(I),\quad\quad \Phi\Psi=1\otimes\Delta\circ\Delta(I),\eqno(9.5)$$
  $$\Psi\Phi\Psi=\Psi,\quad\quad \Phi\Psi\Phi=\Phi\eqno(9.6)$$
   $$\Phi\Delta\otimes 1\circ\Delta(\omega)=1\otimes\Delta\circ\Delta(\omega)\Phi,\eqno(9.7)$$
  $$ \Psi1\otimes\Delta\circ\Delta(\omega)=\Delta\otimes 1\circ\Delta(\omega)\Psi.\eqno(9.8)$$}\medskip
  
  \noindent{\sl Proof}    To show $(9.7)$ and $(9.8)$
 we pass to the predual ${\mathcal C}$, and the corresponding   properties, which,  in Sweedler notation, respectively read as
 $$\Phi(a_1, b_1, c_1)(a_2b_2)c_2=a_1(b_1c_1)\Phi(a_2, b_2, c_2),\eqno(9.9)$$
 $$\Psi(a_1, b_1, c_1)a_2(b_2c_2)=(a_1b_1)c_1\Psi(a_2, b_2, c_2),\eqno(9.10)$$
 for $a$, $b$, $c\in{\mathcal C}$. 
 Set  $a=v^\lambda_{\phi,\psi}$, $b=v^\mu_{\xi, \eta}$, $c=v^\nu_{\chi,\zeta}$. 
  We write,  $$\Phi(v^\lambda_{\phi,\psi}, v^\mu_{\xi, \eta}, v^\nu_{\chi,\zeta})=(\phi\xi\chi, T\psi\eta\zeta)_{p,3}$$ and look for a solution $T$ of $(9.9)$.
Consider orthonormal bases 
 $\psi_r$, $\eta_s$, $\zeta_t$ of the appropriate truncated  tensor powers of $V$. The left hand side  becomes
 $$(\phi\xi\chi, T\psi_r\eta_s\zeta_t)_{p,3}[\psi_r^*\eta_s^*\zeta_t^* p_3\otimes p_3\psi\eta\zeta],$$
 where we have written $p_3$ for $p_{h_\lambda+h_\mu+h_\nu}$.
 But 
 $$(\phi\xi\chi, T\psi_r\eta_s\zeta_t)_{p,3}\psi_r^*\eta_s^*\zeta_t^*=((\psi_r\eta_s\zeta_t, T^+\phi\xi\chi)_{p,3}\psi_r\eta_s\zeta_t)^+=$$
 $$(p_1\otimes p_1'\otimes p_1''\circ  T^+\phi\xi\chi)^+=\phi^*\xi^*\chi^*\circ T\circ p_1\otimes p_1'\otimes p_1''$$
 where $p_1=p_{h_\lambda}$, $p_1'=p_{h_\mu}$, $p_1''=p_{h_\nu}$ and $+$ is the adjoint operator with respect to the triple product form.
 Notice that on each factor we are using Kirillov-Wenzl inner product, hence $p_1$, $p_1'$, $p_1''$ are selfadjoint. Thus the left hand side of $(9.9)$ is the class of
 $$(\phi^*\xi^*\chi^*\circ T\circ p_1\otimes p_1'\otimes p_1''\circ p_3)\otimes p_3\psi\eta\zeta.\eqno(9.11)$$
 On the other hand, $b_1c_1=\sum_{\gamma, i} [\xi^*\chi^*\circ S_{\gamma, i}\otimes S_{\gamma, i}^*\eta_s\zeta_t]$, hence the  right hand side gives
  $$(\phi^*\xi^*\chi^*\circ 1\otimes  S_{\gamma, i}\circ p_{h_\lambda+h_\gamma})\otimes (p_{h_\lambda+h_\gamma}\circ 1\otimes S_{\gamma, i}^*\psi_r\eta_s\zeta_t)(\psi_r\eta_s\zeta_t, T\psi\eta\zeta)_{p,3}=$$
$$(\phi^*\xi^*\chi^*\circ 1\otimes  S_{\gamma, i}\circ p_{h_\lambda+h_\gamma})\otimes (p_{h_\lambda+h_\gamma}\circ 1\otimes S_{\gamma, i}^*\circ p_1\otimes p_1'\otimes p_1''\circ T\psi\eta\zeta).\eqno(9.12)$$
Hence $(9.11)$ and $(9.12)$ coincide for $T=q_{h_\lambda, h_\mu, h_\nu}\circ p_{h_\lambda+h_\mu+h_\nu}$ thanks to Lemma 5.3.
In a similar way, $(9.10)$ leads to the equation 
$$(\phi^*\xi^*\chi^*\circ T'\circ p_1\otimes p_1'\otimes p_1''\circ 1\otimes S_{\gamma, i}\circ p_{h_\lambda+h_\gamma})\otimes(p_{h_\lambda+h_\gamma}\circ 1\otimes S_{\gamma, i}^*\psi\eta\zeta)=$$
$$\phi^*\xi^*\chi^*\circ p_3\otimes(p_3\circ p_1\otimes p_1'\otimes p_1''\circ T'\psi\eta\zeta),$$
which is solved for $T'=p_{h_\lambda+h_\mu+h_\nu}\circ q_{h_\lambda, h_\mu, h_\nu}$.
The previous lemma shows that $\Psi\Phi$ and $\Phi\Psi$ are the functionals induced respectively by
  $$p_3\circ q\circ p_1\otimes p_1'\otimes p_1''\circ q\circ p_3=p_3,$$  
  $$q\circ p_3\circ p_1\otimes p_1'\otimes p_1''\circ p_3\circ q=  q\circ p_3\circ q,$$ and it is easy to check that these coincide with the right hand sides of $(9.5)$. One similarly shows $(9.6)$.

\medskip

 \noindent{\bf 9.5. Theorem} {\sl The element $\Phi\in\widehat{\mathcal C}^{\otimes 3}$ is a counital  $3$-cocycle:
 $$1\otimes 1\otimes\Delta(\Phi)\Delta\otimes 1\otimes 1(\Phi)=I\otimes\Phi1\otimes\Delta\otimes 1(\Phi)\Phi\otimes I,\eqno(9.13)$$
 $$\varepsilon\otimes 1\otimes 1(\Phi)=1\otimes \varepsilon\otimes 1(\Phi)=1\otimes 1\otimes\varepsilon(\Phi)=P.$$ }\medskip
 
 \noindent{\sl Proof} We compute the elements   $T_1$, $T_2$, $T_3\in\widehat{\mathcal C}^{\otimes 4}$  
 corresponding to  $\Delta\otimes 1\otimes 1(\Phi)$, $1\otimes\Delta\otimes 1(\Phi)$, and $1\otimes 1\otimes\Delta(\Phi)$.  
 We first write  
 $$\Phi(v^\lambda_{\phi,\psi},  v^\mu_{\xi,\eta}, v^\nu_{\chi,\zeta})=\phi^*\xi^*\chi^*(1_{h_\lambda}\otimes S_{\gamma, i}\circ p_{h_\lambda+h_\gamma}\circ 1_{h_\lambda}\otimes S_{\gamma, i}^*\circ p_{h_\lambda+h_\mu+h_\nu}\psi\eta\zeta)$$
 with $S_{\gamma,i}\in(V_\gamma, V^{\underline{\otimes}(h_\mu+h_\nu)})$.
 For $T_1$ we need to replace the first variable by a product in ${\mathcal C}$. Explicitly,
 for $S_{\sigma, j}\in(V_\sigma, V^{\underline{\otimes} (h_\lambda+h_\tau)})$,
 $$\Delta\otimes 1\otimes 1(\Phi)(v^\lambda_{\phi,\psi}, v^\tau_{\alpha,\beta}, v^\mu_{\xi,\eta}, v^\nu_{\chi,\zeta})=\Phi((v^\lambda_{\phi,\psi}v^\tau_{\alpha,\beta}), v^\mu_{\xi,\eta}, v^\nu_{\chi,\zeta})=$$
 $$\phi^*\alpha^*\xi^*\chi^*(S_{\sigma, j} \otimes 1_{h_\mu+h_\nu}\circ1_{h_\sigma}\otimes S_{\gamma, i}\circ p_{h_\sigma+h_\gamma}\circ 1_{h_\sigma}\otimes S_{\gamma, i}^*\circ p_{h_\sigma+h_\mu+h_\nu}\circ S_{\sigma, j}^* \otimes 1_{h_\mu+h_\nu}\psi\beta\eta\zeta).$$
Taking into account Lemma 5.1 a),
$$T_1=1_{h_\lambda+h_\tau}\otimes S_{\gamma, i}\circ p_{h_\lambda+h_\tau+h_\gamma}\circ 1_{h_\lambda+h_\tau}\otimes S_{\gamma, i}^*\circ p_{h_\lambda+h_\tau+h_\mu+h_\nu}.$$
One also finds
$$T_2=1_{h_\lambda}\otimes [S_{\sigma', j'}\otimes 1_{h_\nu}\circ S_{\gamma', i}]\circ   p_{h_\lambda+h_{\gamma'}}\circ 1_{h_\lambda}\otimes S_{\gamma', i'}^*\circ p_{h_\lambda+h_{\sigma'}+h_\nu}\circ 1_{h_\lambda}\otimes S_{\sigma', j'}^*\otimes 1_{h_\nu},$$
where $S_{\gamma', i'}\in(V_{\gamma'}, V^{\underline{\otimes}(h_{\sigma'}+h_\nu)})$, $S_{\sigma', j'}\in(V_{\sigma'}, V^{\underline{\otimes}(h_{\tau}+h_\mu)})$, and
$$T_3=1_{h_\lambda+h_\tau}\otimes S_{\sigma'', j''}\circ 1_{h_\lambda}\otimes S_{\gamma'', i''}\circ p_{h_\lambda+h_{\gamma''}}\circ 1_{h_\lambda}\otimes S_{\gamma'', i''}^*\circ p_{h_\lambda+h_\tau+h_{\sigma''}}\circ 1_{h_\lambda+h_\tau}\otimes S_{\sigma'', j''}^*,$$                
with $S_{\gamma'', i''}\in(V_{\gamma''}, V^{\underline{\otimes}(h_{\tau}+h_{\sigma''})})$, $S_{\sigma'', j''}\in(V_{\sigma''}, V^{\underline{\otimes}(h_{\mu}+
h_\nu)})$. Verification of the cocycle condition $(9.13)$ involves rather long formulas. We shall thus ease notation dropping most indices and writing for example $p_2$ for expressions of the form $p_{h_\tau+h_\gamma}$.

 We   compare  the matrices corresponding to both  sides of $(9.13)$ taking into account Lemma 9.3,  Lemma 5.3 and Lemma 5.1 a) again. For the  left hand side we have
 $$T_3\circ p_{1}\otimes p_{1}\otimes p_{1}\otimes p_{1}\circ T_1=T_3T_1=$$
 $$1\otimes 1\otimes S_{\sigma''}\circ1\otimes S_{\gamma''}\circ p_2\circ1\otimes S_{\gamma''}^*\circ p_3\circ 1\otimes 1\otimes S_{\sigma''}^*\circ 1\otimes 1\otimes S_\gamma\circ p_3\circ 1\otimes 1\otimes S_\gamma^*\circ p_4=$$
 $$1\otimes 1\otimes S_{\sigma''}\circ1\otimes S_{\gamma''}\circ p_2\circ1\otimes S_{\gamma''}^* \circ 1\otimes 1\otimes S_{\sigma''}^*\circ p_4,$$
 while the right hand side gives, if $T$ is the matrix defining $\Phi$,
 $$(I\otimes T)\circ p_1^{\otimes 4}\circ T_2\circ p_1^{\otimes 4}\circ(T\otimes I)=$$
 $$1\otimes 1\otimes S_\gamma\circ 1\otimes p_2\circ 1\otimes 1\otimes S_\gamma^*\circ 1\otimes p_3\circ p_1^{\otimes 4}\circ 1\otimes (S_{\sigma'}\otimes 1\circ S_{\gamma'})\circ p_2\circ$$
 $$ 1\otimes S_{\gamma'}^*\circ p_3\circ  1\otimes S_{\sigma'}^*\otimes 1\circ p_1^{\otimes 4} \circ 1\otimes S_{\gamma'''}\otimes 1\circ p_2\otimes 1\circ 1\otimes S_{\gamma'''}^*\otimes 1\circ p_3\otimes 1=$$
 $$1\otimes 1\otimes S_\gamma\circ 1\otimes p_2\circ 1\otimes 1\otimes S_\gamma^*  \circ 1\otimes (S_{\sigma'}\otimes 1\circ S_{\gamma'})\circ p_2\circ$$
 $$ 1\otimes S_{\gamma'}^*\circ [p_3 \circ (p_2\circ 1\otimes S_{\sigma'}^*\circ p_3)\otimes 1] \circ 1\otimes 1\otimes 1\otimes p_1=$$
  in the last  equality  we have   deleted the extra idempotents $1\otimes p_3$, the first     $p_1^{\otimes 4}$ altogether, and the first three factors $p_1^{\otimes 3}$ of the second copy of the same idempotent, and moved the fourth to the far right. This has also allowed to use $S_{\sigma'}^*\circ S_{\gamma'''}=\delta_{\sigma', \gamma'''}$. Next write the term in square brackets as
 $  (p_2\circ 1\otimes S_{\sigma'}^*\circ p_3)\otimes 1\circ p_4$ and after elimination of  $p_2$ and $p_3$, the above term becomes
  $$1\otimes 1\otimes S_\gamma \circ  1\otimes (p_2\circ  1\otimes S_\gamma^*\circ S_{\sigma'}\otimes 1\circ S_{\gamma'})\circ p_2\circ$$
 $$ 1\otimes S_{\gamma'}^*\circ 
1\otimes S_{\sigma'}^* \otimes 1\circ p_4 \circ 1\otimes 1\otimes 1\otimes p_1=$$
   $$1\otimes 1\otimes S_\gamma \circ  1\otimes S'_{\gamma'}\circ p_2\circ1\otimes {S'_{\gamma'}}^*\circ 1\otimes 1\otimes S_\gamma^*\circ p_4 \circ 1\otimes 1\otimes 1\otimes p_1$$
where $S'_{\gamma'}=p_2\circ  1\otimes S_\gamma^*\circ S_{\sigma'}\otimes 1\circ S_{\gamma'}$ is another orthonormal system of isometries.
  Hence the two matrices induce the same functional.
 \medskip
 
 \noindent{\bf 9.6. Remarks} 
 a) Notice that not only
 the associators $\Phi$, $\Psi$ are   quasi-invertible elements of a rather special kind, in that they are idempotent, but also there is   a simple relation between    iterated coproducts associated so
 $$\Delta_{\rm left}^{(n)}:=\Delta\otimes 1_{n-1}\circ\dots\circ\Delta\otimes 1\circ \Delta$$ and
arbitrary  iterated coproducts $\Delta^{(n)}$ of order $n$, i.e.   coproducts that  can be obtained
 as compositions ${\widehat{\mathcal C}}\to {\widehat{\mathcal C}}\otimes {\widehat{\mathcal C}}\to\dots\to {\widehat{\mathcal C}}^{\otimes n+1}$ where the connecting maps  
 ${\widehat{\mathcal C}}^{\otimes j}\to{\widehat{\mathcal C}}^{\otimes j+1}$ can be  an arbitrary translates $1_r\otimes\Delta\otimes 1_{j-r-1}$ of $\Delta$. Set $P_n=\Delta_{\rm left}^{(n)}(I)$, hence in particular  $P_1=P$, $P_2=\Psi\Phi$.
 Then
 $$\Delta_{\rm left}^{(n)}(\omega)=
 P_n\Delta^{(n)}(\omega)P_n$$
 for all possible choices of $\Delta^{(n)}$.
 This relation can be derived with computations similar to those  appearing in the proof  of Theorem 9.4, and also in the following Prop. 9.7.
 b) In a similar way, the associators are explicitly related to non-unitality of $\Delta$ by $\Phi=1\otimes\Delta\circ\Delta(I)\Delta\otimes1\circ\Delta(I)$, and similarly for $\Psi$.
\bigskip

 \noindent{\it d) $R$-matrices}\bigskip

\noindent We may introduce quasi-invertible $R$-matrices ${\mathcal R}, {\mathcal R}_1\in\widehat{\mathcal C}\otimes\widehat{\mathcal C}$ appealing to the braided symmetry of the quotient category ${\mathcal F}_\ell$.
 More precisely, set
  $${\mathcal R}(v^\lambda_{\phi,\psi}, v^\mu_{\xi,\eta})=(\phi\otimes \xi, \Sigma_{k,h}p_{k+h}\Sigma_{h,k}Rp_{h+k}\psi\otimes \eta)_{p,2},$$
    $${\mathcal R}_1(v^\lambda_{\phi,\psi}, v^\mu_{\xi,\eta})=(\phi\otimes \xi, p_{k+h}R^{-1}\Sigma_{k, h}p_{k+h}\Sigma_{h,k}\psi\otimes\eta)_{p,2},$$
  with $\Sigma_{h,k}: V^{\otimes h}\otimes V^{\otimes k}\to V^{\otimes k}\otimes V^{\otimes h}$ the transposition operator.\medskip
  
  \noindent{\bf 9.7. Proposition} {\sl The following relations hold,
  $${\mathcal R}\Delta(I)={\mathcal R}= \Delta^{\text{op}}(I){\mathcal R},\quad\quad {\mathcal R}_1\Delta^{\text{op}}(I)={\mathcal R}_1=\Delta(I){\mathcal R}_1\eqno(9.14)$$
   $${\mathcal R}\Delta(\omega){\mathcal R}_1= \Delta^{\text{op}}(\omega),\quad\quad{\mathcal R}_1 \Delta^{\text{op}}(\omega){\mathcal R}=\Delta(\omega),\eqno(9.15)$$
  where $ \Delta^{\text{op}}$ is the opposite coproduct of $\widehat{{\mathcal C}(G,\ell)}$.}\medskip
  
  \noindent{\sl Proof} We check $(9.15)_1$. Let $T=(T_n)$ and $T'=(T'_n)$ define functionals $\omega_T$, $\omega_{T'}$ on $\widehat{\mathcal C}\otimes\widehat{\mathcal C}$.
  Then a relation $\omega_T\Delta(\tau)\omega_{T'}= \Delta^{\text{op}}(\tau)$, for all $\tau\in\widehat{\mathcal C}$ is equivalent to $\omega_T(a_1, b_1)a_2b_2\omega_{T'}(a_3, b_3)=ba$
  for all $a$, $b\in{\mathcal C}$. For $a=v^\lambda_{\phi,\psi}$, $b=v^\mu_{\xi,\eta}$, this relation reads
  $$(\phi^*\xi^*\circ T\circ p_h\otimes p_k\circ p_{h+k})\otimes(p_{h+k}\circ p_h\otimes p_k\circ T'\psi\eta)=\xi^*\phi^*\circ p_{k+h}\otimes p_{k+h}\eta\psi.$$   as class elements of ${\mathcal C}$.
  For $T=\Sigma_{k,h}p_{k+h}\Sigma_{h,k}Rp_{h+k}$, $T'=p_{k+h}R^{-1}\Sigma_{k, h}p_{k+h}\Sigma_{h,k}$ the left hand side becomes
  $$(\xi^*\phi^*p_{k+h}\Sigma_{h,k}Rp_{h+k})\otimes(p_{h+k}(\Sigma_{h,k}R)^{-1}p_{k+h}\eta\psi).$$
  But $p_{k+h}\Sigma_{h,k}R p_{h+k}$ and 
  $p_{h+k}(\Sigma_{h,k}R)^{-1} p_{k+h}$ are braiding arrows and  inverses of one another in  the quotient category, hence they cancel out.\medskip
   
 We next briefly discuss  a relation between the original quantum group $U_q({\mathfrak g})$ and $\widehat{{\mathcal C}(G,\ell)}$. There is a natural map
 $$\pi: U_q({\mathfrak g})\to \widehat{{\mathcal C}(G,\ell)}$$
 taking an element $a\in U_q({\mathfrak g})$ to the functional
 $$\pi(a)(v^\lambda_{\phi,\psi})=(\phi, a\psi), \quad v^\lambda_{\phi,\psi}\in M_\lambda,\quad \lambda\in \Lambda_\ell.\eqno(9.16)$$
 \medskip
 
 \noindent{\bf 9.8. Proposition} {\sl $\pi$ is a surjective homomorphism of  $^*$--algebras satisfying
 $$P\pi\otimes\pi(\Delta(a))=\Delta(\pi(a))=\pi\otimes\pi(\Delta(a))P, \quad a\in U_q({\mathfrak g}).\eqno(9.17)$$
 
 }\medskip
 
 \noindent{\sl Proof} The map $\pi$ corresponds to $\oplus_{\lambda\in\Lambda_\ell} V_\lambda$ under the identification
 $\widehat{\mathcal C}\simeq\oplus_{\lambda\in\Lambda_\ell}{\mathcal L}(V_\lambda)$, and this shows the homomorphism property.  Surjectivity is a consequence of irreducibility of the $V_\lambda$.
 The last statement is a straightforward computation.
  \medskip
  
  We use the above link to resume a compatibility relation between coproduct and adjoint of $\widehat{\mathcal C}$.\medskip
  
  \noindent{\bf 9.9. Proposition} {\sl We have
  $$\Delta(\omega)^*=\pi\otimes\pi(\overline{R})\Delta(\omega^*)\pi\otimes\pi(\overline{R})^{-1},\quad\quad\omega\in\widehat{\mathcal C}.\eqno(9.18)$$}\medskip
  
  \noindent{\sl Proof} We first show the relation for $\omega=I$. Recalling  that $\Delta(I)=P$, we have
  $$P^*(v^\lambda_{\phi,\psi}, v^\mu_{\xi,\eta})=\overline{P((v^\lambda_{\phi,\psi})^*, (v^\mu_{\xi,\eta})^*)}=$$
  $$\overline{P(v^\lambda_{\psi, \phi}, v^\mu_{\eta, \xi})}=\overline{(\psi\eta, p_{h+k}\phi\xi)}_{p,2}=$$
  $$(p_{h+k}\phi\xi, \psi\eta)_{p,2}=(\phi\xi, \overline{R}p_{h+k}\overline{R}^{-1}\phi\xi)_{p,2}=$$
  $$\pi\otimes\pi(\overline{R})P\pi\otimes\pi(\overline{R})^{-1}(v^\lambda_{\phi,\psi}, v^\mu_{\xi,\eta}).$$
 Taking into account the relations $ \Delta^{\text{op}}(a)=\overline{R}\Delta(a)\overline{R}^{-1}$ and $\Delta(a)^*= \Delta^{\text{op}}(a^*)$ in $U_q({\mathfrak g})$ for $a\in U_q({\mathfrak g})$,
 it follows that
  $$\Delta(\pi(a))^*=(\pi\otimes\pi(\Delta(a)) P)^*=$$
  $$\pi\otimes\pi(\overline{R})P\pi\otimes\pi(\overline{R})^{-1}\pi\otimes\pi(\overline{R}\Delta(a^*)\overline{R}^{-1})=$$
  $$\pi\otimes\pi(\overline{R})\Delta(\pi(a^*))\pi\otimes\pi(\overline{R})^{-1}.$$\bigskip
  
  \noindent{\it e)  Comparing the groupoids associated to different sections  }\bigskip
  
\noindent  Let $s: {\mathcal C}(G,\ell)\to{\mathcal D}(V,\ell)$ and  $s': {\mathcal C}(G,\ell)\to{\mathcal D}(V,\ell)$ be different sections of the quotient map ${\mathcal D}(V,\ell)\to{\mathcal C}(G,\ell)$.
  Correspondingly, we upgrade the algebra $\widehat{{\mathcal C}(G,\ell)}$ to quasi-coassociative quantum groupoids in two different ways, 
  $$(\Delta,   \Phi, \Psi, {\mathcal R}, {\mathcal R}_1) \quad \text{ and }\quad  (\Delta', \Phi', \Psi', {\mathcal R}', {\mathcal R}'_1),$$
  (we have noticed before that $S$ is independent of the section.)
  We claim that these are related by a twisting procedure in the sense of Drinfeld induced, again,  by  special  quasi-invertible elements $\widehat{\mathcal C}\otimes\widehat{\mathcal C}$.
    We sketch  the intertwining relation between the coproducts, but we refrain from giving details for the remaining   relations.
 
 \medskip
  
  \noindent{\bf 9.10. Proposition} {\sl Set $P=\Delta(I)$, $P'=\Delta'(I)$.
Then for  $\omega\in\widehat{\mathcal C}$,
  $$P\Delta'(\omega)P=\Delta(\omega),\quad\quad P'\Delta(\omega)P'=\Delta'(\omega).$$}\medskip
  
  \noindent{\sl Proof}
  To show for example the first relation, we first write $P$ in a different way, as follows. Let $s$ and $s'$ be defined by two choices of irreducible summands $V_\lambda$ and $V'_\lambda$ both of highest weight $\lambda$    in $V^{\underline{\otimes} h_\lambda}$ and 
  $V^{\underline{\otimes} h'_\lambda}$, respectively.
For each $\lambda\in\Lambda_\ell$ there is a unitary intertwiner $U_\lambda: V'_\lambda\to V_\lambda$, unique up to a scalar multiple 
  by an element of ${\mathbb T}$. We associate the element  $F\in\widehat{\mathcal C}\otimes\widehat{\mathcal C}$ defined by
  $$F(v^\lambda_{\phi', \psi'}, v^\mu_{\xi',\eta'})=(U_\lambda\phi'\otimes U_\mu\xi',  \, p_{h_\lambda+h_\mu}\circ U_\lambda\psi'\otimes U_\mu\eta')_{p,2},$$
  for $\phi'$, $\psi'\in V'_\lambda$, $\xi'$, $\eta'\in V'_\mu$.
  On one hand we may use the defining identifications in ${\mathcal C}(G,\ell)$ and  write, for $\phi$, $\psi\in V_\lambda$, %%$\xi$, $\eta\in V_\mu$,
  $v^\lambda_{\phi, \psi}=v^\lambda_{U_\lambda^*\phi, U_\lambda^*\psi}$, and this shows that
  $F=P.$
  On the other, we can use the computations of the first part of the proof of Prop. 9.7 and derive the relation $F\Delta'(\omega)F=\Delta(\omega)$.
\bigskip

\noindent{\it f) The tensor $^*$--equivalence ${\mathcal F}_\ell\to \text{Rep}_V \widehat{{\mathcal C}(G,\ell)}.$
  }\bigskip

\noindent In this subsection we show   that if  ${\mathcal C}(G,\ell)$ is cosemisimple then the smallest full tensor subcategory $ \text{Rep}_V \widehat{{\mathcal C}(G,\ell)}$
of the representation category of $\widehat{{\mathcal C}(G,\ell)}$ containing the fundamental representation is a tensor $C^*$--category equivalent 
to the fusion category ${\mathcal F}_\ell$.\medskip

  There is a known way to associate a tensor category to any weak quasi-Hopf algebra \cite{MS}, that we briefly recall for $\widehat{\mathcal C}$. Let $\text{Rep}(\widehat{{\mathcal C}})$ be the category of unital representations of   $\widehat{\mathcal C}$ on f.d. vector spaces.   We have already noticed that $P=\Delta(I)$ is an idempotent of $ \widehat{{\mathcal C}}\otimes  \widehat{{\mathcal C}}$ commuting with all elements in the image of $\Delta$. Hence if $u$ and $v$
 are objects of $\text{Rep}(\widehat{\mathcal C})$   then $\omega\in\widehat{{\mathcal C}}\to u\otimes v\circ\Delta(\omega)$ is   a representation of the same algebra on the tensor product space with support the image of the idempotent $u\otimes v(P)$.
 We   define the tensor product object, still denoted $u\otimes v$, to be  the unital representation of $\widehat{\mathcal C}$ obtained restricting the operators $u\otimes v\circ\Delta(\omega)$ to that subspace.
 The   map $\Phi_{u,v,w}$ given by the restriction of $u\otimes v\otimes w(\Phi)$ to the space of $(u\otimes v)\otimes w$ is an invertible intertwiner from $(u\otimes v)\otimes w$ to $u\otimes(v\otimes w)$, by Theorem 
 9.4, which, together with Theorem  9.5, guarantees that
   $\text{Rep}(\widehat{\mathcal C})$ becomes a tensor category with these associativity morphisms. We next consider the full subcategory $\text{Rep}_h(\widehat{{\mathcal C}})$ of $\text{Rep}(\widehat{{\mathcal C}})$ with objects $^*$--representations 
   of $\widehat{\mathcal C}$ on non-degenerate Hermitian spaces. Following a procedure similar to that of the non-semisimple case, if $u$ and $v$ are two  objects of $\text{Rep}_h(\widehat{{\mathcal C}})$ then
   $u\otimes v$ is still a $^*$--homomorphism, and hence an object of $\text{Rep}_h(\widehat{{\mathcal C}})$, provided we endow the corresponding vector space with the restriction    of Wenzl's Hermitian form induced by $\overline{R}$. This can be seen with Prop. 9.9. Hence  $\text{Rep}_h(\widehat{{\mathcal C}})$  is a tensor $^*$--category. 
  Now if $u$ and $v$ are Hilbert space $^*$--representations then we can infer that $u\otimes v$ is still a Hilbert space representation provided  $u\otimes v(P)$ is a positive operator with respect to Wenzl's form, and this may not always be the case. But there are sufficient objects of $\text{Rep}_h(\widehat{{\mathcal C}})$ for which the positivity condition holds.
  Consider the map $\hat{V}: \widehat{\mathcal C}\to{\mathcal B}(V)$ defined by
  $$(\xi, \hat{V}(\omega)\eta)=\omega([\xi^*\otimes\eta]),\quad\quad\omega\in \widehat{\mathcal C}, \quad \xi,\eta\in V.$$  It is easy to see that this  is a unital representation, and in fact a Hilbert space $^*$--representation. Let $\text{Rep}_V(\widehat{\mathcal C})$ be the
  the smallest full tensor subcategory of $\text{Rep}_h(\widehat{\mathcal C})$ containing $\hat{V}$. The objects of $\text{Rep}_V(\widehat{\mathcal C})$ are the trivial representation, $\hat{V}$ and representations of the form
$\omega\to \hat{V}\otimes \dots\otimes \hat{V}\circ\Delta^{(n)}(\omega)$, $n\geq 1$, where $\Delta^{(n)}$ denotes an iteration of  $n$ translates of the coproduct   as in Remark 9.6. 
\medskip
   
   \noindent{\bf 9.11. Theorem} {\sl    $\text{Rep}_V(\widehat{\mathcal C})$ is a tensor $C^*$--category.}\medskip
   
   \noindent{\sl Proof} The representation $\omega\to \hat{V}\otimes \dots\otimes \hat{V}\circ\Delta^{(n)}(\omega)$ acts on the subspace of $V^{\otimes n}$ corresponding to the idempotent
  $\hat{V}\otimes \dots\otimes \hat{V}\circ\Delta^{(n)}(I)$. We need to show that Wenzl's inner product of $V^{\otimes n}$ is positive on that subspace. We limit ourselves to computing these idempotents in the example  
   $\Delta^{(2)}=1\otimes\Delta\circ\Delta$. In Sweedler notation, $\Delta(I)=\omega_1\otimes\omega_2$, $\Delta^{(2)}(I)=\omega_1\otimes\omega_{2,1}\otimes\omega_{2,2}$. Hence
   $$(\phi_1\phi_2\phi_3, \hat{V}\otimes\hat{V}\otimes\hat{V}\circ\Delta^{(2)}(I)\psi_1\psi_2\psi_3)_p=$$
   $$(\phi_1,\hat{V}(\omega_{1})\psi_1)(\phi_2,\hat{V}(\omega_{2,1})\psi_2)(\phi_3,\hat{V}(\omega_{2,2})\psi_3)=$$
   $$\omega_1([\phi_1^*\otimes\psi_1])\omega_{2,1}([\phi_{2}^*\otimes\psi_2])\omega_{2,2}([\phi_3^*\otimes\psi_3])=$$
   $$\omega_1([\phi_1^*\otimes\psi_1])\omega_2([\phi_2^*\phi_3^*\circ p_2\otimes p_2\psi_2\psi_3])=$$
   $$\omega_1([\phi_1^*\otimes\psi_1])\omega_2([\phi_2^*\phi_3^*\circ S_\alpha\otimes S_\alpha^*\psi_2\psi_3])=$$
   $$\varepsilon([(\phi_1^*\phi_2^*\phi_3^*\circ 1\otimes S_\alpha\circ p_{1+n_\alpha})\otimes(p_{1+n_\alpha}\circ 1\otimes S_\alpha^*\psi_1\psi_2\psi_3)])=$$
   $$(\phi_1\phi_2\phi_3, 1\otimes S_\alpha\circ p_{1+n_\alpha}\circ 1\otimes S_\alpha^*\psi_1\psi_2\psi_3)_p,$$
   where $S_\alpha\in(V_{\lambda_\alpha}, V^{\underline{\otimes} 2})$ are the isometries describing  decomposition of $p_2$ into irreducibles with respect to  the chosen section $s$, and $V_{\lambda_\alpha}$ is a summand of $V^{\otimes h_{\lambda_\alpha}}$, as before. Hence  $\hat{V}\otimes\hat{V}\otimes\hat{V}\circ\Delta^{(2)}(I)=\sum 1\otimes S_\alpha\circ p_{1+h_{\lambda_\alpha}}\circ 1\otimes S_\alpha^*$, and this is a positive operator w.r.t. Wenzl's inner product.\medskip
   
   We finally discuss the relation between the original (strict) fusion category ${\mathcal F}_\ell$ and $\text{Rep}_V(\widehat{\mathcal C})$.  
Let $X=V^{\underline{\otimes} n}$ be regarded as an object of ${\mathcal F}_\ell$, meaning that $X$ is  the truncated $U_q({\mathfrak g})$--submodule of $V^{\otimes n}$, endowed with  Hilbert space structure, cf. Theorem 5.4.
We associate to $X$ the map $$\hat{X}: \widehat{{\mathcal C}(G,\ell)}\to{\mathcal B}(X),\quad\quad (\phi, \hat{X}(\omega)\psi)=\omega([\phi^*\otimes\psi]),$$
for $\phi$, $\psi\in V^{\underline{\otimes} n}$, $\omega\in \widehat{{\mathcal C}(G,\ell)}$. This formula extends the previously introduced $\hat{V}$ to all objects of ${\mathcal F}_\ell$, and, as before, 
it is easily seen that $\hat{X}$ is a unital $^*$--representation of $\widehat{{\mathcal C}(G,\ell)}$ on $V^{\underline{\otimes} n}$. One has $\hat{X}=\hat{V}\otimes\dots\otimes\hat{V}\circ\Delta^{(n-1)}_{\text{left}}$, hence $\hat{X}$ is an object of $\text{Rep}_V(\widehat{\mathcal C})$.
 \medskip
 
 \noindent{\bf  9.12. Theorem} {\sl The map   $X\to\hat{X}$ and acting identically on arrows is a tensor  $^*$--equivalence $${\mathcal E}: {\mathcal F}_\ell\to\text{Rep}_V(\widehat{\mathcal C}).$$}\medskip
 
 \noindent{\sl Proof}
 An arrow $T\in (X, X')$ in ${\mathcal F}_\ell$ is an intertwiner of the corresponding modules of $U_q({\mathfrak g})$. But  it also lies in the arrow space $(\hat{X}, \hat{X'})$ of $\text{Rep} \widehat{{\mathcal C}(G,\ell)}$,
 as 
 $$(\phi, T\hat{X}(\omega)\psi)=(T^*\phi, \hat{X}(\omega)\psi)=\omega([(T^*\phi)^*\otimes\psi])=$$
 $$\omega([\phi^*\circ T\otimes\psi])=\omega([\phi^*\otimes T\psi])=(\phi,\hat{X'}(\omega)T\psi),$$   
 hence ${\mathcal E}$ is a functor between the stated categories, and it 
 is easy to see that it is actually a $^*$--functor.  We verify that ${\mathcal E}$ has full image. Routine computations show that      formula $(9.16)$ defining  the surjection $\pi: U_q({\mathfrak g})\to\widehat{{\mathcal C}(G, \ell)}$
 extends to all representations, in that $\pi(a)([\phi^*\otimes \psi])=(\phi, a\psi)$ holds for all $X$. Now if $T\in(\hat{X}, \hat{X'})$ then $T\hat{X}(\pi(a))=\hat{X'}(\pi(a))T$, and expliciting this relation yields $T$ as an intertwiner of the corresponding $U_q({\mathfrak g})$-representations. Since it is obviously essentially surjective, we are left to construct a tensor structure for ${\mathcal E}$ making it into a tensor functor. Explicitly, since ${\mathcal F}_\ell$ is strict while $\text{Rep}_V(\widehat{\mathcal C})$ is not,
  we look for   natural isomorphisms ${\mathcal E}_{X, Y}\in(\hat{X}\otimes\hat{Y}, \widehat{X\otimes Y})$ such that 
 $${\mathcal E}_{W, X\otimes Y}\circ1_{\hat{W}}\otimes{\mathcal E}_{X,Y}\circ\Phi_{\hat{W},\hat{X},\hat{Y}}={\mathcal E}_{W\otimes X, Y} \circ{\mathcal E}_{W, X}\otimes 1_{\hat{Y}}.\eqno(9.19)$$
 For 
 $W=V^{\underline{\otimes} m}$,
 $X=V^{\underline{\otimes} n}$, $Y=V^{\underline{\otimes} r}$, we write
 $${\mathcal E}_{X,Y}=\hat{V}^{\otimes(n+r)}(Q_{n,r})$$ and we are reduced to look for 
  a quasi-invertible $Q_{n,r}\in\widehat{\mathcal C}^{\otimes n+r}$ satisfying the intertwining relation
 $$Q_{n,r}\Delta^{(n-1)}_{\text{left}}\otimes \Delta^{(r-1)}_{\text{left}}(\Delta(\omega))=\Delta^{(n+r-1)}_{\text{left}}(\omega)Q_{n,r}.$$
 It is solved by $$Q_{n,r}=\Delta^{(n+r-1)}_{\text{left}}(I)\Delta^{(n-1)}_{\text{left}}\otimes \Delta^{(r-1)}_{\text{left}}(\Delta(I))$$
 and naturality easily follows. 
 The tensorial structure $(9.19)$ corresponds to
$$Q_{m, n+r} I_m\otimes Q_{n,r}\Delta^{(m-1)}_{\text{left}}\otimes \Delta^{(n-1)}_{\text{left}}\otimes\Delta^{(r-1)}_{\text{left}}(\Phi)=Q_{m+n, r}Q_{m, n}\otimes I_r.\eqno(9.20)$$ 
Computations analogous to those of Theorem 9.5 give, in the same short notation, $Q_{n,r}=\omega_{T_{n,r}}$ with 
   $$T_{n,r}=
 p_{n+r}\circ 1_n\otimes S_\gamma\circ p_{n+1}\circ 1_{n}\otimes S_\gamma^*$$ and this expression shows that $Q_{n,r}$ is quasi-invertible.
 Furthermore $\Delta^{(m-1)}_{\text{left}}\otimes\Delta^{(n-1)}_{\text{left}}\otimes\Delta^{(r-1)}_{\text{left}}(\Phi)$ corresponds to
 $$S_{\sigma}\otimes S_{\sigma'}\otimes S_{\sigma''}\circ 1\otimes S_\gamma\circ p_2\circ 1\otimes S_{\gamma}^*\circ p_3\circ S_{\sigma}^*\otimes S_{\sigma'}^*\otimes S_{\sigma''}^*$$
 and we next check validity of the desired relation $(9.20)$.
At the left hand side we obtain
\[
\begin{split}
p_{m+n+r}\circ 1_m\otimes S_{\gamma}\circ p_{m+h_{\gamma}}\circ 1_m\otimes S_{\gamma}^*\circ p_m\otimes 1_{n+r}\circ 1_m\otimes p_{n+r}\circ \\
\circ 1_{m+n}\otimes S_{\delta}\circ 1_m\otimes p_{n+h_{\delta}}\circ 1_{m+n}\otimes S_{\delta}^*\circ S_{\sigma}\otimes S_{\sigma'}\otimes S_{\sigma''}\circ \\
\circ 1_m\otimes S_{\alpha}\circ p_{h_{\sigma}+h_{\alpha}}\circ 1\otimes S_{\alpha}^*\circ p_{h_{\sigma}+h_{\sigma'}+h_{\sigma''}}\circ S_{\sigma}^*\otimes S_{\sigma'}^*\otimes S_{\sigma''}^* 
\end{split}
\]
which equals $$p_{m+n+r}\circ 1_m\otimes S_{\sigma'}\otimes S_{\sigma''}\circ p_{m+h_{\sigma'}+h_{\sigma''}}\circ 1_m\otimes S_{\sigma'}^*\otimes S_{\sigma''}^*$$
by repeated use of  Lemma 5.1 and Lemma 5.3. The right hand side becomes
\[
\begin{split}
p_{m+n+r}\circ 1_{m+n}\otimes S_{\sigma''}\circ p_{m+n+h_{\sigma''}}\circ 1_{m+n}\otimes S_{\sigma''}^*\circ p_{m+n}\otimes p_r\circ \\
\circ 1_m\otimes S_{\sigma'}\otimes 1_r\circ p_{m+h_{\sigma'}}\otimes 1_r\circ 1_m\otimes S_{\sigma'}^*\otimes 1_r
\end{split}
\]
which in turn equals $$p_{m+n+r}\circ 1_{m+n}\otimes S_{\sigma''}\circ[p_{m+n+h_{\sigma''}}\circ(p_{m+n}\circ 1_m\otimes S_{\sigma'}\circ p_{m+h_{\sigma'}})\otimes 1_{h_{\sigma''}}]\circ 1_m\otimes S_{\sigma'}^*\otimes S_{\sigma''}^*.$$
It is easy to see that the expression in the square brackets can be rewritten as
$$(p_{m+n}\circ 1_m\otimes S_{\sigma'}\circ p_{m+h_{\sigma'}})\otimes 1_{h_{\sigma''}}\circ p_{m+h_{\sigma'}+h_{\sigma''}}.$$
Substituting it in the right hand side formula and using the usual properties of projections we get the desired identity.

 \bigskip

  \section{Cosemisimplicity in type $A$ case}

  As already mentioned, in this section we formulate a sufficient condition for cosemisimplicity of ${\mathcal C}(G,\ell)$, that we next verify for $G={\rm SU}(N)$. Our condition appeals to the existence of a Haar 
  functional  on 
  ${\mathcal C}(G,\ell)$, but more is needed. In particular, we shall use the associative filtration $\widetilde{\mathcal C}_k$ constructed in  Sect. 8.  We should remark that   space limit has prevented us from      studying   cosemisimplicity of ${\mathcal C}(G,\ell)$ for other Lie types. A positive answer would obviously suffice to extend our main result to other types  provided   ${\mathfrak g}\neq E_8$. \bigskip

 \noindent{\it a) A sufficient condition for cosemisimplicity}\bigskip
 
\noindent In this subsection $G$ is still general. It turns out useful to tackle  the cosemisimplicity problem at the level of the filtration $\widetilde{\mathcal C}_k$, as this is a better behaved structure, in that it is provided with a  multiplication. Of course, this involves the question of
non-triviality of this filtration, or, more precisely,
whether the image   $\widetilde{M}_\lambda^k$
of $M_\lambda$ in $\widetilde{\mathcal C}_k$  under the quotient map ${\mathcal C}_k\to\widetilde{\mathcal C}_k$ is a matrix coalgebra for some $k\geq n$ and sufficiently many $\lambda\in\Lambda_\ell$.
\medskip

Linear independence can be easily settled.\medskip

\noindent{\bf 10.1. Proposition} {\sl
 \begin{itemize}
      \item[{\rm a)}] The subcoalgebras ${M}_\lambda$,  are linearly independent  in ${\mathcal C}$ as $\lambda$ varies in $\Lambda_\ell$,
      \item[{\rm b)}]       similarly, $\widetilde{M}_\lambda$ are linearly independent in $\widetilde{\mathcal C}_k$ as $\lambda\in\Lambda_\ell^k$, for all $k$.

    \end{itemize}}\medskip
    
    \noindent{\sl Proof} We show b), the proof of a) is easier and can be done along similar lines.
    Let $V_{\lambda, n}$ denote the isotypic submodule of  $p_n(V^{\otimes n})$ of type $V_\lambda$, with orthogonal complement
   $V_{\lambda, n}^\perp$. The subspaces $V_{\lambda, n}^*\otimes V_{\lambda, n}^\perp$, $(V_{\lambda, n}^\perp)^*\otimes V_{\lambda, n}$ and $V_{\lambda, n}^*\otimes V_{\lambda, n}$, for $n=0,\dots, k$, are  linearly independent in 
   ${\mathcal D}_k$. Let $W_\lambda$ denote  their span. Consider the projection $E_\lambda:{\mathcal D}_k\to W_\lambda$ with      complement   $\oplus_{n=0}^k({V_{\lambda, n}}^\perp)^*\otimes V_{\lambda, n}^\perp$. The main point is that $\widetilde{\mathcal I}_k$ is stable under $E_\lambda$. This can be seen noticing  that $\widetilde{\mathcal I}_k$
 is linearly spanned by   $(V_{\lambda, n}^\perp)^*\otimes V_{\lambda, n}$, $V_{\lambda, n}^*\otimes V_{\lambda, n}^\perp$, $\phi^*\otimes A\psi-\phi^*A\otimes\psi$, with $\psi\in V_{\lambda, n}$, and $\phi\in  V_{\lambda, m}$ or $\psi\in {V_{\lambda, n}}^\perp$ and $\phi \in {V_{\lambda, m}}^\perp$; and also by $V_{\lambda, n}^*\otimes Z V^{\otimes n}$, $({V_{\lambda, n}}^\perp)^*\otimes Z V^{\otimes n}$, $(V^{\otimes n})^*Z^*\otimes {V_{\lambda, n}}^\perp$, $(V^{\otimes n})^*Z^*\otimes {V_{\lambda, n}}$.
 Hence if a sum $T=\sum T_\mu$  of elements $T_\mu\in V_\mu^*\otimes V_\mu$ with $\mu\in\Lambda_\ell^k$ lies in 
   $\widetilde{\mathcal I}_k$ then each  $E_\lambda(T)=T_\lambda$ does as well.
 \medskip

Cosemisimplicity in the generic case  was studied  by means of the Haar functional.  We look for a generalisation of that approach to the present setting. Notice that if one can establish that $\widetilde{M}_\lambda^k$ is a matrix coalgebra in $\widetilde{\mathcal C}_k$ then $M_\lambda$ is a matrix coalgebra in ${\mathcal C}(G,\ell)$ as well, by dimension count. We shall verify cosemisimplicity in this stronger form.
\medskip

\noindent{\bf 10.2. Definition} A linear functional $h$ on $\widetilde{\mathcal C}_k$ is said to be a Haar functional if $h(I)=1$ and $h$ annihilates the subcoalgebras $\widetilde{M}^k_\lambda$ for  $\lambda\in\Lambda_\ell^k\backslash\{0\}$.\medskip

Obviously a Haar functional on $\widetilde{\mathcal C}_k$ is unique.  Furthermore, if $\widetilde{\mathcal C}_k$ admits a Haar functional then so does $\widetilde{\mathcal C}_h$ for $h<k$.\medskip

   \noindent{\it Notation}  For  a given $\lambda\in\Lambda_\ell$, let $\text{deg}(\lambda)$ denote the smallest  integer $k$ such that $\lambda\in\Lambda_\ell^k$, or, in other words,  such that $V_\lambda$ is a summand of $V^{\otimes k}$. Furthermore, set $$m({\mathfrak g},\ell):=\text{max}\{\text{deg}(\lambda), \lambda\in\Lambda_\ell\}.$$  \medskip

  \noindent{\bf 10.3. Definition}  We will say that the pair $({\mathfrak g},\ell)$ satisfies the   cosemisimplicity condition if
   \begin{itemize}
      \item[{\rm 1)}]  there is $\tilde{m}\geq m({\mathfrak g},\ell)$ such that $\widetilde{\mathcal C}_{\tilde{m}}$  admits a Haar functional,
      \item [{\rm 2)}] every  $\lambda\in\Lambda_\ell$ has a conjugate 
$\overline{\lambda}\in\Lambda_\ell$  satisfying $\text{deg}(\lambda)+\text{deg}(\overline{\lambda})\leq \tilde{m}$.
\end{itemize}
  
  \medskip
  
  \noindent{\bf 10.4. Lemma} {\sl Let $\lambda\in\Lambda_\ell^k$ have a conjugate   $\overline{\lambda}\in\Lambda_\ell^h$  such that $\widetilde{\mathcal C}_{h+k}$ admits a Haar functional.
Then $\widetilde{M}_\lambda^k$ is a matrix coalgebra in $\widetilde{\mathcal C}_k$.
}\medskip
 
 \noindent{\sl Proof} 
  Let $V_\lambda$ be a summand of 
   $V^{\otimes n}$, $n\leq k$.
  Let $r: {\mathbb C}\to V_{\overline{\lambda}}\otimes V_\lambda$ be defined as in the proof of Theorem 6.7. The composed arrow $r^*r: {\mathbb C}\to V_{\overline{\lambda}}\otimes V_\lambda\to {\mathbb C}$ is nonzero since $\lambda\in\Lambda_\ell$ \cite{Andersen, Wenzl}. In particular,  the trivial submodule defined by $r$ is a summand of $V_{\overline{\lambda}}\otimes V_\lambda$.  But       $(1-p_{h+n})V_{\overline{\lambda}}\otimes V_\lambda$ can not contain a trivial submodule, as otherwise it would be a summand, by multiplicity count. Hence $(1-p_{h+n})\circ r=0$. This shows that $r\in p_{h+n}V^{\otimes h+n}$. If a linear combination $x=\sum_{i,j} \mu_{i,j} e^{\lambda}_{i,j}=0$ vanishes  in $\widetilde{\mathcal C}_k$ then $h(ax)=0$ for all $a\in\widetilde{\mathcal C}_h$, where $h$ is   a Haar functional for $\widetilde{\mathcal C}_{h+k}$.
Now  computations  
 analogous to those of the generic case  show that $\mu_{i,j}=0$.
  \medskip

\noindent{\bf 10.5. Theorem} {\sl    If $({\mathfrak g},\ell)$ satisfies the  cosemisimplicity condition 10.3 then 
 \begin{itemize}
      \item[{\rm a)}]
$\widetilde{M}_\lambda^k$ is a matrix coalgebra in $\widetilde{\mathcal C}_k$ for $k=\text{deg}(\lambda)$, 
      \item[{\rm b)}]
  $M_\lambda$ is a matrix coalgebra in ${\mathcal C}(G,\ell)$, for all $\lambda\in\Lambda_\ell$, hence ${\mathcal C}(G,\ell)$ is cosemisimple:
 $${\mathcal C}(G,\ell)=\bigoplus_{\lambda\in\Lambda_\ell} M_\lambda.$$
 \end{itemize}  }\medskip
 
 \noindent{\sl Proof} The proof follows from Prop. 10.1 and Lemma 10.4. \bigskip

 \noindent{\it b) The case $G={\rm SU}(N)$}\bigskip

\noindent   The rest of the section is dedicated to the proof of the following theorem, which     concludes the main result of the paper.
\bigskip

\noindent{\bf 10.6. Theorem} {\sl  If ${\mathfrak g}={\mathfrak sl}_N$ then $({\mathfrak g}, \ell)$, satisfies the   cosemisimplicity condition 10.3 for all $N\geq2$ and $\ell\geq N+1$ with 
$m({\mathfrak g},\ell)=(N-1)(\ell-N)$ and  $\tilde{m}:=m({\mathfrak g},\ell)+\ell-1$.}\medskip

   We start fixing   notation of type $A_{N-1}$ root systems \cite{Humphreys}. Consider ${\mathbb R}^{N}$ with the usual euclidean inner product, and let $e_1,\dots, e_{N}$ be the canonical orthonormal basis.
Consider the subspace $E\subset {\mathbb R}^{N}$ of elements $\mu_1e_1+\dots+\mu_{N} e_{N}$ such that
$\mu_1+\dots+\mu_{N}=0$.
 The $A_{N-1}$ root system is 
 $\Phi=\{e_i-e_j, i\neq j\}$,   the simple roots are $ \alpha_i=e_{i}-e_{i+1},$ and the fundamental weights are $\omega_i=e_1+\dots+e_i-\frac{i}{N}e$, where
  $e:=e_1+\dots+e_{N}$ and $ i=1,\dots, N-1$.
The weight lattice and  the dominant Weyl chamber  of $(E, \Phi)$ are respectively  $$\Lambda=\{\lambda=\lambda_1e_1+\dots+\lambda_{N-1}e_{N-1}-\frac{\lambda_1+\dots+\lambda_{N-1}}{N}e, \quad \lambda_i\in{\mathbb Z}\},$$
 $$\Lambda^+=\{\lambda\in\Lambda: \lambda_1\geq\lambda_2\geq\dots\geq\lambda_{N-1}\geq0\}.$$
 The highest root is $\theta=e_1-e_{N}$. Direct computations show that  
 $$\Lambda_\ell=\{\lambda\in\Lambda^+: \lambda_1\leq \ell-N\}, \quad \overline{\Lambda_\ell}=\{\lambda\in\Lambda^+: \lambda_1\leq \ell-N+1\}.$$
 We have $V=V_{\omega_1}$, the vector representation, its weights are   
$$\gamma_1=\omega_1,\quad \gamma_i=e_i-\frac{1}{N}e,\quad i=2, \dots, N.$$
In the next lemmas we shall make   use of the  decomposition into irreducibles 
$$V_\lambda\otimes V\simeq \oplus V_{\lambda+\gamma_i},\quad  \lambda\in\Lambda_\ell,$$
 in the category ${\mathcal T}_\ell$ of tilting modules, where the sum is extended to all $i$ such that $\lambda+\gamma_i$ is dominant, (Theorem 4.4, c).)
We derive two simple consequences.\medskip

 \noindent{\bf 10.7. Lemma} {\sl For any $\lambda\in\Lambda_\ell$,
the negligible submodule $N_\lambda$  of $V_\lambda\otimes V$ is non-zero  if and only if
$\lambda_1=\ell-N$, and one has $N_\lambda\simeq V_{\lambda+\omega_1}$.}\medskip

\noindent{\sl Proof} The indicated submodule is isomorphic to the  sum corresponding to $\lambda+\gamma_i\in\overline{\Lambda_\ell}\setminus\Lambda_\ell$, which is realised if and only if $i=1$ and $\lambda_1=\ell-N$.\medskip

\noindent{\bf 10.8. Lemma} {\sl  $m({\mathfrak sl}_N,\ell)=(N-1)(\ell-N)$ .}\medskip

\noindent{\sl Proof}
Let $\lambda\in\Lambda_\ell$ be determined by non negative integers $\lambda_1,\dots$, $\lambda_{N-1}$ as above,
and let us identify $\lambda$ with $(\lambda_1,\dots,\lambda_{N-1})$. 
The dominant weight with  coordinates all equal to $\ell-N$ lies in $\Lambda_\ell$,   it is a summand of $V^{\otimes (N-1)(\ell-N)}$ and    this is the smallest possible power.
We need to show that every module $V_\lambda$ with $\lambda\in\Lambda_\ell$ is a summand of some
$p_nV^{n}$ with $n\leq (N-1)(\ell-N)$. Notice that
 $(1,0,\dots,0),\dots$, $(\lambda_1, 0,\dots, 0)$,
 $(\lambda_1, 1, 0, \dots,0),\dots$,  $(\lambda_1, \lambda_2, \dots, 0),$  $\dots,$ $(\lambda_1,\dots,\lambda_{N-1})$ is a sequence of $\lambda_1+\lambda_2+\dots+\lambda_{N-1}$ dominant weights of $\Lambda_\ell$ starting with $\omega_1$ and obtained from one another by adding a weight of $V$. The fusion rules  then show that $V_\lambda$ is  a summand of $p_nV^n$, where
 $n=\lambda_1+\lambda_2+\dots+\lambda_{N-1}\leq {(N-1)}(\ell-N)$.\medskip

We next derive   information on the   negligible   summands of $V^{\otimes n}$, including the non-canonical ones, for the bounded values of $n$. The following Lemma plays a crucial role for the Haar functional. \medskip   
 
 \noindent{\bf  10.9. Lemma} {\sl 
 No negligible summand of $V^{\otimes n}$, with $n$ up to $\tilde{m}=m({\mathfrak sl}_N,\ell)+\ell-1$
 contains the trivial module among the  successive factors of its Weyl filtrations.}\medskip

\noindent{\sl Proof}  
A negligible summand of $V^{\otimes n}$ is isomorphic to  a   summand of $N_n=(1-p_n)V^{\otimes n}$. Furthermore
  the inductive procedure described in Sect. 5 shows that $N_n$    is in turn spanned
 by the summands 
 $$N(p_rV^{\otimes r}\otimes V)\otimes V^{\otimes (n-r-1)}, \quad r=1,\dots, n-1,\eqno(10.1)$$   
 where
 $N(p_rV^{\otimes r}\otimes V)$ is the     canonical negligible summand of $p_rV^{\otimes r}\otimes V$, hence we are reduced to show 
 the statement for these modules. 
 
 Now on one hand
 $N(p_rV^{\otimes r}\otimes V)$ is completely reducible and the dominant weights of the irreducible components are of the form
 $\lambda+\omega_1=(\ell-N+1,\lambda_2,\dots,\lambda_{N-1})$ by Lemma 10.7.

 On the other, the dominant weights appearing in the Weyl filtrations of $(10.1)$ are the same as  those appearing in the irreducible decomposition
 of the corresponding module at the level of the semisimple category $\text{Rep}({\mathfrak g})$, see Prop. 3 and Remark 2 in \cite{Sawin}.
 
Hence we are reduced to show that  the smallest integer $t$ such that
$V_{\lambda+\omega_1}\otimes V^{\otimes t}$ contains the trivial module in 
$\text{Rep}({\mathfrak g})$ satisfies $$t+r+1> (N-1)(\ell-N)+\ell-1.$$
We compute $t$. For a general dominant weight $\mu=(\mu_1,\dots,\mu_{N-1}),$ the shortest path
to the trivial module is obtained as follows. If $\mu_{N-1}>0$ we consider 
the path $$\mu+\gamma_N, \quad \mu+2\gamma_N,\quad \dots, \quad\mu+\mu_{N-1}\gamma_N$$ which  lowers $\mu$ to $$\mu'=(\mu_1-\mu_{N-1},\dots, \mu_{N-2}-\mu_{N-1}, 0),$$
and we have thus used $\mu_{N-1}$ powers of $V$. We need no such power if $\mu_{N-1}=0$.
We   proceed in the same way for the $N-2$ coordinate and the new weight $\mu'$, but we   now need to follow a longer    path, due to vanishing of the last coordinate, and the shortest is
 $$\mu'+\gamma_{N-1},\quad \mu'+\gamma_{N-1}+\gamma_N, \quad  \mu'+2\gamma_{N-1}+\gamma_N,\quad \mu'+2\gamma_{N-1}+2\gamma_N,\dots$$
 using $2(\mu_{N-2}-\mu_{N-1})$ more powers of $V$. Continuing in this way, we find
 $$t=\mu_{N-1}+2(\mu_{N-2}-\mu_{N-1})+3(\mu_{N-3}-\mu_{N-2})+\dots+(N-1)(\mu_1-\mu_2).
 $$
 Taking into account  the fact that in $\text{Rep}({\mathfrak sl}_N)$, the dominant weights $\mu$ appearing in $V^{\otimes r+1}$ satisfy $\mu_1+\dots+\mu_{N-1}\leq r+1,$ we easily get, for   $\mu=\lambda+\omega_1$,
 $$t+r+1\geq t+\mu_1+\dots+\mu_{N-1}=N\mu_1=(N-1)(\ell-N)+\ell>(N-1)(\ell-N)+\ell-1.$$
 which finally gives  the desired estimate.
  \medskip

\noindent{\bf 10.10. Remark} The  proof   also shows that $N(p_rV^{\otimes r}\otimes V)=0$ for  $r<\ell-N$.
 \medskip

  \noindent{\bf 10.11. Corollary} {\sl Let $n\leq \tilde{m}$.
 \begin{itemize}
  \item[{\rm a)}] 
  $e_n$ is a central element of $(V^{\otimes n}, V^{\otimes n})$,
   \item[{\rm b)}] $e_n\circ1_{V^q}\otimes (1-p_j)\otimes 1_{V^u}=0$, $q+j+u=n$.
   \end{itemize}
 
 }\medskip

\noindent{\sl Proof} a) The multiplicity of the trivial representation in $p_nV^n$ is the same as that of the classical case, by the previous lemma,  hence $e_n$ is the specialisation of   a central intertwiner of the generic case. b) 
$e_n\circ1_{V^q}\otimes (1-p_j)\otimes 1_{V^u}\circ p_n=0$  by Lemma 5.3, hence
$$e_n\circ1_{V^q}\otimes (1-p_j)\otimes 1_{V^u}=e_n\circ1_{V^q}\otimes (1-p_j)\otimes 1_{V^u}\circ (1-p_n)=$$
$$1_{V^q}\otimes (1-p_j)\otimes 1_{V^u}\circ e_n\circ (1-p_n)=0$$
by a).
\medskip

 We next consider the $\tilde{m}$-th term, $\widetilde{\mathcal C}_{\tilde{m}}$, of the associative filtration $\widetilde{\mathcal C}_k$ corresponding to ${\mathcal C}({\rm SU}(N), \ell)$.    For convenience we recall 
  that $\widetilde{\mathcal C}_{\tilde{m}}={\mathcal D}_{\tilde{m}}/{\widetilde{\mathcal J}_{\tilde{m}}}$, where ${\widetilde{\mathcal J}_{\tilde{m}}}$ is spanned by elements of the form $\phi\otimes A\circ\psi-\phi\circ A\otimes\psi$, together with $\phi\otimes Z\circ\psi'$, and $\phi'\circ 
Z'\otimes\psi$,
 where 
 $$A\in(V^{\underline{\otimes}m}, V^{\underline{\otimes}n}),\quad\quad  Z, \,Z'^*=p_{q+j+u}\circ1_{V^q}\otimes (1-p_j)\otimes 1_{V^u}$$ with $m$, $n$, $q+j+u\leq \tilde{m}$.
  \medskip

We   define the linear functional
 $$h: {\mathcal D}_{\tilde{m}}\to {\mathbb C},$$
 setting 
 $$h(\phi\otimes\psi)=\phi(e_n\psi),\quad \phi\otimes\psi\in (V^n)^*p_n\otimes p_n V^n,\quad n\leq \tilde{m}. $$
 where $e_n\in(p_nV^n, p_nV^n)$ is the orthogonal projection onto the isotypical component of the trivial representation.\medskip
 
 \noindent{\bf  10.12. Theorem} {\sl The functional $h$ annihilates ${\widetilde{\mathcal J}_{\tilde{m}}}$. Hence  it gives rise
 to a Haar functional on $\widetilde{\mathcal C}_{\tilde{m}}$.}\medskip
 
 \noindent{\sl Proof}
The functional $h$ obviously annihilates elements $\phi\otimes A\psi-\phi A\otimes\psi$. Furthermore it also annihilates elements of the form $\phi\otimes Z\psi'$, $\phi' Z'\otimes\psi\in {\widetilde{\mathcal J}_{\tilde{m}}}$, by Corollary 10.11, b). The rest is now clear.\medskip 

We finally verify the needed upper bound for  $\text{deg}(\lambda)+\text{deg}(\overline{\lambda})$ for all $\lambda\in\Lambda_\ell$. We are interested in property c) of the following proposition.
\medskip

\noindent{\bf 10.13. Proposition} {\sl If $\lambda=(\lambda_1,\dots,\lambda_{N-1})\in\Lambda^+$ then
\begin{itemize}
\item[{\rm a)}] $\text{deg}(\lambda)=\lambda_1+\dots+\lambda_{N-1}$,
\item[{\rm b)}]  $\overline{\lambda}=(\lambda_1, \lambda_1-\lambda_{N-1}, \lambda_1-\lambda_{N-2},\dots,\lambda_1-\lambda_2)$,
 \item[{\rm c)}] $\text{deg}(\lambda)+\text{deg}(\overline{\lambda})=N\lambda_1\leq \tilde{m}$ for $\lambda\in\Lambda_\ell$.
\end{itemize}}
 
 \noindent{\sl Proof} Properties a) and b) are classical, and c) follows from an easy computation, since $\lambda_1\leq(\ell-N)$ for $\lambda\in\Lambda_\ell$ and $\tilde{m}=(N-1)(\ell-N)+\ell-1$. We sketch a proof for completeness. a) follows from the fusion rules for the powers of $V$, recalled at the beginning of this section. b)  Let $w_0$ be longest element of the Weyl group. For ${\mathfrak sl}_N$, this is the permutation group ${\mathbb P}_{N}$ and $w_0$ is the permutation reversing the order of $(e_1,\dots, e_N)$. Then 
 $$\overline{\lambda}=-w_0\lambda=$$
 $$-(\lambda_1e_N+\lambda_2 e_{N-1}+\dots+\lambda_{N-1}e_2)+
 \frac{\lambda_1+\dots+\lambda_{N-1}}{N}e=$$
 $$(\lambda_1, \lambda_1-\lambda_{N-1}, \lambda_1-\lambda_{N-2},\dots,\lambda_1-\lambda_2).$$\medskip
 
  \appendix{}
 \section{\\ On the fusion  $C^*$-category of type $A$}
   
In this appendix we  illustrate  a concrete realisation of the fusion $C^*$--category ${\mathcal F}_\ell$ associated to ${\mathfrak sl}_N$ in the sense of Theorem 5.4 and we identify a single  generator, the quantum determinant element $S$.  Most of the results are   due to \cite{KW}, but we have not been able to find a reference where   $S$ and the conjugation structure, Prop. A.5, are made explicit. The main aspect is the use of the Kirillov-Wenzl inner product.

 The fusion rules seen in the previous section  show that   the trivial object of ${\mathcal F}_\ell$  is contained with multiplicity $1$ in  
 $V^{\underline{\otimes} N}$, hence there is  a nonzero intertwiner, the quantum  determinant element,
 $$S\in(\iota, V^{\underline{\otimes} N})$$
  unique up to a scalar multiple.  
  \medskip

  \noindent{\bf A.1. Theorem} {\sl The quotient category
${\mathcal F}_\ell$ of type $A_{N-1}$ is generated by the quantum determinant $S$ as a tensor $^*$--category.}\medskip

  \noindent{\sl Proof}  We have  that
  $(V^{\underline{\otimes} s}, V^{\underline{\otimes} t})\neq 0\quad\Rightarrow s=t+hN$
for some $h\in{\mathbb Z}$, by the fusion rules. If  for example $h>0$   we can   write an element  $B\in (V^{\underline{\otimes} s}, V^{\underline{\otimes} t})$  in the form
$$ B=S^*\otimes 1_t\circ S^*\otimes 1_{t+N}\circ\dots \circ S^*\otimes 1_{t+hN}\circ Y , \quad \text{with} \quad Y\in (V^{\underline{\otimes} s}, V^{\underline{\otimes} s}).$$
 Indeed, we may choose 
 $Y=S\otimes 1_{t+hN}\circ\dots\circ S\otimes 1_t\circ B $
 which  lies in the claimed space as
 $$Y=S\otimes 1_{t+hN}\circ\dots\circ S\otimes 1_t\circ p_t\circ B  =p_s\circ S\otimes 1_{t+hN}\circ\dots\circ S\otimes 1_t\circ p_t\circ B.$$
  by Lemma 5.1. On the other hand, by Schur-Weyl duality at roots of unity \cite{Du} the intertwining space $(V^{\underline{\otimes} s}, V^{\underline{\otimes} s})$  is
linearly spanned by the representation $\underline{\varepsilon}(b)=p_n\varepsilon(b)p_n$ of the braid group ${\mathbb B}_s$ obtained reducing the  translates $1_w\otimes\varepsilon\otimes 1_{w'}$ of the basic braiding operator $\varepsilon\in(V^{\otimes 2}, V^{\otimes 2})$ by $p_n$.
Since the $p_n$ themselves lie in the image of ${\mathbb C}{\mathbb B}_s$,
we are reduced  to show that the basic braiding operator derived from the $R$-matrix acting on $(V^{\otimes 2}, V^{\otimes 2})$ lies in the tensor $^*$--subcategory generated by $S$, and this will follow from  Prop. A.5. \medskip

If for example   $\ell\geq 2N$ then tensor powers $V^{\otimes k}$ with $k$ up to $N$ are completely reducible, with multiplicities given as in the classical case. In particular,  the $k$-th fundamental representation $V_{\omega_k}$, $k=1,\dots, N-1$, is   contained in $V_{\omega_{k-1}}\otimes V$, and hence in $V^{\otimes k}$,  with multiplicity $1$.

Set 
$$\varepsilon=\Sigma R\in(V^{\otimes 2}, V^{\otimes 2}),\quad\quad g=q^{1/N}\varepsilon.$$
We recall that $g$ gives rise to a representation of the Hecke algebra $H_n(q)$ of type $A$ on $V^{\otimes n}$. Indeed $g$
has the following spectral decomposition
on $V^{\otimes 2}$,
$$g=q[2\omega_1]-q^{-1}[\omega_2],$$
hence  the translates  $g_i=1_{i-1}\otimes g\otimes 1_{n-i-1}$  satisfy the defining relations $(g_i-q)(g_i+\frac{1}{q})=0$ of  $H_n(q)$. In particular, $\varepsilon$ is determined by the single projection
$[\omega_2]$.
We thus need to show that $[\omega_2]$   is in turn   an arrow of the tensor $^*$--category generated by $S$.  
Showing this claim, as we next see, involves determining the conjugate equations for the submodule $V_{\omega_2}$ of $V^{\otimes 2}$. Since it will  not be more effort,   we shall  identify the conjugate equations    for all the fundamental representations $V_{\omega_n}$.

Consider the following sequence of central elements of $\alpha_{-n}\in H_n(q)$ inductively defined by $\alpha_{-1}=I$ and
$$\alpha_{-n-1}:=\sigma(\alpha_{-n})-q^{-1}g_1\sigma(\alpha_{-n})+\dots+(-1)^n q^{-n}g_n\dots g_1\sigma(\alpha_{-n}),$$
with $\sigma: H_n(q)\to H_{n+1}(q)$   the homomorphism taking $g_i$ to $g_{i+1}$.  
In particular,   $\alpha_{-2}$ is a scalar multiple of $[\omega_2]$.  
The next two lemmas recall certain known facts \cite{KW, Wenzl}, see also \cite{P}. 
 
\medskip

\noindent{\bf A.2.  Lemma} {\sl

\begin{itemize}
   \item[{\rm a)}]  $g_i \alpha_{-n}=\alpha_{-n}g_i=-q^{-1}\alpha_n$, $i=1,\dots, n-1$,
    \item[{\rm b)}]  $\alpha_{-n}^2=\lambda_n \alpha_{-n}$,  with $\lambda_n:=q^{-n(n-1)/2}[n]_q!$.
     \end{itemize}
    }\medskip

Hence if   $\ell>n$, we may consider the idempotent $E_{-n}:=\lambda_n^{-1}\alpha_{-n}$.
 We next identify its  range.
The Hecke algebra generator $g=q^{1/N}\varepsilon$
  acts on a suitable orthonormal basis $\psi_i$, $i=1,\dots, N$ of the Hilbert space of $V=V_{\omega_1}$ as  
$$g\psi_i\otimes\psi_i=q\psi_i\otimes\psi_i,$$
$$g \psi_i\otimes\psi_j=\psi_j\otimes\psi_i,\quad i>j$$
$$g\psi_i\otimes\psi_j=\psi_j\otimes\psi_i+(q-\frac{1}{q})\psi_i\otimes\psi_j,\quad i<j.$$
The quantum determinant $S\in(\iota, V^{\otimes N})$ is given by
$$ S=\sum_{p\in{\mathbb P}_N }(-q)^{i(p)}\psi_{p(1)}\otimes\dots\otimes\psi_{p(N)}.$$
More generally, we also introduce the following antisymmetrized elements   in $V^{\otimes n}$ for $n\leq N$ spanning the submodule $V_{\omega_n}$ of $V^{\otimes n}$.  Consider indices $i_1<i_2<\dots<i_n$ with $i_j\in\{1,\dots, N\}$ and set, for $\underline{i}=(i_1,\dots, i_n)$,
$$ S_{\underline{i}}=\sum_{p\in{\mathbb P}_n }(-q)^{i(p)}\psi_{i_{p(1)}}\otimes\dots\otimes\psi_{i_{p(n)}}.$$
It will be useful to introduce also  corresponding  antisymmetrized elements over $q^{-1}$,
$$ \tilde{S}_{\underline{i}}=\sum_{p\in{\mathbb P}_n }(-q^{-1})^{i(p)}\psi_{i_{p(1)}}\otimes\dots\otimes\psi_{i_{p(n)}}.$$
However, some caution is needed, as  the $\tilde{S}_{\underline{i}}$ do {\it not} span a submodule of $V^{\otimes n}$.
\medskip

\noindent{\bf A.3. Lemma} {\sl   We have that
\begin{itemize}
\item[{\rm a)}]  $\alpha_{-n}=0$ for $n>N$, 
   \item[{\rm b)}]
   $E_{-n}$ is the idempotent of $(V^{\otimes n}, V^{\otimes n})$ corresponding to   $V_{\omega_n}$   for $n\leq N-1$
and
  to the trivial submodule for $n=N$.
  \end{itemize}

       }\medskip

\noindent{\sl Proof} This follows from    
$\alpha_{-n}\psi_{i_1}\otimes\dots\otimes\psi_{i_n}=0$ if two indices repeat, while  
$q^{n(n-1)}\alpha_{-n}\psi_{i_{p^{-1}(1)}}\otimes\dots\otimes\psi_{i_{p^{-1}(n)}}=(-q)^{i(p)}S_{\underline{i}}$
for any $\underline{i}=i_1<\dots<i_n$ and $p\in{\mathbb P}_n$.
 \medskip

A main aspect for the values $q\in{\mathbb T}$
  is that the $^*$--structure   relies on the Kirillov-Wenzl inner product on tensor product spaces. 
We next identify an orthonormal basis of $V_{\omega_n}$ in the correct Hilbert space. \medskip

\noindent{\bf A.4. Proposition} {\sl  Consider indices $\underline{i}=i_1<\dots< i_n$, $\underline{j}=j_1<\dots<j_n$. Then

\begin{itemize}
\item[{\rm a)}]  
     $E_{-n} $ is a selfadjoint projection of $(V^{\otimes n}, V^{\otimes n})$ for $n\leq N$,
   \item[{\rm b)}]
$\overline{R}^{(n)}S_{\underline{i}}=q^{n(n-1)/2}\tilde{S}_{\underline{i}}.$
   \item[{\rm c)}] $(S_{\underline{i}}, S_{\underline{j}})=0$ for $\underline{i}\neq\underline{j}$,
   \item[{\rm d)}] $\|S_{\underline{i}}\|^2=[n]_q!$.
\end{itemize} }\medskip

\noindent{\sl Proof} 
a) Selfadjointness of $E_{-n}$ was shown in \cite{P} for 
 the involution introduced in  \cite{W_Hecke}, which, however, does coincide with that induced by Kirillov-Wenzl inner product,
as they both make  the spectral idempotent $E_{-2}=[\omega_2]$ selfadjoint.
b) Recall from Sect. 2 that
$$\overline{R}^{(n)}=R^{(n)}\Theta^{(n)}$$ and that 
$$R^{(n)}=R^{(n-1)}\otimes 1\circ \Delta^{(n-2)}\otimes 1(R)=R^{(n-1)}\otimes 1\circ (R_{1\, n}\dots R_{n-1\, n}).$$
Hence 
$$R^{(n)}= R_{12}\circ\dots\circ(R_{1\, n-1}\dots R_{n-2\, n-1})\circ (R_{1\, n}\dots R_{n-1\, n})$$
It turns out convenient to decompose the permutation $\Sigma_n$ as a product of transpositions
$$\Sigma_n=\Sigma_{n-1\,n}\circ\dots\circ(\Sigma_{23}\dots\Sigma_{n-1\,n})\circ(\Sigma_{12} \Sigma_{23}\dots\Sigma_{n-1\,n}).$$
Distributing these transpositions among the $R_{i\, j}$ gives
$$R^{(n)}=\Sigma_n\circ\varepsilon_{n-1}\circ \dots(\varepsilon_2\dots\varepsilon_{n-1})\circ(\varepsilon_1\dots\varepsilon_{n-1})$$
An  easy computation shows that 
$$g_i S_{\underline{i}}=-q^{-1}S_{\underline{i}},\quad i=1,\dots n-1$$
and this implies  
$$q^{n(n-1)/2N}R^{(n)}S_{\underline{i}}=
\Sigma_n\circ g_{n-1}\circ \dots(g_2\dots g_{n-1})\circ(g_1\dots g_{n-1})S_{\underline{i}}=$$
$$(-q^{-1})^{n(n-1)/2}\Sigma_nS_{\underline{i}}=$$
$$ (-q^{-1})^{n(n-1)/2}\sum_{p\in{\mathbb P}_n}(-q)^{i(p)}\psi_{i_{p(n)}}\otimes\dots\otimes\psi_{i_{p(1)}}=$$
$$\sum_{p\in{\mathbb P}_n}(-q^{-1})^{i(p)}\psi_{i_{p(1)}}\otimes\dots\otimes\psi_{i_{p(n)}}=\tilde{S}_{\underline{i}}.$$
On the other hand, recall from Sect. 3 that  $\Theta^{(n)}$ acts as scalar multiplication by
$$q^{n(\omega_1+2\rho,\omega_1)-(\omega_n+2\rho,\omega_n)/2}$$ on  the   submodule $V_{\omega_n}$ of $V^{\otimes n}$.
A straightforward computation gives
$(\omega_n,\omega_n+2\rho)= n+nN-n^2-\frac{n^2}{N}$
and  the conclusion follows.
c) follows from b) and d)  from explicit computation of the Kirillov-Wenzl norm,
$$\|S_{\underline{i}}\|^2=(S_{\underline{i}}, \overline{R}^{(n)}S_{\underline{i}})_p=$$
$$q^{n(n-1)/2}\sum_{p,p'\in{\mathbb P}_n}((-q)^{i(p)}\psi_{i_{p(1)}}\otimes\dots\otimes
\psi_{i_{p(n)}}, (-q^{-1})^{i(p')}\psi_{i_{p'(1)}}\otimes\dots\otimes\psi_{i_{p'(n)}})_p=$$ 
$$q^{n(n-1)/2} \sum_{p\in{\mathbb P}_n} (q^{-2})^{i(p)}=[n]_q!$$
where  inner products are assumed antilinear in the first variable.\medskip

Consider a partition of $N=m+n$, and fix  indices $\underline{i}=i_1<\dots< i_n$ in $\{1,\dots, N\}$.
We denote by $\underline{j}=j_1<\dots<j_m$ the indices obtained from $\{1,\dots,N\}$ after removing $i_1,\dots, i_n$. We shall refer to $\underline{j}$ as   conjugate to $\underline{i}$.
Obviously, $\underline{i}$ is conjugate to $\underline{j}$ as well.   
Fix  permutations  $s\in{\mathbb P}_m$ and $r\in{\mathbb P}_n$ and  define $p\in{\mathbb P}_N$ by
$$p(1)=j_{s(1)},\dots, p(m)=j_{s(m)},$$
$$p(m+1)=i_{r(1)},\dots, p(m+n)=i_{r(n)}.$$  
 Every $p\in{\mathbb P}_N$ can   be uniquely written in this form. Denoting by $N({\underline{i}})$ the number of pairs $(j_h, i_k)$ with $j_h> i_k$, we have
$i(p)=i(r)+i(s)+N(\underline{i})$.  We   thus have
$$S=
\sum_{\underline{i}}(-q)^{N(\underline{i})} S^{(m)}_{\underline{j}}\otimes S^{(n)}_{\underline{i}},$$
where we have introduced upper indices clarifying  the corresponding  tensor power of $V$. 
Exchanging the roles of $m$ and $n$, we can also write
$$S=
\sum_{\underline{h}}(-q)^{N(\underline{h})} S^{(n)}_{\underline{k}}\otimes S^{(m)}_{\underline{h}},
$$
where $\underline{h}$ and $\underline{k}$ are conjugate indices as well. Similar relations hold for $\tilde{S}$, where $q$ is replaced by $q^{-1}$. These decompositions turn out useful  for  the next result.
\medskip

\noindent{\bf A.5. Proposition} {\sl  The following conjugate equations hold in the  tilting tensor $^*$--category ${\mathcal T}_\ell$ associated to ${\mathfrak sl}_N$,
$$S^*\otimes 1_{V^{\otimes m}}\circ 1_{V^{\otimes m}}\otimes S=(-1)^{mn}[m]_q![n]_q!\,E_{-m}$$
$$S^*\otimes 1_{V^{\otimes n}}\circ 1_{V^{\otimes n}}\otimes S=(-1)^{mn}[m]_q![n]_q!\,E_{-n},$$
where $m+n=N$. Hence $V_{\omega_m}$ and $V_{\omega_n}$ are conjugate  of each other.  }\medskip

\noindent{\sl Proof} It suffices to show the first equation. Evaluating the left hand side on a vector $\psi\in V^{\otimes m}$ gives 
$$ \sum_{\underline{h}}(-q)^{N(\underline{h})}(S, \psi\otimes S_{\underline{k}}^{(n)})S_{\underline{h}}^{(m)} = $$
$$\sum_{\underline{h}}(-q)^{N(\underline{h})}(\overline{R}^{(N)}S, \psi\otimes S_{\underline{k}}^{(n)})_p\,S_{\underline{h}}^{(m)} =$$
$$q^{-N(N-1)/2}\,\sum_{\underline{h}}(-q)^{N(\underline{h})}(\tilde{S}, \psi\otimes S_{\underline{k}}^{(n)})_p\,S_{\underline{h}}^{(m)} =$$
$$q^{-N(N-1)/2}\,\sum_{\underline{h},\underline{i}}(-q)^{N(\underline{h})}(-q)^{N({\underline{i}})}(\tilde{S}_{\underline{j}}^{(m)}\otimes\tilde{S}_{\underline{i}}^{(n)}, \psi\otimes S_{\underline{k}}^{(n)})_p\,S_{\underline{h}}^{(m)}=$$
$$q^{-N(N-1)/2}\,\sum_{\underline{h}} (-q)^{N(\underline{h})}(-q)^{N({\underline{k}})} (\tilde{S}_{\underline{h}}^{(m)},\psi)_p(\tilde{S}_{\underline{k}}^{(n)}, S_{\underline{k}}^{(n)})_p\,S_{\underline{h}}^{(m)}=$$
$$ (-1)^{mn}[n]_q!\,\sum_{\underline{h}} ({S}_{\underline{h}}^{(m)},\psi)  S_{\underline{h}}^{(m)} =$$
$$(-1)^{mn}[n]_q![m]_q!\, E_{-n}\psi. $$
where we have  successively used selfadjointness of $\overline{R}^{(N)}$ with respect to the tensor product form, Prop. 3.6, $N(\underline{h})+N({\underline{k}})=mn$, and the fact that 
$\{([m]_q!)^{-1/2}\,S^{(m)}_{\underline{h}}\}$ is an orthonormal basis of $V_{\omega_m}$ with respect to the Kirillov-Wenzl inner product.\medskip

\noindent{\it Acknowledgements} C.P. would like to thank Martijn Caspers, Huichi Huang, and Thomas Timmermann for extending to her the invitation to the conference `Quantum Groups and Operator Algebras' held in M\"unster in May 2014, which constituted a strong motivation  for her to work.
She is also grateful to  Gabriella B\" ohm for fruitful conversations   on  that occasion. We   thank   Fabio Cipriani and St\'ephane Vassout for inviting us to the 
 conference `Advances in Noncommutative Geometry',  Paris 2015. 
We acknowledge support from GREFI-GENCO and Sapienza Universit\`a di Roma, AST fundings.

\end{document}